\documentclass [leqno,12pt]{article}
\usepackage[utf8]{inputenc} 
\usepackage[T1]{fontenc} 
\usepackage{amsmath}
\usepackage{amsfonts}
\usepackage{amsthm}
\usepackage{color}
\usepackage{array,multirow,makecell}
\usepackage[table]{xcolor}

\usepackage{graphics,graphicx,psfrag, bbm}

\newlength{\oldparindent}
\newcommand{\cL}{{\mathbb {L}}}

\newcommand{\bpf}{\begin{preuve}}
\newcommand{\epf}{ \end{preuve} \medskip}

\newcommand{\benum}{\begin{enumerate}}
\newcommand{\eenum}{\end{enumerate}}

\newcommand{\bitem}{\begin{itemize}}
\newcommand{\eitem}{\end{itemize}}

\newcommand{\brmq}{\begin{rmq}}
\newcommand{\ermq}{\end{rmq}}

\newcommand{\brmqs}{\begin{rmqs}}
\newcommand{\ermqs}{\end{rmqs}}

\newcommand{\bapp}{\begin{application}}
\newcommand{\eapp}{\end{application}}

\newcommand{\bapps}{\begin{applications}}
\newcommand{\eapps}{\end{applications}}

\newcommand{\bdefi}{\begin{definition}}
\newcommand{\edefi}{\end{definition}}

\newcommand{\beq}{\begin{equation}}
\newcommand{\eeq}{\end{equation}}

\def\bpm{\begin{pmatrix}}
\def\epm{\end{pmatrix}}

\newcommand{\bcas}{\begin{cases}}
\newcommand{\ecas}{\end{cases}}

\newcommand{\bex}{\begin{exemp}}
\newcommand{\eex}{\end{exemp}}

\newcommand{\bexs}{\begin{exemps}}
\newcommand{\eexs}{\end{exemps}}

\newcommand{\bprop}{\begin{proposition}}
\newcommand{\eprop}{\end{proposition}}

\newcommand{\bthm}{\begin{theoreme}}
\newcommand{\ethm}{\end{theoreme}}

\newcommand{\bcor}{\begin{corollaire}}
\newcommand{\ecor}{\end{corollaire}}

\newcommand{\blem}{\begin{lemme}}
\newcommand{\elem}{\end{lemme}}

\newcommand{\beqna}{\begin{eqnarray}}
\newcommand{\eeqna}{\end{eqnarray}}

\newcommand{\beqnas}{\begin{eqnarray*}}
\newcommand{\eeqnas}{\end{eqnarray*}}

\newcommand{\cA}{{\mathcal A}}

\definecolor{green}{rgb}{0,.7,.2}
\definecolor{orange}{rgb}{0.9,.5,0}

\newcommand{\LL}{{\rm L}}

\def\det{{ \rm{det}}}  

\def\Id{{\rm{Id}}} 

\def\cA{{\mathcal A }}

\def\cC{{\mathcal C}}
\def\cD{{\mathcal D}}

\def\cH{{\mathcal  H}}
\def\cI{{\mathcal  I}}

\def\cL{{\mathcal L }}

\def\cO{{\mathcal  O}}

\def\cT{{\mathcal  T}}

\def\cV{{\mathcal V}}

\def\bbC{{\mathbb{C}}}

\newcommand{\bbN}{{\mathbb {N}}}
\newcommand{\bbR}{{\mathbb {R}}}
\newcommand{\bbS}{{\mathbb {S}}} 

\newcommand{\bbZ}{\mathbb {Z}}

\newtheorem{theoreme}{Theorem}[section]
\newtheorem{lemme}[theoreme]{Lemma}
\newtheorem{definition}[theoreme]{Definition}
\newtheorem{proposition}[theoreme]{Proposition}
\newtheorem{corollaire}[theoreme]{Corollary}

\newenvironment{exemp}{\noindent{\bf Example. --- }}{\par}
\newenvironment{exemps}{\noindent{\bf Examples}\benum}{\eenum\par}
\newtheorem{rmq}[theoreme]{Remark}
\newtheorem{rmqs}[theoreme]{Remarks}

\newenvironment{preuve}{\noindent{\it Proof. --- }}
{\hfill\rule{1.3mm}{2mm}\par} 
\newenvironment{application}{\noindent{\bf Application. --- }}{\par}
\newenvironment{applications}{\noindent{\bf Applications. --- 
}\benum}{\eenum\par}

\theoremstyle{definition}

\author{Dominique Bakry, Xavier Bressaud \footnote{Institut de Math\'ematiques de Toulouse, Universit\'e Paul Sabatier, 118 
route de Narbonne, 31062 Toulouse, France}}
\date{\today}

\date{\today}

\makeatletter \renewcommand{\@oddfoot}{\sl \small
 \hfil \thepage\hfil \today}
\renewcommand{\@oddhead}{\sl \small
 \hfil preprint under construction}
 \makeatother

\begin{document}

\title{Diffusions with polynomial eigenvectors via finite subgroups of $O(3)$}

\maketitle

 \abstract{We provide new examples of diffusion operators in dimension 2 and 3 which have orthogonal polynomials as eigenvectors. Their construction rely on the finite subgroups of $O(3)$ and their invariant polynomials.}
\medskip


\medskip

\section{Introduction}

We investigate in this paper new examples of bounded domains in $\bbR^2$ on which there exists a probability measure $\mu$ with an orthonormal  basis  of $\cL^2(\mu)$  such that the elements of this basis are eigenvectors of a diffusion operator. To determine such a basis, one needs first to define a valuation (a definition for the degree)  for a polynomial in two variables. The complete determination of all the possible such domains in $\bbR^2$ has been carried in~\cite{BOZ2013}, under the restriction that the valuation is the usual one (that is the degree of the monomial $x^py^q$ is $p+q$). We shall show in this paper that relaxing this requirement on the valuation leads to many new models.  We have no claim to exhaustivity, and for the moment have no clue about a possible scheme which would lead to a complete classification for the general valuation. However, the domains that we describe here all share some common algebraic properties that we want to underline. 

The construction of these domains rely  mainly on the study of finite subgroups of $O(3)$, and are in particular   related to the Platonic polyhedra.  It relies on the study of invariant polynomials for subgroups of $O(3)$. The analysis of these invariants also lead to the construction of new polynomial models in dimension 3.

\section{ Orthogonal polynomials and diffusion operators}

  The short description of  diffusion operators that we present below is inspired from~\cite{bglbook}, and we refer the reader to it for further details. 

Diffusion operators are second order differential operators with no zero order terms, and  are central in the study of diffusion processes, solutions of stochastic differential equations, Riemannian geometry, classical partial differential equations, potential theory, and many other areas.
When they have smooth coefficients, they may be described in some open subset $\Omega$ of $\bbR^d$ as 
\beq\label{diffusion.gal} \LL(f) = \sum_{ij} g^{ij}(x) \partial^2_{ij} f+ \sum_i b^i(x)\partial_i (f),\eeq
where the symmetric matrix $(g^{ij})(x)$ is everywhere non negative (the operator $\LL$ is said to be semi-elliptic). 
We are mainly interested here in the case where this operator is symmetric with respect to some probability measure $\mu$, that is when, for any smooth functions $f,g$, compactly supported in $\Omega$, one has 
\beq\label{IPP} \int_\Omega f \LL(g) d\mu = \int_\Omega g \LL(f) d\mu.\eeq
We then say that $\mu$ is a reversible measure for $\LL$, which reflects the fact that, in a probabilistic context, the associated stochastic process has a law which is invariant under time reversal, provided that the law at time $0$ of the process is $\mu$.

When $\mu$ has a smooth positive density $\rho$ with respect to the Lebesgue measure, this symmetry property translates immediately in 
\beq\label{eq.density} b^i (x)= \sum_j \partial_j g^{ij}(x) + \sum_j g^{ij} \partial_j \log \rho,\eeq
which shows a fundamental relation between the coefficients of $\LL$ and the measure $\mu$, and allows in general to completely determine $\mu$ up to some normalizing constant.

Let us introduce the carré du champ operator $\Gamma$. For this, we suppose that we have in $\cL^2(\mu)$ some dense algebra $\cA$ of functions  which is stable under the operator $\LL$, and contains the constant functions. Then, for $(f,g)\in \cA$,  we define 
\beq\label{def.Gamma} \Gamma(f,f) = \frac{1}{2}( L(fg)-f\LL(g)-g\LL(f)).\eeq

If $\LL$ is given by equation~\eqref{diffusion.gal}, and when the elements of $\cA$ are at least $\cC^2$,   it turns out that
$$\Gamma(f,g)=\sum_{ij} g^{ij} \partial_if \partial_j g,$$ so that $\Gamma$ describes in fact the second order part of $\LL$.  The semi-ellipticity of $\LL$ translates into the fact that $\Gamma(f,f)\geq 0$, for any $f\in \cA$.  
If we apply formula~\eqref{IPP} with $g=1$, we observe that $\int_\Omega \LL f d\mu=0$ for any $f\in \cA$. Then, applying~\eqref{IPP} again, we see immediately that, for any $(f,g)\in \cA$ 
\beq\label{IPP2} \int_\Omega f\LL(g) d\mu = -\int_\Omega \Gamma(f,g) d\mu,\eeq so that the knowledge of $\Gamma$ and $\mu$ describes entirely the operator $\LL$. We call such a triple $(\Omega, \Gamma, \mu)$ a Markov triple, although we should also add the algebra $\cA$.  Thanks to~\eqref{diffusion.gal}, we see that $\LL(x^i)= b^i$ and $\Gamma(x^i, x^j)= g^{ij}$. The operator $\Gamma$ is called the co-metric, and in our system of coordinates is described by a matrix $\Gamma= \big(\Gamma(x^i, x^j)\big)= (g^{ij})$.

In our setting, we shall always chose $\cA$ to be the set of polynomials. Under the conditions that we shall describe  below, we may as well extend $\cA$ to be the set of the restrictions to $\Omega$ of the smooth functions defined in a neighborhood of $\Omega$, but this extension is  useless in what follows.

The fact that $\LL$ is a second order differential operator translates into the  change of variable formulas. Whenever $f=(f_1, \cdots, f_n)\in \cA^n$, and whenever $\Phi(f_1, \cdots, f_n)\in \cA$, for some smooth function $\Phi  : \bbR^n\mapsto \bbR$, then
\beq\label{chain.rule.L}\LL(\Phi(f))= \sum_i \partial_i \Phi(f) \LL(f_i) + \sum_{ij} \partial^2_{ij} \Phi(f) \Gamma(f_i,f_j)\eeq
and also
\beq\label{chain.rule.Gamma} \Gamma(\Phi(f),g) = \sum_i \partial_i \Phi(f) \Gamma(f_i,g).\eeq
When $\cA$ is the algebra of polynomials, this applies in particular for any polynomial function $\Phi$. Indeed, in this context, properties~\eqref{chain.rule.L} and~\eqref{chain.rule.Gamma} are equivalent.

 As long as polynomials are concerned, it may be convenient to use complex coordinates. That is, for a pair of variables $(x,y)$, consider $z= x+iy$ and $\bar z= x-iy$,  using linearity and bilinearity to extend $\LL$ and $\Gamma$ to $z$ and $\bar z$, for example setting $\LL(z) = \LL(x)+ i\LL(y)$, $\Gamma(z, z) = \Gamma(x,x)-\Gamma(y,y)+ 2i\Gamma(x,y)$. Then one may compute $\LL\big( P(z,\bar z)\big)$ and $\Gamma\big(P(z, \bar z), Q(z, \bar z)\big)$ for any pair of polynomials  $P$ and $Q$ in the variables $(z, \bar z)$ using the change of variable formulas~\eqref{chain.rule.L} and \eqref{chain.rule.Gamma}.  One may then come back to the original variables $x$ and $y$ setting $x= (z+\bar z)/2$, $y= (z-\bar z)/(2i)$.

Moreover, we shall restrict our attention to the case where the matrix $(g^{ij})$ is everywhere positive definite, that is when the operator $\LL$ is elliptic.
In this situation, one may expect $\LL$ to have a self adjoint extension (not unique in general), and then look for a spectral decomposition for this self adjoint extension.  We may expect then that the spectrum is discrete, and look for the eigenvectors.

It is quite rare that one may exhibit explicitly  any eigenvalue or eigenvector, and this makes the analysis of such operators quite hard. However, a good situation is when there is a complete $\cL^2(\mu)$ basis formed of polynomial eigenvectors, in which case one may have explicit computation for the eigenvalues and  expect a good description of the eigenvectors (recurrence formulas, generating functions, etc). These polynomials being eigenvectors of a symmetric operator are orthogonal whenever the eigenvalues are different, and this leads to a family of orthogonal polynomials for the invariant measure $\mu$.
 
Unfortunately, this situation does not appear quite often. In dimension 1 for example, up to affine transformations, there are only 3 cases, corresponding to the Jacobi, Laguerre and Hermite polynomials, see for example~\cite{BakryMazet03}.

\begin{enumerate}
\item The Hermite case corresponds to the Gaussian measure  $\frac{e^{-x^2/2}}{\sqrt{2\pi}}\,dx$ on $\bbR$ and to the  Ornstein--Uhlenbeck operator  
$$\LL_{OU}=\frac{d^2}{dx^2}-x\frac{d}{dx}.$$  The Hermite polynomial $H_n$ of degree $n$ satisfy $\LL_{OU} P_n= -nP_n$.

\item The Laguerre polynomials operator correspond to the measure $\mu_a(dx)=C_ax^{a-1}e^{-x}\,dx$ on $(0, \infty)$ , $a>0$, 
and to the Laguerre operator $$L_a=x\frac{d^2}{dx^2}+ (a-x)\frac{d}{dx}.$$ The Laguerre  polynomial $L^{(a)}_n$ with  degree $n$ satisfies $\LL_a L^{(a)}_n= -nL_n^{(a)}$.

\item The Jacobi polynomials  correspond to the measure $\mu_{a,b}(dx)=C_{a,b}(1-x)^{a-1}(1+x)^{b-1}\,dx$ on $(-1,1)$, $a,b>0$ and to the Jacobi operator
$$\LL_{a,b}=(1-x^2)\frac{d^2}{dx^2}- \big(a-b+(a+b)x\big)\frac{d}{dx}.$$
The Jacobi polynomial  $(J^{(a,b)}_n)_n$ with degree $n$  satisfy $$\LL_{a,b}J_n^{(a,b)}= - n(n+a+b-1)J_n^{(a,b)}.$$
\end{enumerate}

In this paper, we  concentrate on probability measures on  bounded domains $\Omega\subset\bbR^d$. For such measures, the set of polynomials is dense in $\cL^2(\mu)$, and we want to construct bases of $\cL^2(\mu)$ formed with polynomials. There is not an unique choice for such a basis. 

First, we choose a valuation. That is, choosing some positive integers $a_1, \cdots, a_d$, we decide that the degree of a monomial 
$x_1^{p_1}\cdots x_d^{p_d}$ is $a_1p_1+ \cdots a_dp_d$. Then, the degree of a polynomial is the maximum of the degrees of its monomials.

This being done, for $n\in \bbN$, we look at the finite dimensional vector space $\cH_n$ of polynomials with degrees less than $n$. One has 
$\cH_{n}\subset \cH_{n+1}$, and $\cup_n \cH_n$ is the vector space of polynomials. It is dense in $\cL^{2}(\mu)$. Then, a polynomial basis is a choice, for any $n$,   of an orthonormal basis is the orthogonal complement of $\cH_{n+1}$ in $\cH_n$.

Our problem is then to describe for which open bounded subsets $\Omega\subset\bbR^d$, one may find a probability measure $\mu$ on it with positive density $\rho(x)$ with respect to the Lebesgue measure, and an elliptic diffusion operator $\LL$ on $\Omega$ such that such a polynomial basis for $\mu$ is made of eigenvectors for $\LL$, for some given valuation. We restrict our attention to those sets $\Omega$ with piecewise smooth boundary.  Let us call  such a set $\Omega$ a polynomial set, and the triple $(\Omega, \Gamma, \mu)$ a polynomial model.

We recall here some of the results of~\cite{BOZ2013}, where the same structure  is described only for the usual valuation (that is when all the integers $a_i$ are equal to 1), but easily extended to the general valuation case. We have

\bprop\label{prop.poly.model.gal} Choose a  valuation $\deg$  described as above by some integer parameters $(a_1, \cdots, a_d)$, and  let $(\Omega, \Gamma, \mu)$  be a polynomial model in $\bbR^d$. Then,  with $\LL$ described by equation~ \eqref{diffusion.gal},
\benum

\item\label{prop1.pnt1} for $i= 1, \cdots, d$, $b^i$ is a polynomial with $\deg(b^i)\leq a_i$;

\item \label{prop1.pnt2} for $ i, j= 1, \cdots, d$, $g^{ij}$ is a polynomial with $\deg(g^{ij})\leq a_i+a_j$;

\item \label{prop1.pnt3} the boundary $\partial\Omega$ is included in the algebraic set $\{\det(g^{ij})=0\}$;

\item \label{prop1.pnt4}  if $P_1\cdots P_k=0$ is the reduced equation of the boundary $\partial\Omega$ (see remark~\ref{rmk.reduced.bdry.eq} below), then, for each $q= 1, \cdots k$, each $i= 1, \cdots d$, one has 
\beq\label{boundary.eq}\Gamma( \log P_q, x_i)= L_{i,q},\eeq where 
$L_{i,q}$ is a polynomial with $\deg(L_{i,q})\leq a_i$;

\item\label{prop1.pnt5}  all the measures $\mu_{\alpha_1, \cdots, \alpha_k}$ with densities $C_{\alpha_1, \cdots, \alpha_k}|P_1|^{\alpha_1}\cdots |P_k|^{\alpha_k}$ on $\Omega$, where the $\alpha_i$ are such that the  density is  is integrable on $\Omega$, are such that $(\Omega, \Gamma, \mu_{\alpha_1, \cdots, \alpha_k})$ is a polynomial model;

\item\label{prop1.pnt6}  when the degree of $P_1\cdots P_k$ is equal to the degree of $\det(g^{ij})$ there are no other measures.

\eenum

Conversely, assume that some bounded domain $\Omega$ is such that the boundary $\partial\Omega$ is included in an algebraic surface and has reduced equation $P_1\cdots P_k=0$. Assume moreover that there exists a matrix $(g^{ij}(x))$ which is positive definite in $\Omega$ and such that each component $g^{ij}(x)$ is a polynomial with degree at most $a_i+a_j$.  Let $\Gamma$ denote the associated  carré du champ operator. Assume moreover that equation~\eqref{boundary.eq} is satisfied for any $i= 1, \cdots, d$ and any $q= 1,\cdots, k$, with $L_{i,q}$ a polynomial with degree at most $a_i$.

 Let  
$(\alpha_1, \cdots, \alpha_k)$ be such  that  the $|P_1|^{\alpha_1}\cdots |P_k|^{\alpha_k}$ is integrable  on $\Omega$ with respect to the Lebesgue measure, and denote $\mu_{\alpha_1, \cdots, \alpha_k} (dx)= C_{\alpha_1, \cdots, \alpha_k}P_1^{\alpha_1}\cdots P_k^{\alpha_k}dx$, where $C_{\alpha_1, \cdots, \alpha_k}$ is the normalizing constant such that $\mu_{\alpha_1, \cdots, \alpha_k}$ is a probablity measure.

 Then  $(\Omega, \Gamma,\mu_{\alpha_1, \cdots, \alpha_k})$ is a polynomial model.

\eprop
 Before giving a sketch of the proof of Proposition~\ref{prop.poly.model.gal}, let us make  a few remarks.
\brmq\label{rmk.reduced.bdry.eq} We say that $P_1\cdots P_k=0$ is the reduced equation of $\partial\Omega$ when
\benum
\item The polynomials $P_i$ are not proportional to each other.
\item For $i= 1, \cdots k$, $P_i$ is an irreducible polynomial, both in the real and the complex field.

\item For each $i= 1, \cdots, k$, there exists at least one regular point of the boundary $\partial\Omega$ such that $P_i(x)= 0$.

\item For each regular point $x\in \partial\Omega$, there exist a neighborhood $\cV$  and of $x$  and some $i$ such that 
$\partial\Omega\cap \cV= \{P_i(x)=0\}\cap \cV$.

\eenum

In particular, this does not mean that any point satisfying $P_i(x)=0$ for some $i$ belongs to $\partial\Omega$.

\ermq

\brmq The determination of the polynomial models therefore amounts to the determination of the domains $\Omega$ with an algebraic boundary, with the property that  the reduced equation of $\partial\Omega$ is such that the set of equations~\eqref{boundary.eq} has a non trivial (and even positive definite concerning $(g^{ij})$) solution, for $g^{ij}$ and $L_{i,q}$. Looking at the form of these equations, given the reduced equation of $\partial\Omega$,  they appear as a linear homogeneous  equation in the coefficients of the polynomials $g^{ij}$ and of the polynomials $L_{i,k}$. Unfortunately, there are in general much more equations to be satisfied that unknowns, and this requires very strong constraints on the polynomials appearing in the reduced equation of the boundary.

\ermq

\brmq The set of equations~\eqref{boundary.eq}, which are central in the study of polynomial models, may be reduced to less equations, when $k>1$.  Indeed, if we set $P= P_1\cdots P_k$, it reduces to 
\beq\label{boundary.eq2}\Gamma(x_i,\log P)= L_i, \quad \deg(L_i)\leq a_i.\eeq To see this, assume that this last  equation holds with some polynomial $L_i$.  Then  on the regular part of the boundary described by $\{P_q(x)=0\}$,  we have    $\Gamma(x_i,  P_q)= 0$, since 
$$\Gamma(x_i, P_q)= P_q(L_i-\sum_{l\neq q} \frac{\Gamma(x_i,P_l)}{P_l}).$$
Therefore, $P_q$ divides $\Gamma(x_i, P_q)$.

\ermq

\bpf (Of Proposition~\ref{prop.poly.model.gal}).

We shall be a bit sketchy in the details, all  the arguments being borrowed from~\cite{BOZ2013}. 
Let $\cH_n$ be the finite dimensional vector space of polynomials  $P$ such that $\deg(P)\leq n$. From the definition of a polynomial model, $\LL(\cH_n)\subset \cH_n$.  In the representation~\eqref{diffusion.gal} of $\LL$, we have $b^i= \LL(x_i)$ and $g^{ij}= \Gamma(x_i,x_j)$. Therefore, $b^i\in \cH_{a_i}$ and, from the representation~\eqref{def.Gamma} of $\Gamma$, $g^{ij} \in \cH_{a_i+a_j}$.  This gives items~\ref{prop1.pnt1} and~\ref{prop1.pnt2}.

Now, since $\LL$ has polynomial eigenvectors, for any pair $(P,Q)$ of polynomials, we have $$\int P\LL(Q) d\mu= \int Q\LL(P) d\mu.$$ Since the coefficients $g^{ij}$ and $b^i$ are bounded on $\Omega$ with bounded coefficients, this identity may be extended to any pair $(f,g)$ of smooth functions compactly supported in $\bbR^d$ (not necessary  with support in $\Omega$). Looking at this for  smooth functions compactly supported in $\Omega$ leads to
 equation~\eqref{eq.density}, which is equivalent to the symmetry property for such functions. Furthermore, applying this symmetry property to a pair of smooth function compactly supported in a neighborhood of a regular point of the boundary, and  using Stokes formula,  this implies in fact that, for any $i= 1, \cdots, d$, $\sum_j g^{ij} n_j=0$ at any point of the boundary, where $(n_i)$ is the normal vector to the boundary at that point. Therefore, this normal vector is in the kernel of the matrix $g$ at any regular point of the boundary, which implies that $\det(g)=0$ at such a point. This gives item~\ref{prop1.pnt3}.
 
 We now know that the boundary is included in the algebraic set $\{\det(g)=0\}$, and we may look at the reduced equation for it, say $P_1\cdots P_k= P= 0$.  Let $x$ be a regular point of the boundary and $\cV$ a neighborhood of it such that $\partial\Omega\cap \cV= \{P_q=0\}\cap \cV$, for some $q= 1, \cdots, k$. In $\cV$,  the normal vector $(n_i)$  to the boundary is parallel to $\partial_i P_q$, so that we also have for all $i$, 
 $\sum_j g^{ij} \partial_j P_q= 0$ on $\{P_q=0\}\cap \cV$. But $\sum_j g^{ij} \partial_j P_q$ is a  polynomial, which vanishes in $\cV$ on the zeros of $P_q$. This implies (since $P_q$ is complex irreducible) that  
 \beq\label{boundary.eq2}\sum_j g^{ij} \partial_j P_q= L_{i,q} P_q,\eeq where $L_{i,q}$ is a polynomial, the degree of which is less than $a_i$ since $\deg( \sum_j g^{ij} \partial_j P_q)\leq \deg(P_q)+ a_i$.   Then,  equation~\eqref{boundary.eq} is just a rephrasing of~\eqref{boundary.eq2}. This gives item~\ref{prop1.pnt4}.

 If we now apply equation~\eqref{boundary.eq} and look at the value of $b^i$ given by formula~\eqref{eq.density}, we see that, when the measure is $\mu_{\alpha_1, \cdots, \alpha_k}$,
 $$b^i= \sum_i \partial_j g^{ij} + \sum_{k}  \alpha_k L_{i,k},$$ and therefore is a polynomial with $\deg(b^i)\leq a_i$.
 
 Therefore, for every $n\in \bbN$,  the associated operator maps $\cH_n$ into $\cH_n$. Moreover, the boundary equation~\eqref{boundary.eq} shows that for any pair of smooth functions compactly supported in $\bbR^d$, for the associted operator $\LL_{\alpha_1, \cdots, \alpha_k}$
 $$\int f \LL_{\alpha_1, \cdots, \alpha_k} (h) d\mu_{\alpha_1, \cdots, \alpha_k}= \int h \LL_{\alpha_1, \cdots, \alpha_k} (f) d\mu_{\alpha_1, \cdots, \alpha_k},$$ and this in particular applies for polynomials. Therefore, the operator $\LL_{\alpha_1, \cdots, \alpha_k}$ is symmetric on the finite dimensional space $\cH_n$, and this allows to construct a basis of eigenvectors for $\LL_{\alpha_1, \cdots, \alpha_k}$ made of orthogonal polynomials. This gives item~\ref{prop1.pnt5}.
 
 The last item~\ref{prop1.pnt6}, that we shall not use in this note,  is more technical, and relies on the fact that, looking at equation~\eqref{eq.density}, any density measure $\rho$ is such that $\partial_i\log\rho$ is a rational function, with singularities concentrated  on $\{\det(g)=0\}$, and   degree $-a_i$. We refer to~\cite{BOZ2013}, where the arguments are developed, and which furthermore provides  a complete  description of the possible measures in the case where the reduced equation of $\partial\Omega$ is not $\{\det(g)=0\}$.
 
 The proof of the reverse part of Proposition~\ref{prop.poly.model.gal} is just a rephrasing of that of item~\ref{prop1.pnt5}.

\epf

From Proposition~\ref{prop.poly.model.gal}, the important data are the set $\Omega$ (open subset of $\bbR^d$, bounded with piecewise smooth boundary given by an algebraic reduced equation $P_1\cdots P_k=0$), and the operator $\Gamma$, given by polynomial functions $(g^{ij})$, elliptic in $\Omega$,   satisfying the degree condition~\ref{prop1.pnt2},   and the boundary equation~\ref{prop1.pnt4}. 

To fix the ideas, we provide a few definitions

\bdefi~
\benum
\item A polynomial domain $\Omega\subset \bbR^d$  is a bounded open set in $\bbR^d$ with boundary $\partial\Omega$ included in some algebraic surface with reduced equation $\{P(x)=0\}$, and such that there exists some valuation $\{a_1, \cdots, a_d\}$ and some elliptic co-metric $\Gamma= (g^{ij})$ on $\Omega$ with $\deg(g^{ij}) \leq a_i+a_j$ satisfying the boundary equation~\eqref{boundary.eq2}.

\item A polynomial system $(\Omega, \Gamma)$ is given by a polynomial domain $\Omega$ together with the associated co-metric $\Gamma$.

\item A polynomial model is a triple $(\Omega, \Gamma, \mu),$ where  $(\Omega, \Gamma)$  is a polynomial system, and $\mu$ is  probability measure on $\Omega$  with smooth density $\rho$ such that $\Gamma(x^i, \log \rho)= S_i$, with $\deg(S_i) \leq a_i$.

\eenum

\edefi

By definition, to each polynomial domain corresponds at least one polynomial system (there may indeed be many different co-metrics $\Gamma$ associated with the same domain $\Omega$, see~\cite{BOZ2013}). Moreover, we saw that to any polynomial system corresponds many polynomial models.

\brmq The valuation is not unique. Beyond the trivial change $(a_1, \cdots, a_d)\mapsto (ca_1, \cdots, ca_q)$, the same polynomial model (or system) may correspond to various valuations. We shall make no effort to provide the lowest ones since in general a good choice is provided by a simple look at the co-metric $\Gamma$.

\ermq

 In~\cite{BOZ2013}, a complete description of all polynomial models is provided when the chosen  valuation is the natural one (we give this description in Section~\ref{sec.BOZmodels} at the end of the paper for completeness). This description relies in an essential way on the fact that for the natural valuation,  the  problem is affine invariant, that is that a polynomial domain $\Omega$ is transformed into another one through affine transformations. This allows for an analysis of the boundary equation, and to the classification of algebraic curves  in the plane for which the boundary equation~\eqref{boundary.eq} has a non trivial solution, through the analysis of the singular points of the curve and its dual.

 This affine invariance is lost for other valuations, since an affine transformation no longer maps the set $\cH_n$ of polynomials with degree at most $n$ into itself. This paves the way for the construction of new  models. In what follows, we shall mainly concentrate on the construction of polynomial systems in dimension 2. These two dimensional models also provides new 3-d models through the use of 2-fold covers of our 2-d models. These two fold covers already appear in~\cite{BOZ2013}. But even for the 2-d models which already appear there ($\Omega_{11}$ and $\Omega_{13}$ of Section~\ref{sec.cube}, e.g.), some two-fold coverings appear as new. The reason is that in~\cite{BOZ2013}, only the models with natural valuation are considered. Here, even though the 2-d models  may be considered with the usual valuation, this is no longer the case for their coverings.

\section{Constructing  polynomial systems\label{sec.new.models}}

  A generic way for the construction of polynomial models in dimension $d$ is to consider some other  symmetric diffusion operator   $\LL$ (often in higher dimension)  and look for functions $(X_1, \cdots, X_d)$ such that, setting $X= (X_1, \cdots, X_d)$
$$\LL(X_i)= B^i(X), ~\Gamma(X_i,X_j)= G^{ij}(X),$$ where $B^i$ and $G^{ij}$ are some smooth functions. Then according to formula~\eqref{chain.rule.L},
$$\LL(\Phi(X))= \hat \LL (\Phi) (X),$$ where 
$$\hat \LL \Phi = \sum_{ij} G^{ij}(X) \partial^2_{ij} + \sum_i B^i(X)\partial_i.$$
This new operator $\hat \LL$ has as reversible measure $\hat \mu$ which is the image of the reversible measure $\mu$  of $\LL$ under $X$. This is often a good way to identify image measures, through equation~\eqref{eq.density}. Then, $\hat \LL$ corresponds to a new triple $(\hat \Omega, \hat \Gamma, \hat \mu)$, where $\hat \Omega$ is the image $X(\Omega)$,  $\hat\Gamma= (G^{ij})$ and $\hat \mu$ is the image of $\mu$.

\bdefi~
\benum
\item When we have such functions $(X_1, \cdots, X_d)$ such that  $\Gamma(X_i, X_j)= G^{ij}(X)$, we say that $(X_1, \cdots, X_d)$ form a closed system for $\Gamma$.

\item  If moreover $\LL(X_i)= B^i(X)$, we say that we have a closed system for $\LL$. 
\eenum
\edefi 

It may   happen that for some specific polynomial model  $(\Omega, \Gamma, \mu)$ and some functions  $X=(X_1, \cdots, X_d)$,  $X$  is a closed system  for $\Gamma$, but not for $\LL$. 

Now, if $\LL$ itself maps polynomials with degree $n$  into  polynomials with degree $n$ (say with the usual valuation), if $X_i$ is a polynomial with degree $a_i$, and  if $B^i(X)$ and $G^{ij}(X)$ are polynomials in $X$, then $\hat \LL$ provide a next polynomial model with valuation $(a_1, \cdots, a_d)$.

It may  also happen that this transformation $x\mapsto X= (X_1, \cdots, X_d)$ is a diffeomorphism, in which case we do not distinguish between those  two models.  If this diffeomorphism and its inverse are given through polynomial transformations,  and if both are polynomial systems or models,  we say that these systems or models are isomorphic. It is not always easy to see when a model is an image of another one, or when they are isomorphic.

Apart of one specific case  (example~7 of Section~\ref{sec.BOZmodels}),   all the models which appears in~\cite{BOZ2013} may be constructed either from the Euclidean Laplace operator in $\bbR^2$ acting on function invariant under the symmetries of a regular lattice (examples 1, 6 and 11 in Section~\ref{sec.BOZmodels}), or from the spherical Laplace operator on the unit sphere $\bbS^2\subset\bbR^3$ acting on functions which are invariant under some finite subgroup of $O(3)$ (all the other models of Section~\ref{sec.BOZmodels}).  Here, we shall explore in a systematic way all the models that one may construct from the finite subgroups of $O(3)$.    This construction may be  also carried  in higher dimension letting the spherical Laplace operator on $\bbS^{d-1}$ act on polynomials in $\bbR^d$ (that is on the restriction to $\bbS^{d-1}$ of such polynomials).

The spherical Laplace operator on $\bbS^{d-1}$ may be described through its action on linear forms. If $e$ is any vector in the Euclidean space $\bbR^d$, we look at  the associated linear form $e^* : x\mapsto e\cdot x$, and more precisely  to its  restriction to the unit sphere, as a function $\bbS^{d-1}
\mapsto \bbR$. Then,  for the Laplace operator $\LL^{\bbS}$ and its associated carré du champ $\Gamma^{\bbS}$, we have
\beq\label{lapl.sphere}  \LL^{\bbS}(e^*)= -(d-1)e^*, ~\Gamma^{\bbS}(e^*, f^*)= e\cdot f-e^*f^*.\eeq

Therefore, looking at the canonical basis $(e_i)$ of $\bbR^d$, we see that any polynomial in the variables $x_i$  $(= e_i^*)$ is transformed under $\LL^{\bbS}$ into a polynomial with the same (natural)  degree.  Moreover, the spherical Laplace operator commutes with all the elements of $O(d)$. Then, if we are   given any subgroup of $O(d)$ and if we  look at the set of polynomials invariants under the group action,  $\LL^{\bbS}$ will preserve this set. If we  may describe some  polynomial basis for these invariant polynomials, then we  expect to get in such a way a closed  system, and therefore construct new polynomial models.

\section{Invariant polynomials\label{sec.invariants}}

The theory of invariant polynomials has a long history going back to D. Hilbert, E. Noether, etc. It now plays an important role in coding theory and combinatorics (see \cite{Stanley79}).  In what follows, we provide a brief account which is useful for the understanding of our construction method, reducing to the case  of finite subgroups of $O(n)$.  We refer to {\cite{ Smithbook95, SmithAMS97}  for further details.

Given any finite subgroup $G$ of $O(n)$, any element  $g\in G$ acts on the set of linear functionals $(x_1, \cdots, x_n)$.  We may consider its action on homogeneous polynomials in the variables $(x_1, \cdots, x_n)$ and look for invariant polynomials, that are homogeneous polynomials which are invariant under the group action. If one denotes by $d_n$ the dimension of the vector space of invariant polynomials with homogeneous degree $n$, then  Molien's formula  allows to compute the Hilbert sum $F(G,t)= \sum_n d_n t^n$ through 
\beq \label{eq.mollien} F(G,t) = \frac{1}{|G|} \sum_{g\in G} \frac{1}{\det(\Id-tg)}.\eeq

Moreover, the set of invariant polynomials may be represented as follows. First, there exist $n$ algebraically independent polynomials $(\theta_1, \cdots, \theta_n)$, called primary invariants, and some other invariant polynomials $(\eta_1, \cdots, \eta_k)$ (the number of them may depend on the choice of the $\theta_i$),  called secondary invariants, such that any invariant may be written as 
$$P_0(\theta_1, \cdots, \theta_n)+ \sum_{i=1}^k \eta_i P_i(\theta_1, \cdots, \theta_n),$$ where $P_i$ are polynomials (in the variables $(\theta_1, \cdots, \theta_n)$). Moreover, each $\eta_i$ satisfies some monic polynomial equation in the variables $\theta=(\theta_j)$, that is satisfies an  algebraic identity of the form
$$\eta_i^{p_i}+ \eta_i^{p_i-1} Q_{i,1}(\theta)+ \cdots +Q_{i, p_i}(\theta)=0,$$
where $Q_{i,k}(\theta)$ are polynomials in the variables $(\theta_1, \cdots, \theta_n)$. These algebraic relations are called syzygies.

Furthermore, there are only primary generators if and only if the group $G$ is generated by reflections, that is when  $G$ is a Coxeter group.

In order to construct polynomial systems, we then consider finite subgroups of $O(n)$,  compute their invariants (primary and secondary when they exist),  look at their restriction to the unit sphere (that is consider those polynomials modulo $\sum_i x_i^2-1$). They are no longer homogeneous, and, since $\sum_i x_i^2$ may always be considered as a primary invariant, we may reduce to $n-1$ primary invariants, plus some number of secondary invariants. We then let the spherical Laplace operator act on them. Since the spherical Laplace operator commutes with rotations, it preserves the set of invariant polynomials. Moreover, it maps polynomials with degree $n$ into polynomials with the same degree. Constructing such polynomial systems amounts then to  choose some family  $(\zeta_1, \cdots, \zeta_p)$ of invariants, and look for $\Gamma(\zeta_i, \zeta_j)$, expecting that it may be written as $G^{ij}(\zeta_1, \cdots, \zeta_p)$ (that is to provide a closed system for $\Gamma$). Then, the extra condition on the degrees will be automatically satisfied, where the valuation is defined through $a_i=\deg(\xi_i)$.  The difficulty then is to produce such a closed system $(\zeta_i)$  of algebraically independent polynomials. When such happens, we produce a polynomial system which is an image of the starting Laplace operator.

In all the examples in dimension $3$, one may always chose 2 primary invariant $(\theta_1, \theta_2)$ to produce a closed system (this is no longer true in higher dimension, see Section~\ref{sec.further}). Moreover, when one adds one secondary invariant $\eta$, we always obtain a closed system with 3 variables $(\theta_1,\theta_2, \eta)= (\zeta_1, \zeta_2, \zeta_3)$.  Now, it turns out  that, if one forgets about the algebraic relations $Q(\theta_1, \theta_2, \eta)=0$, and consider the polynomials $G^{ij}(\zeta_i, \zeta_j)= \Gamma(\zeta_i, \zeta_j)$ as a polynomial co-metric in dimension $3$,  it provides a polynomial model on a bounded domain in $\bbR^3$ which has the surface  $Q(\zeta_1, \zeta_2, \zeta_3)=0$ (the syzygy) as a  part of its boundary. Although the first construction with just the primary invariants $(\theta_1, \theta_2)$ is not surprising (all our groups are sub-groups of Coxeter groups), the second property (construction of 3-dimensional models from 2-dimensional ones through the syzygies) remains quite mysterious. 

Let us show this phenomenon  in dimension 1, on the simpler form of the cyclic group $\bbZ_n$ acting in $\bbR^2$, as a rotation in the complex plane with angle $2\pi/n$. Writing $z= x+iy$, we may choose as primary invariant $X= \Re(z^n)$, and secondary invariant $Y=\Im(z^n)$. We now restrict them to the unit circle $\bbS^1$ and let the spherical Laplace operator act on it. Using  formulas~\eqref{lapl.sphere}, or the computations provided at the beginning of Section~\ref{sec.cyclic}, one sees that 
$$\Gamma(X,X)= n^2(1-X^2), ~\LL(X)= -n^2X,$$ and therefore it provides a closed system for $\LL$ which corresponds (up to the factor $n^2$), to the classical Jacobi operator on $(-1,1)$. Now, if we  add  the variable $Y$, we get again a closed system for $\Gamma$, with co-metric
$$\bpm \Gamma(X,X)& \Gamma(X,Y)\\\Gamma(X,Y)&\Gamma(X,Y)\epm= n^2\bpm 1-X^2& -XY\\-XY& 1-Y^2\epm$$
The determinant of this matrix is $1-X^2-Y^2$, and indeed we have $X^2+Y^2=1$ in our model (this is the syzygy relating $X$ and $Y$). But the metric $\bpm 1-X^2& -XY\\-XY& 1-Y^2\epm$ is a metric on the unit ball $\Omega= \{ 1-X^2-Y^2>0\}\subset \bbR^2$ which corresponds  to the model 2 in the  2-d polynomial models of \cite{BOZ2013} in Section~\ref{sec.BOZmodels}.   The various probability measures for this model have the form $C_a(1-X^2-Y^2)^a dXdY$, whith $a>-1$.   When $a= (d-3)/2$,  for some integer $d\geq 2$, this corresponds to the image of the Laplace operator on $\bbS^d$ through the projection $(x_1, \cdots, x_{d+1})\in \bbR^{d+1}\mapsto (x_1=X, x_2= Y)$.   This measure concentrates when $a\to -1$ to the uniform measure on the boundary $\bbS^1$.  The case that we just described as the image of the Laplace operator on $\bbS^1$ corresponds in this model to a limiting  case when $d\to 1$.  

This is this phenomenon that will remain valid in dimension $2$ in the examples described below, although we will not provide   such  simple geometric interpretation for the various 3 dimensional models constructed from the syzygies.

\section{ Finite subgroups of $O(3)$}

In our context, we shall restrict to finite subgroups of $O(3)$. We first describe them, following \cite{MeyerBurn54}.  There are only five (families of) finite subgroups of $SO(3)$, described by F. Klein, corresponding to the cyclic, dihedral,  tetrahedral, octahedral and icosahedral respectively,  denoted  in what follows as $\cC_n, \cD_n, \cT, \cO, \cI$ respectively.  The groups $\cT, \cO, \cI$ correspond  to the elements of $O(3)$ preserving respectively the tetrahedron, the octahedron or its dual the cube, the icosahedron or its dual the dodecahedron. 

The finite subgroups of $O(3)$ are described in two ways.  The first class is obtained adding the central symmetry $J: x\mapsto -x$ to one of the subgroups of $SO(3)$. If $G$ is such a group, we denote $G_J$ this new group, with $|G_J|= 2|G|$.

The second class is obtained by those groups $G$ of $SO(3)$ which contain a subgroup  $G_1$ of  index 2. A new group  denoted $G_1| G$ is obtained as $G_1 \cup \{ Jg, g\in G\setminus G_1\}$.  This provides the groups $\cT |\cO, \cC_n | \cD_n, \cD_n | \cD_{2n}, \cC_n |\cC_{2n}$, where in the case of the cyclic and dihedral groups, the structure of invariants may depend on the fact that $n$ is odd or even. A complete table of  Molien's formulas is provided in 
\cite{MeyerBurn54} together with the associated list of invariants (with however some error in the secondary invariant for the group $\cI$).

In the following sections, we shall  describe the various invariants, and provide the polynomial models which they produce, both in dimension $2$ with the primary invariants, and then in dimension 3 with the use of the secondary ones and their syzygies.

Among the subgroups of $O(3)$, the following are Coxeter groups: $\cD_{n,J}$ ($n$ even) and $\cD_n | \cD_{2n}$ ($n$ odd) ; $\cC_n | \cD_n$, for all $n$ ; $\cT | \cO$, $\cO_J$ and $\cI_J$.  They will yield the primary invariants and hence a closed system and a model. Among them, some were known: those obtained from $\cC_2|\cD_2$ (coaxial parabolas), $\cD_3 | \cD_{6}$ (the cuspidal cubic with secant), $\cT | \cO$ (the swallow tail), $\cO_J$ (the cuspidal cubic with tangent). But $\cD_{n,J}$ ($n$ even) and $\cD_n | \cD_{2n}$ ($n$ odd) for $n$ larger yield new models involving Tchebychev polynomials, not very surprising ; and $\cI_J$ yields a nice model with an angle based on $\pi/5$ whose existence has been the initial motivation of this work.

Since each other subgroup of $O(3)$ is a subgroup of one of these Coxeter groups, we obtain the higher dimensional models  by adding the secondary invariant as auxiliary variable. In most examples, if the equation of the boundary of the two dimensionnal model yield by the Coxeter group is $P(X,Y) =0$, then  the equation of the boundary of the three dimensionnal models are either of the form $Z^2 - P(X,Y) = 0$,  of the form $X(Z^2 - P(X,Y)) = 0$ or of the form $Z^2 - X P(X,Y) = 0$ and the boundary of the three dimensional domain is either a bounded  two leaves cover of the two dimensional domain, either the same but bounded also by a plane. The case of the groups $\cC_n | \cC_{2n}$ or $\cC_{n,J}$ is special in that  we have more than one secondary invariant; each of them yield a different three dimensional system.

From now on, the operator $\Gamma$ will always be the carré du champ operator of the sphere~$\bbS^2$.

\section{Cyclic and dihedral  groups\label{sec.cyclic}}

Let $(x,y,z)$ be the standard coordinate system in $\bbR^3$. On the unit circle $\{z=0\}\cap \bbS^2$, we choose $n$ equidistant points $(e_i, i= 1 , \cdots n)$. The group $\bbZ_n$ acts on them by circular permutations, which consist of elements of $SO(3)$  with vertical axis and angle $2\pi/n$.

 We first consider the complex function $Z= x+iy$, with its conjugate $\bar Z= x-iy$, and observe that 
 $$\Gamma(Z,Z)= -Z^2, \Gamma(\bar Z, \bar Z)= -\bar Z^2, \Gamma(Z, \bar Z)= 2-Z\bar Z,$$
 and 
 $$\Gamma(z,z)= 1-z^2, \Gamma(z, Z)= -zZ, ~\Gamma(z, \bar Z)= -z\bar Z.$$
 With this in hand, we set 
 $$X_n= \Re(Z^n)= \frac{1}{2}(Z^n+ \bar Z^n), ~Y_n = \Im(Z^n)= \frac{1}{2i} (Z^n-\bar Z^n).$$
 The 3 variables $(z,X_n,Y_n)$ are linked by the relation $X_n^2+Y_n^2= (1-z^2)^n$.

 Then, the table
 $$\Gamma= \bpm \Gamma(z,z))&\Gamma(z,X_n)&\Gamma(z,Y_n)\\& \Gamma(X_n,X_n)& \Gamma(X_n,Y_n)\\ && \Gamma(Y_n,Y_n)\epm$$ is given by 
 \beq\label{Gamma3d.Zn}\bpm 1-z^2& -nzX_n&-nzY_n\\
 & n^2((1-z^2)^{n-1}-X_n^2) &-n^2X_nY_n\\
  & & n^2((1-z^2)^{n-1}-Y_n^2
 \epm
 \eeq

We may chose $\theta_1=z, \theta_2=X_n$ as primary invariants; these are the invariants of the Coxeter group $\cC_n| \cD_n$. 
Hence, consider $(\theta_1, \theta_2)= (z,X_n)$.  From table~\eqref{Gamma3d.Zn},  we see that they form a closed system for $\Gamma$. Let $\Gamma_1$ be the extracted matrix corresponding to the two first lines and columns from $\Gamma$. 
$$\Gamma_1=  \bpm 1-\theta_1^2& -n\theta_1\theta_2\\& n^2((1-\theta_1^2)^{n-1}-\theta_2^2)\epm.$$

Up to the factor $n^2$, the determinant of this matrix is $P(\theta_1, \theta_2)= (1-\theta_1^2)^n-\theta_2^2$, and according to $n$ being odd or even, it has 1 or 2 irreducible factors. Then, it is quite immediate to see that 
$$\Gamma_1(\theta_1, \log P)= -2n\theta_1, \Gamma_1(\theta_2, \log P)= -2n^2 \theta_2.$$
\begin{figure}[!h]
\begin{center}
\includegraphics[width=5cm]{Omega1_3}
\includegraphics[width=5cm]{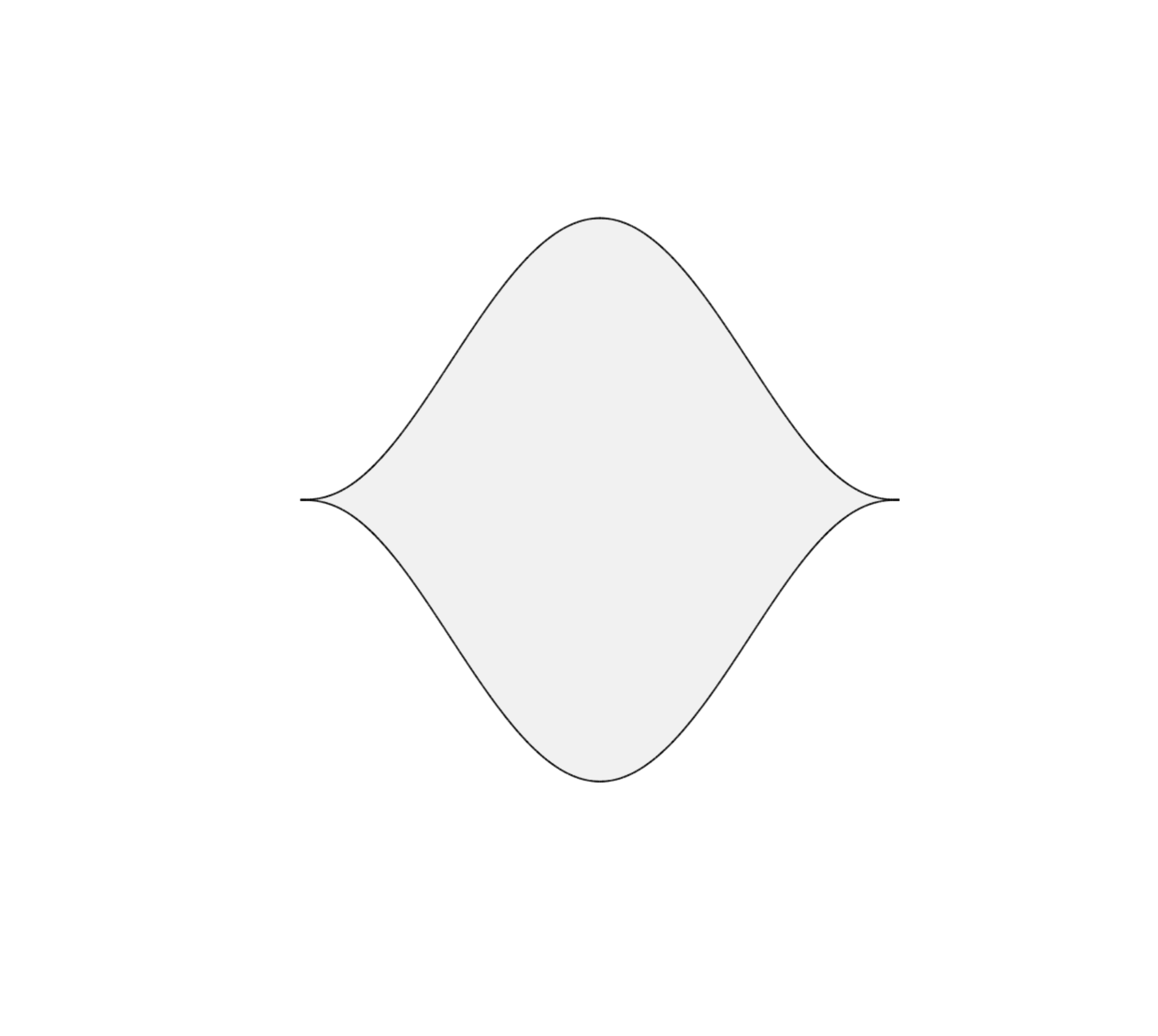}
\includegraphics[width=5cm]{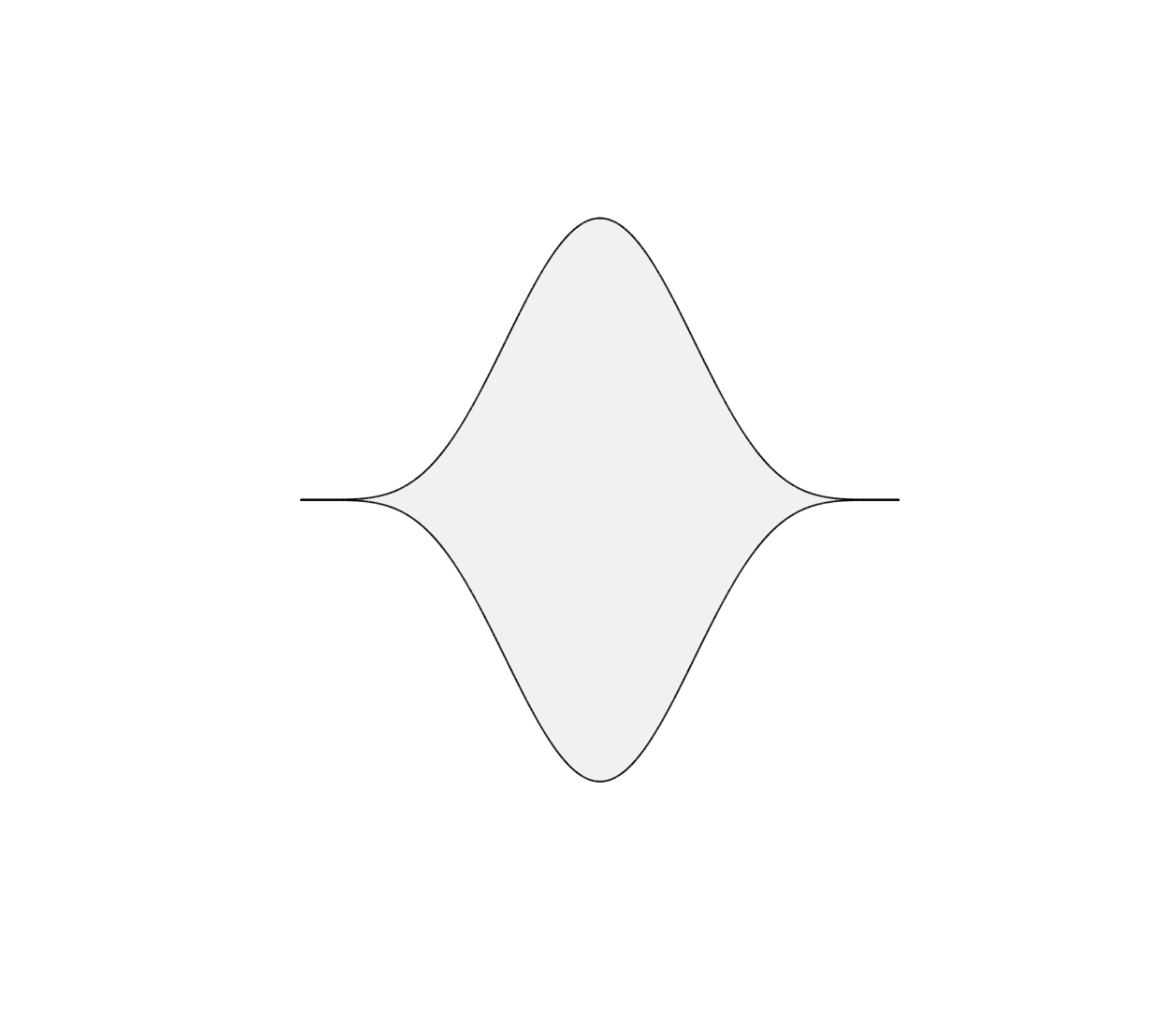}
\caption{$\Omega_1^{(3)}$,   $\Omega_1^{(5)}$ and $\Omega_1^{(11)}$.}
\end{center}
\end{figure}
\begin{figure}[!h]
\begin{center}
\includegraphics[width=6cm]{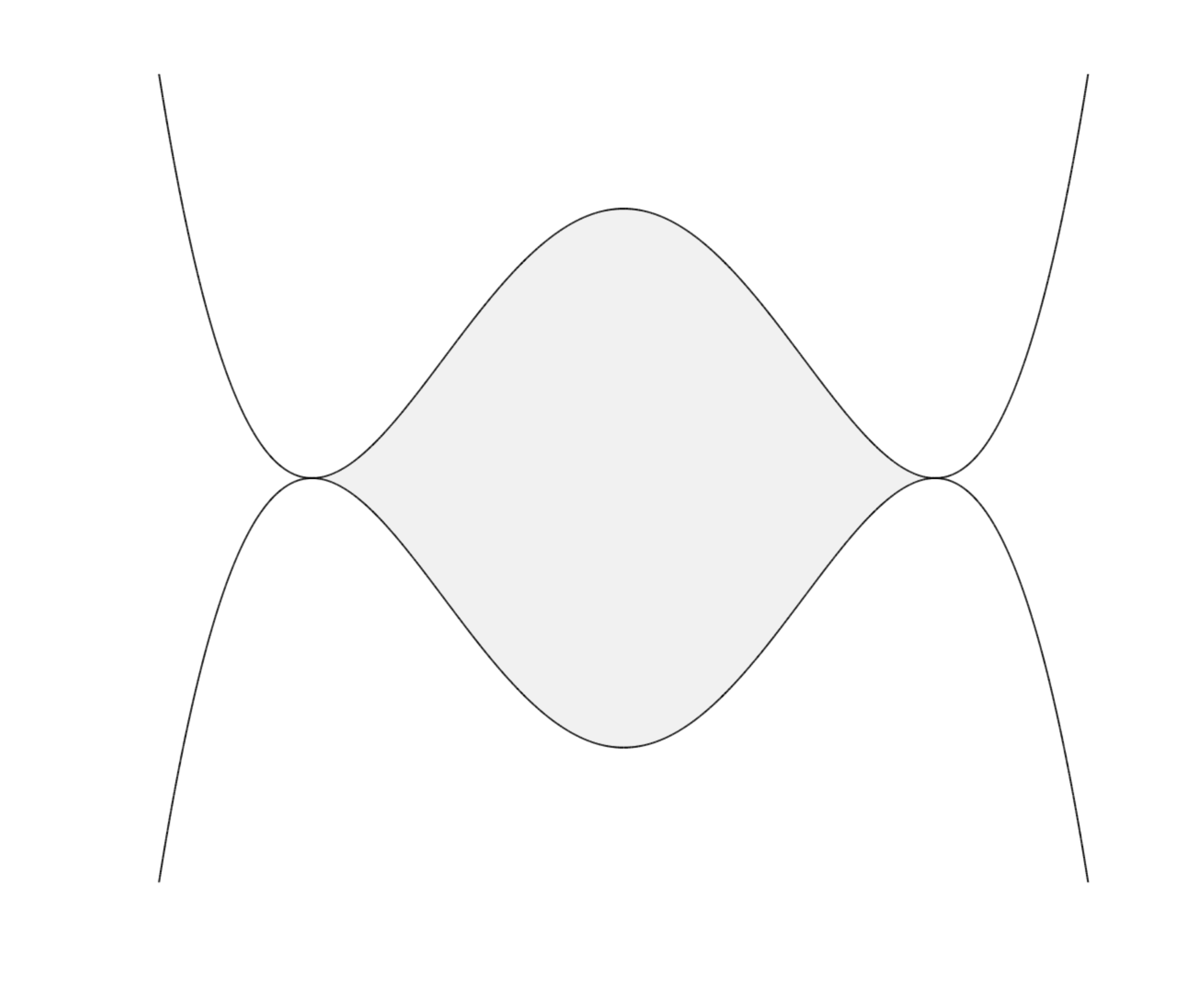}
\caption{$\Omega_1^{(4)}$ to illustrate the even case }
\end{center}
\end{figure}

It satisfies therefore the boundary equation. When $n$ is odd, the  set $\Omega_1^{(n)}= \{ P(\theta_1, \theta_2)>0\}\subset \bbR^2$ is bounded and $(\Omega_1^{(n)}, \Gamma_1)$ provides a  polynomial system.  When $n=2p$ is even, 
$P= P_1P_2$, where $P_1= (1-\theta_1^2)^p-\theta_2$, $P_2=(1-\theta_1^2)^p+\theta_2$. The area   $\Omega_1^{(n)}$ plane with $P_1(\theta_1, \theta_2)>0$, $P_2(\theta_1, \theta_2)>0$ and $\theta_1\in (0,1)$  has  $P_1P_2=0$ as reduced boundary equation, and $(\Omega_1^{(n)}, \Gamma_1)$ is again a polynomial system. In this model, we may chose $\deg(\theta_1)=1$ and $\deg(\theta_2)= n$, which comes from the sphere interpretation, but we may observe that we may as well chose $\deg(\theta_1)=1, \deg(\theta_2)=n-1$. Observe that these domains correspond to the disk if $n=1$ (model 2 in Section~\ref{sec.BOZmodels}) and to the double parabola if $n=2$ (model~4 in Section~\ref{sec.BOZmodels}),  but are new as soon as $n \geq 3$.

We now add the variable $Y_n= \eta$ in the figure, which is our secondary invariant  (observe that the roles of $X_n$ and $Y_n$ are  similar, and we may as well exchange them).  This reflects the symmetries of the cyclic group $\cC_n$. We now have the co-metric 
$$\Gamma_2=\bpm
 1-\theta_1^2& -n\theta_1\theta_2 &-n\theta_1\eta \\& n^2((1-\theta_1^2)^{n-1}-\theta_2^2)& -n^2 \theta_2\eta\\
&&n^2 ((1-\theta_1^2)^{n-1}-\eta ^2)
\epm
$$
whose determinant factorizes as 
$n^4(1-\theta_1^2)^{n-1} ((1-\theta_1^2)^n- \theta_2^2-\eta ^2)$.
Note that the last factor $$P(\theta_1, \theta_2, \eta)=  (1-\theta_1^2)^n- \theta_2^2-\eta ^2$$  reflects the  syzygy relating $\eta$ to $(\theta_1, \theta_2)$.  It is not a surprise that this determinant vanishes identically, since $\Gamma$ is the Gramm matrix of the gradients  (on the sphere) of three functions, and the range of these 3 gradients is at most 2.  Observe also that this syzygy,  which will appear in the boundary of the 3-d system, may be written as $P(\theta_1, \theta_2)-\eta^2$, where $P$ is the boundary equation of the corresponding 2-d system.

But now consider $\Gamma_2$ as a co-metric in $\bbR^3$ on the bounded domain $\Omega^{(n)}_2= \{|\theta_1|<1 ,  P(\theta_1, \theta_2, \eta)>0\} \subset \bbR^3$, which has indeed again reduced equation $P(\theta_1, \theta_2, \eta)=0$.  We may check that
$$\Gamma_2(\theta_1, \log(P))= -2n\theta_1, \Gamma_2(\theta_2, \log(P))= -2n^2\theta_2, \Gamma_2(\eta , \log(P))= -2n^2\eta ,$$ so that 
indeed 
$(\Omega_2^{(n)}, \Gamma_2)$ is a polynomial system in $\bbR^3$, with degrees $\deg(\theta_1)=1, \deg(\theta_2)=n, \deg(\eta)=n$.

\begin{figure}[!h]
\begin{center}
\includegraphics[width=7cm]{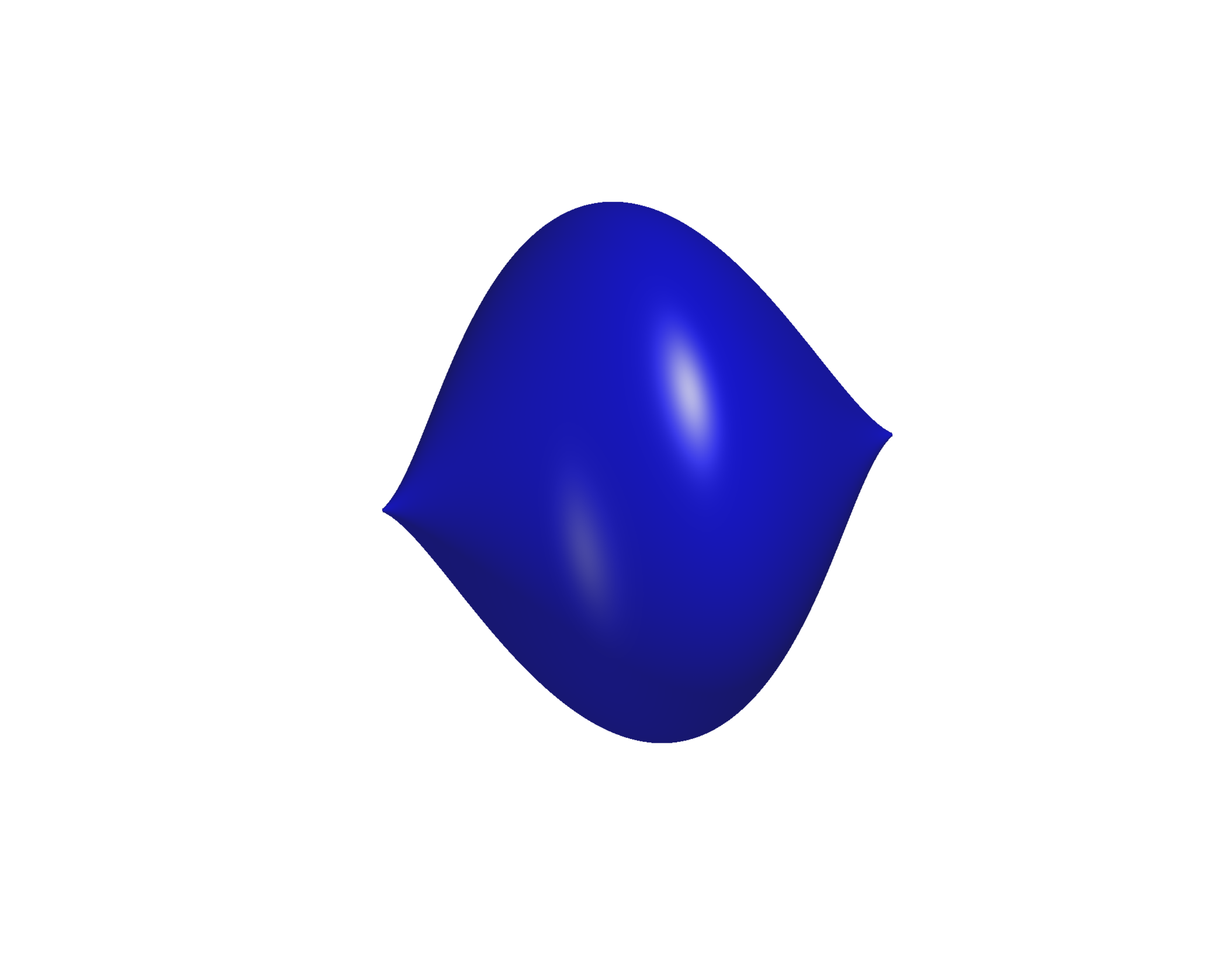}
\includegraphics[width=7.5cm]{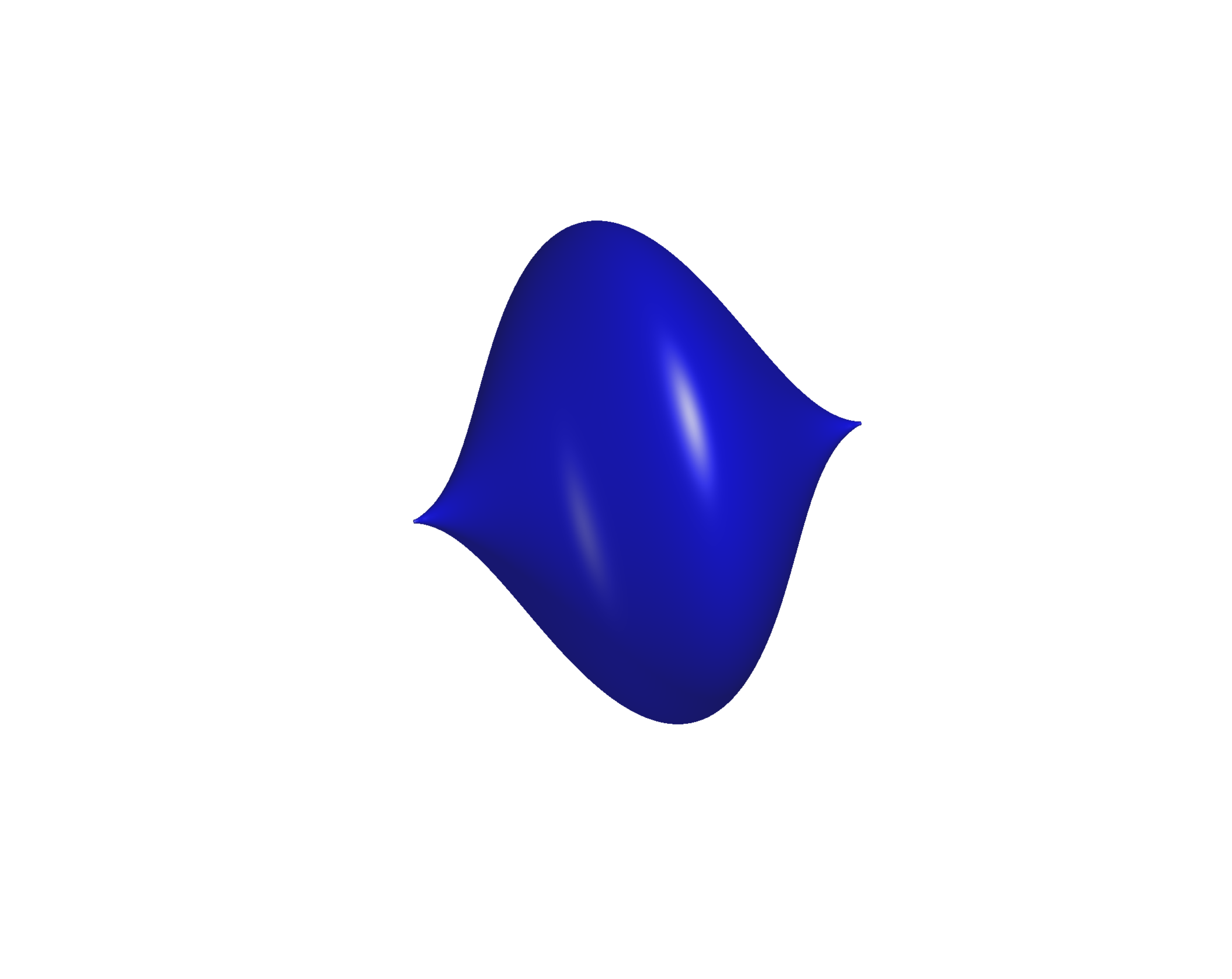}
\caption{Here, the domains $\Omega_2^{(3)}$ and $\Omega_2^{(5)}$.}
\end{center}
\end{figure}

We now study the case where the groups contain the central symmetry.  This corresponds to the new system of primary invariants  $(\theta_1= z^2, \theta_2= X_n)$ associated with the Coxeter group $\cD_{n,J}$ ($n$ even)  or $D_{n}|\cD_{2n}$ ($n$ odd).  
We get 
$$\Gamma_3= \bpm 4 \theta_1(1- \theta_1) & -2n\theta_1 \theta_2\\&n^2((1- \theta_1)^{n-1}- \theta_2^2)\epm$$
which, up to a constant, has determinant $P(\theta_1, \theta_2)=  \theta_1((1-\theta_1)^n- \theta_2^2)= \theta_1P_1(\theta_1, \theta_2)$. It has three irreducible components when $n$ is even and 2 when $n$ is odd. 
Once again 
$$\Gamma_3(\theta_1, \log(\theta_1)) =4(1-\theta_1), \Gamma_3(\theta_2, \log(\theta_1)) =-2n\theta_2,$$
and, 
$$\Gamma_3(\theta_1, \log(P_1)) =-4n\theta_1, \Gamma_3(\theta_2, \log(P_1)) =-2n^2\theta_2,$$

so that the domain $\Omega_3 = \{\theta_1\in (0,1), P_1(\theta_1, \theta_2)>0\}$, which has reduced boundary equation $\theta_1P_1(\theta_1, \theta_2)=0$ is such that 
$(\Omega_3, \Gamma_3)$ is a polynomial system.
\begin{figure}[!h]
\begin{center}
\includegraphics[width=7cm]{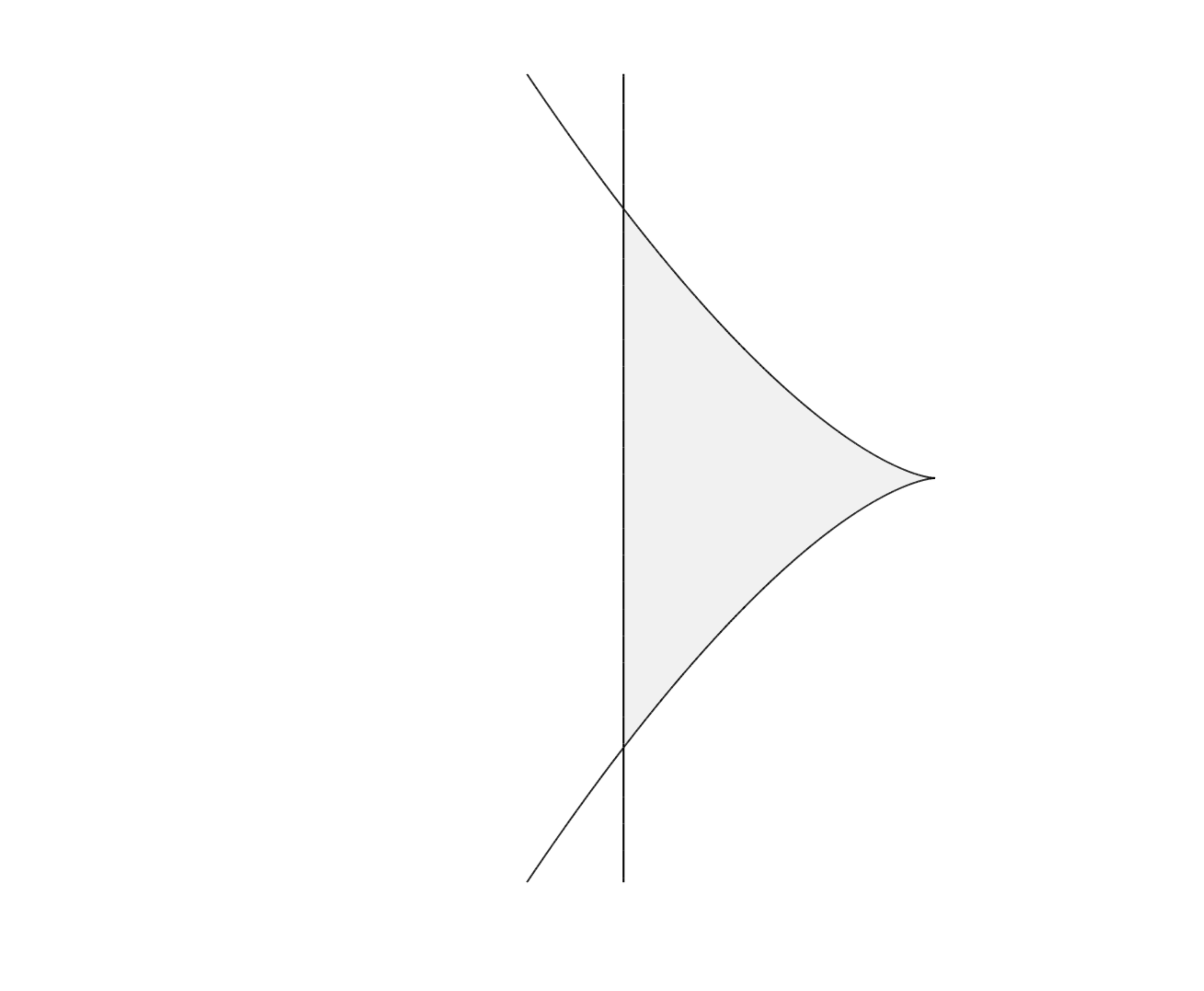}
\includegraphics[width=7cm]{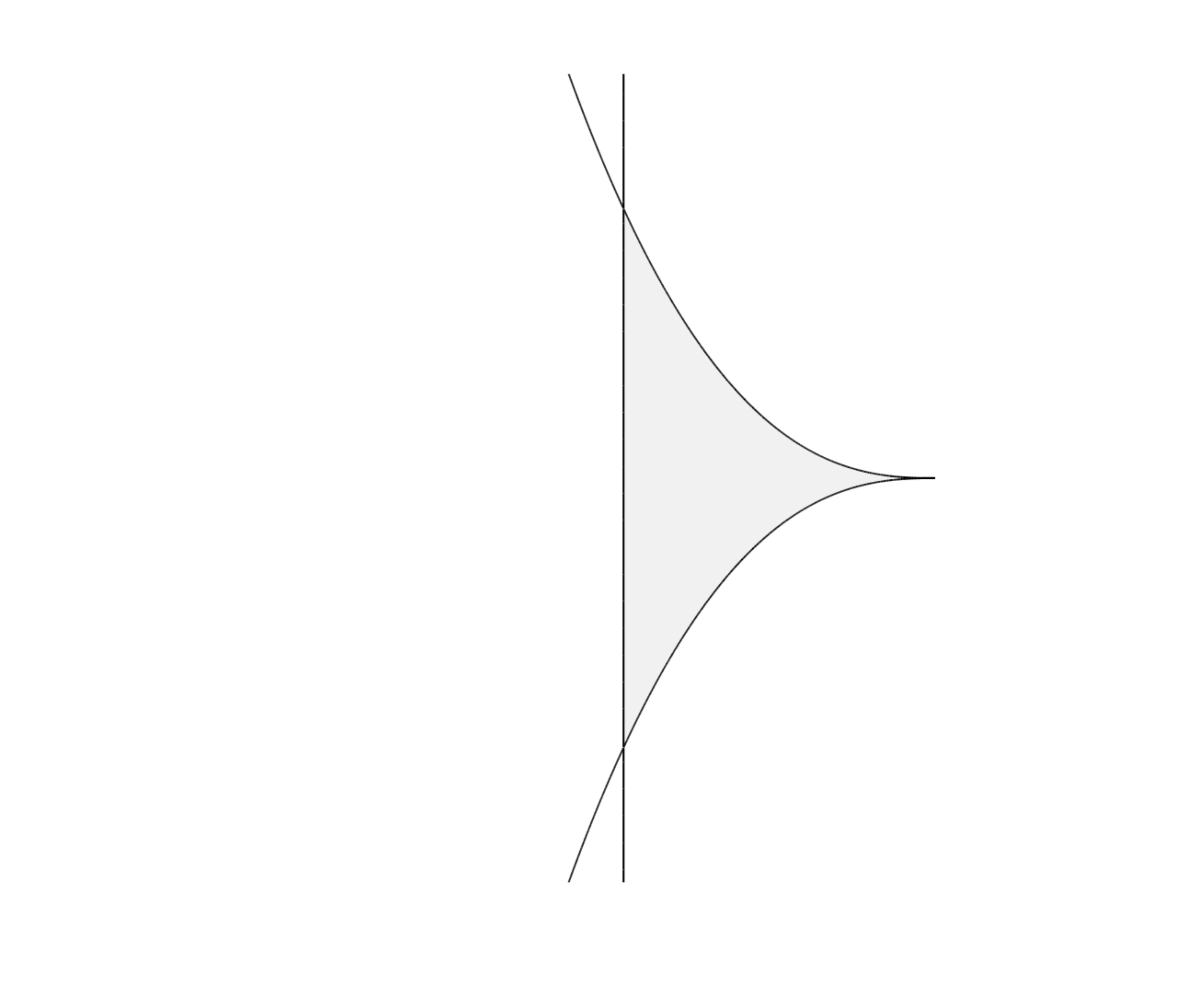}
\caption{The domains $\Omega_3^{(3)}$ and $\Omega_3^{(5)}$ (for $\cC_{3,J}$ and $\cC_{5,J}$)}
\end{center}
\end{figure}
\begin{figure}[!h]
\begin{center}
\includegraphics[width=7cm]{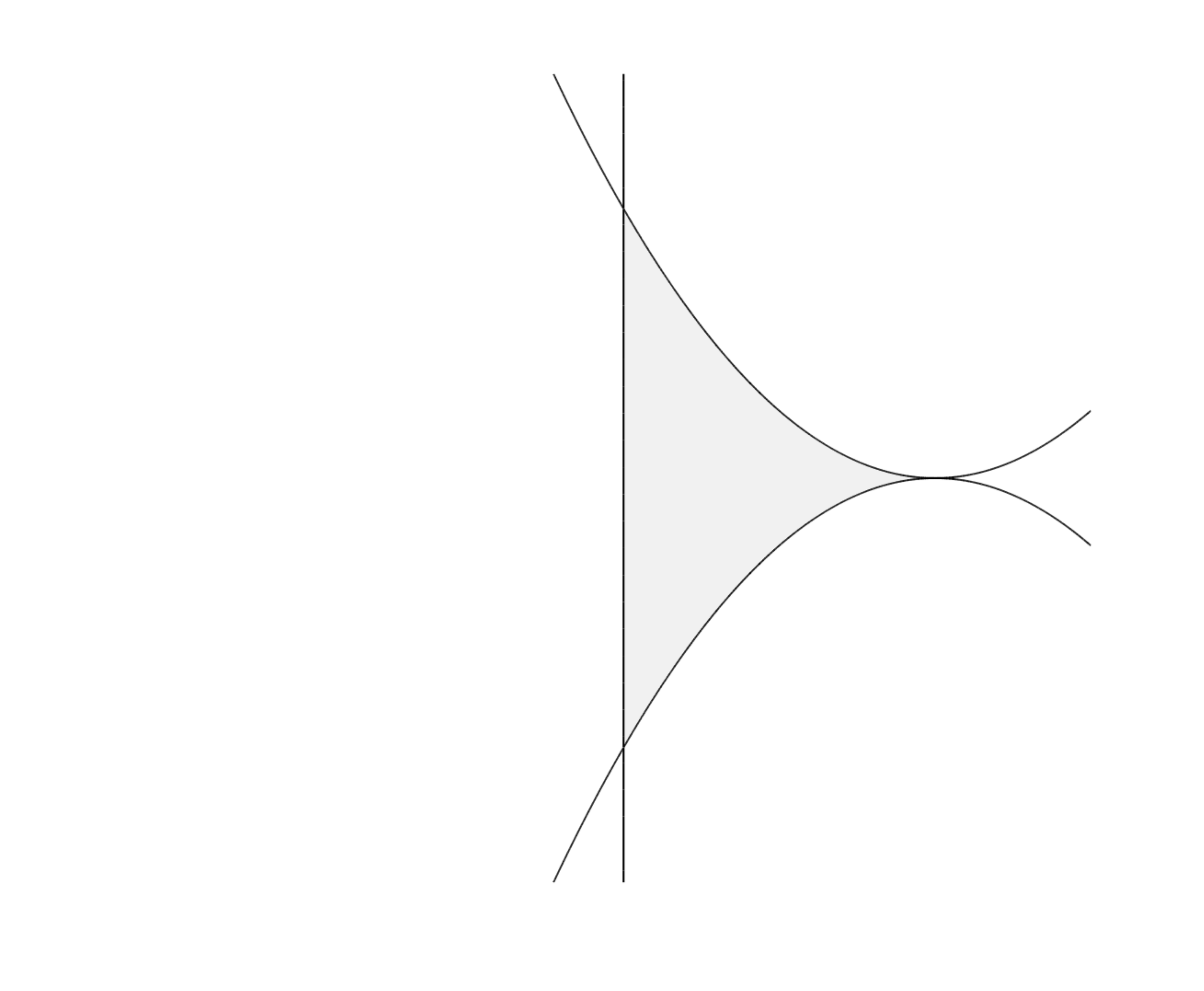}
\caption{The domain $\Omega_3^{(4)}$}
\end{center}
\end{figure}

Let us now add the secondary invariant $\eta= Y_n$, to treat the group $\cC_{n,J}$ when $n$ is even and  $\cC_n|\cC_{2n}$ when $n$ is odd. 
We get a new co-metric 
$$\Gamma_4= \bpm 4 \theta_1(1- \theta_1) & -2n\theta_1 \theta_2& -2n\theta_1\eta\\
&n^2((1-\theta_1)^{n-1}-\theta_2^2)&-n^2\theta_2\eta \\
&& n^2((1-\theta_1)^{n-1}-\eta^2)
\epm
$$
The determinant of this matrix factorizes as 
$$n^4(1-\theta_1)^{n-1}\theta_1((1-\theta_1)^n-\theta_2^2-\eta^2). $$ The factor $P_1(\theta_1,\theta_2, \eta)=(1-\theta_1)^n-\theta_2^2-\eta^2$ represents the relation between $\eta, \theta_1, \theta_2$ (the syzygy).   
The two factors 
$\theta_1, P_1(\theta_1,\theta_2, \eta)=(1-\theta_1)^n-\theta_2^2-\eta^2$ satisfy the boundary equations
$$\Gamma_4(\theta_1, \log(\theta_1))= 4(1-\theta_1), \Gamma_4(\theta_2, \log(\theta_1))= -2n\theta_2,\Gamma_4(\eta ,\log(\theta_1))= -2n\eta$$
$$\Gamma_4(\theta_1, \log(P_1))= -4n\theta_1, \Gamma_4(\theta_2, \log(P_1))=-2n^2\theta_2, \Gamma_4(\eta, \log(P_1))= -2n^2\eta.$$
(This is not true for the factor $1-\theta_1$).  The domain $\Omega_4 \subset \bbR^3$ defined by 
$\theta_1\in (0,1), P_1(\theta_1, \theta_2, \eta)>0$, which has reduced  boundary equation $\theta_1P_1(\theta_1, \theta_2, \eta)=0$, is therefore such that 
$(\Omega_4, \Gamma_4)$ is a polynomial system. Observe that the syzygy, which appears in one component of the boundary of $\Omega_4$, may be written as $P(\theta_1, \theta_2)-\eta^2$, where $P$ is one of the components of the boundary of the corresponding 2-d domain $\Omega_3$.

\begin{figure}[!h]
\begin{center}
\includegraphics[width=12cm]{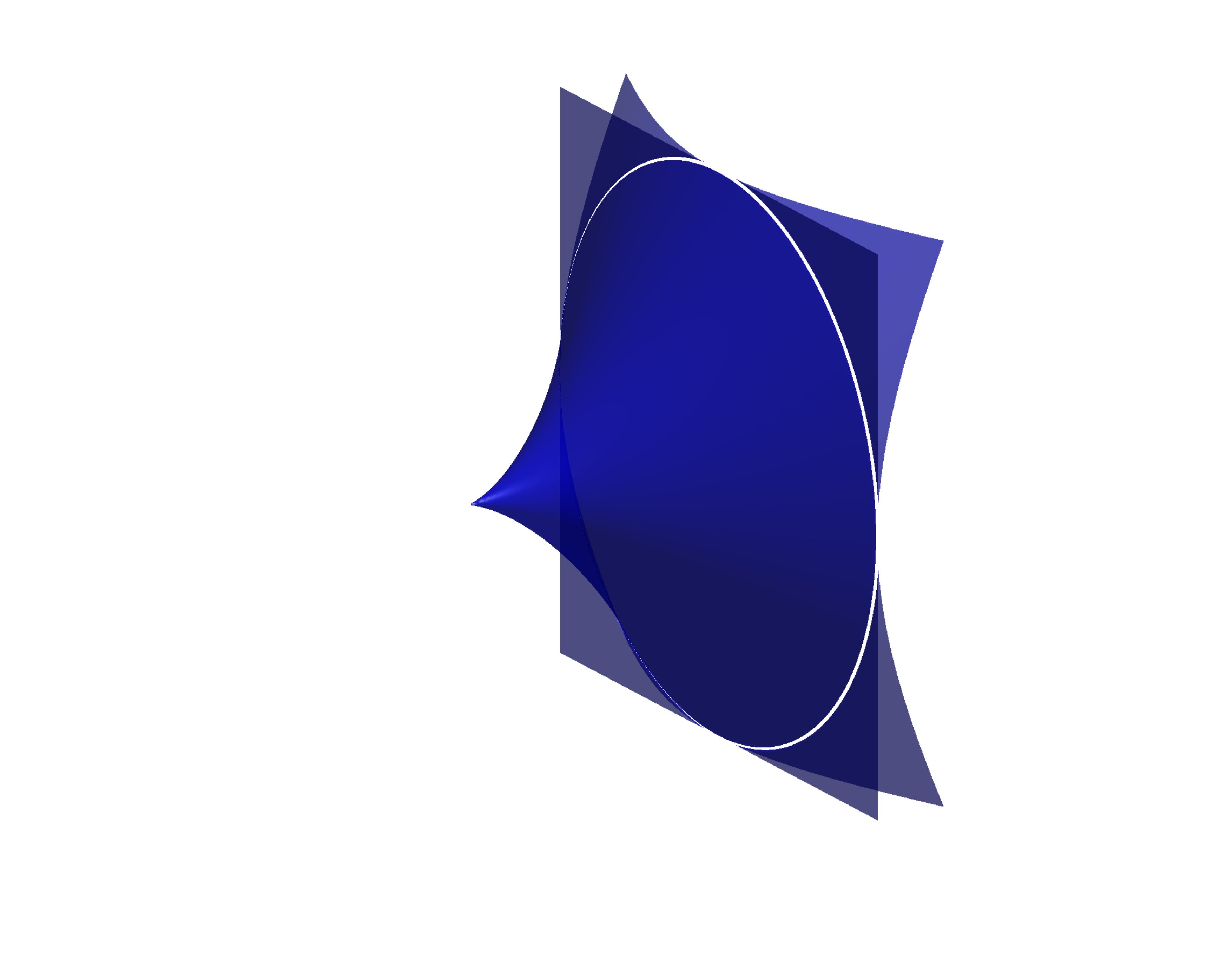}
\caption{The domain $\Omega_4^{(3)}$}
\end{center}
\end{figure}

We may now consider the dihedral group $\cD_n$, which amounts to add to the symmetries of $\cC_n$ the transformation $(x,y,z)\mapsto (x, -y, -z)$. 
It has as primary invariants $\theta_1= z^2, \theta_2= X_n$ as before (corresponding to the 2-d polynomial system $(\Omega_3, \Gamma_3)$), but now the secondary invariant is $\eta= zY_n$. The new co-metric in dimension 3 is then
$$\Gamma_5=\bpm
 4 \theta_1(1- \theta_1) & -2n\theta_1 \theta_2& -2\eta((n+1)\theta_1-1)\\
&n^2((1-\theta_1)^{n-1}-\theta_2^2)&-n(n+1)\theta_2\eta \\
&& (1-\theta_1)^{n-1}(1+(n^2-1)\theta_1)-\theta_2^2-(n+1)^2\eta^2
\epm
$$
The determinant of this metric factorizes as 
$$4n^2(\theta_1(1-\theta_1)^{n}-\theta_1\theta_2^2-\eta^2)((1-\theta_1)^{n-1}((n^2-1)\theta_1-1)-\theta_2^2)= 4n^2P_1P_2,$$ where 
$P_1(\theta_1, \theta_2, \eta)= \theta_1(1-\theta_1)^{n}-\theta_1\theta_2^2-\eta^2$ is the syzygy which relates $\eta$ to $(\theta_1, \theta_2)$. Observe once again the relation between this syzygy and the boundary equation of the corresponding  2-d domain $\Omega_3$.

Once again, we have 
$$\Gamma_5(\theta_1, \log (P_1))= 4(1-(n+1)\theta_1), \Gamma_5(\theta_2, \log (P1))= -2n(n+1)\theta_2, \Gamma_5(\eta, \log(P1))= -2(n+1)^2\eta$$
while the boundary equation is not satisfied for $P_2$. 
In $\bbR^3$, the domain $\Omega_5$ delimited by $\theta_1\in (0,1)$, $P_1>0$ is a bounded domain with reduced boundary equation $P_1=0$, and
$(\Omega_5, \Gamma_5)$ provides a 3-dimensional polynomial system. 

\begin{figure}[!h]
\begin{center}
\includegraphics[width=12cm]{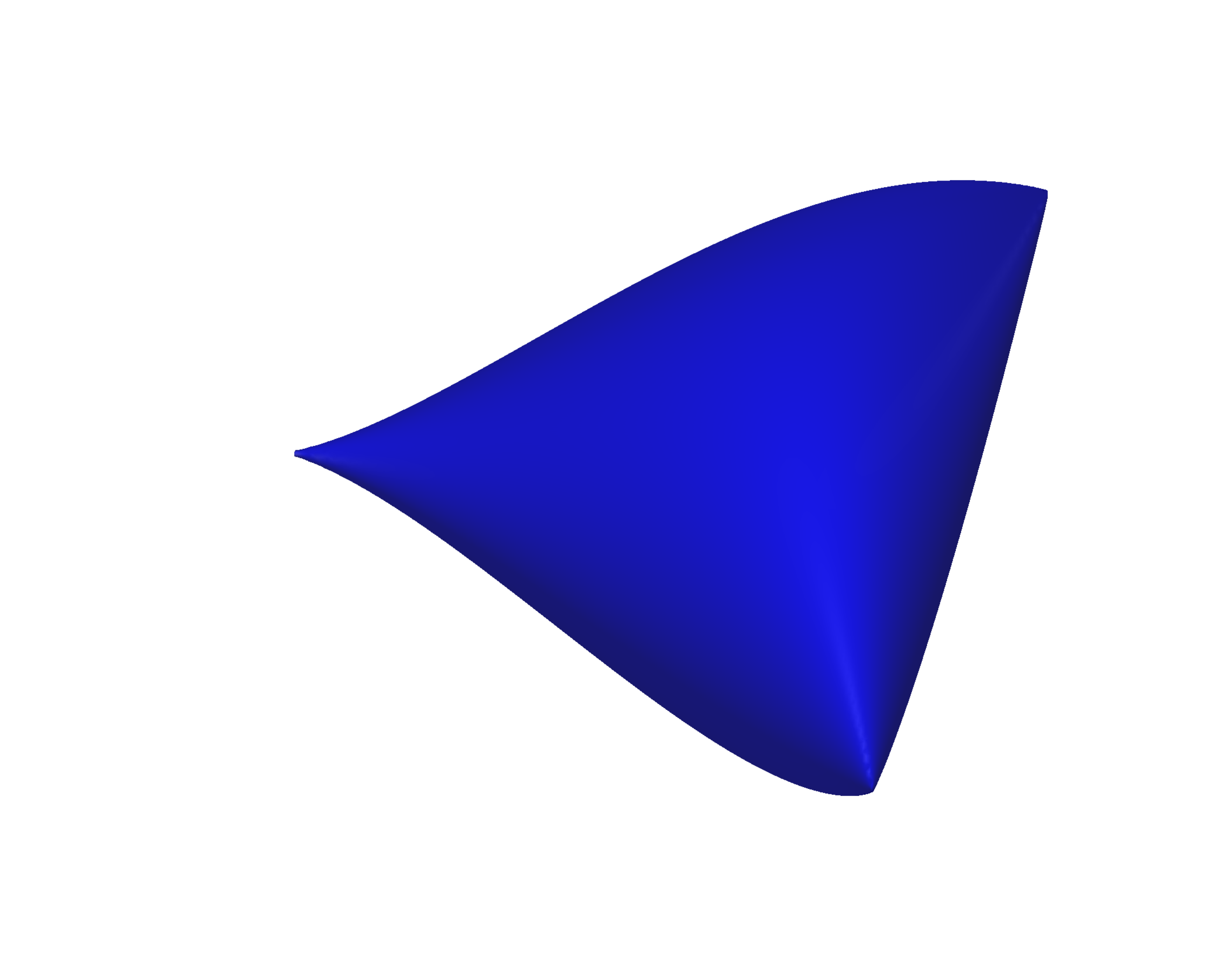}
\caption{The domain $\Omega_5$ for $n=3$}
\end{center}
\end{figure}

The groups $\cD_n|\cD_{2n}$ for $n$ even or $\cD_{n,J}$ for $n$ even have primary invariants $(z^2= \theta_1, X_n^2= \theta_2)$ and secondary invariant $\eta_1= zY_n$ or $\eta_2= zX_n$. However the groups $\cC_n|\cC_{2n}$ ($n$ even) or $\cC_{n,J}$ ($n$ odd) have the same primary invariants and as secondary invariants $zX_n,zY_n, X_nY_n$. It may be worth to observe that $(z, X_n^2, X_nY_n)$ is another form of the invariants for $\cC_{2n}$, since $X_{2n}= X_n^2-Y_n^2= 2X_n^2-(1-z^2)^n$, and $Y_{2n} =2X_nY_n$.  and therefore they do not provide any new model (although they provide them under another form).

We first choose $\theta_1= z^2, \theta_2= X_n^2$ (for which we already know that it corresponds to $(z^2, X_{2n})$ through a change of variables. We then get a co-metric 
$$\Gamma_6= 4 \bpm \theta_1(1-\theta_1) & -n\theta_1\theta_2\\
& n^2\theta_2((1-\theta_1)^{n-1}-\theta_2)
\epm,
$$
which corresponds to a 2-d domain $\Omega_6$ with boundary reduced equation  $\theta_1\theta_2( (1-\theta_1)^n-\theta_2)=0$, which is isomorphic  to the domain $\Omega_3$ when changing $n$ into $2n$ (this model has 3 irreducible components in its boundary equation, and two  of them may be reduced to a line, providing then a simpler form).

\begin{figure}[!h]
\begin{center}
\includegraphics[width=8cm]{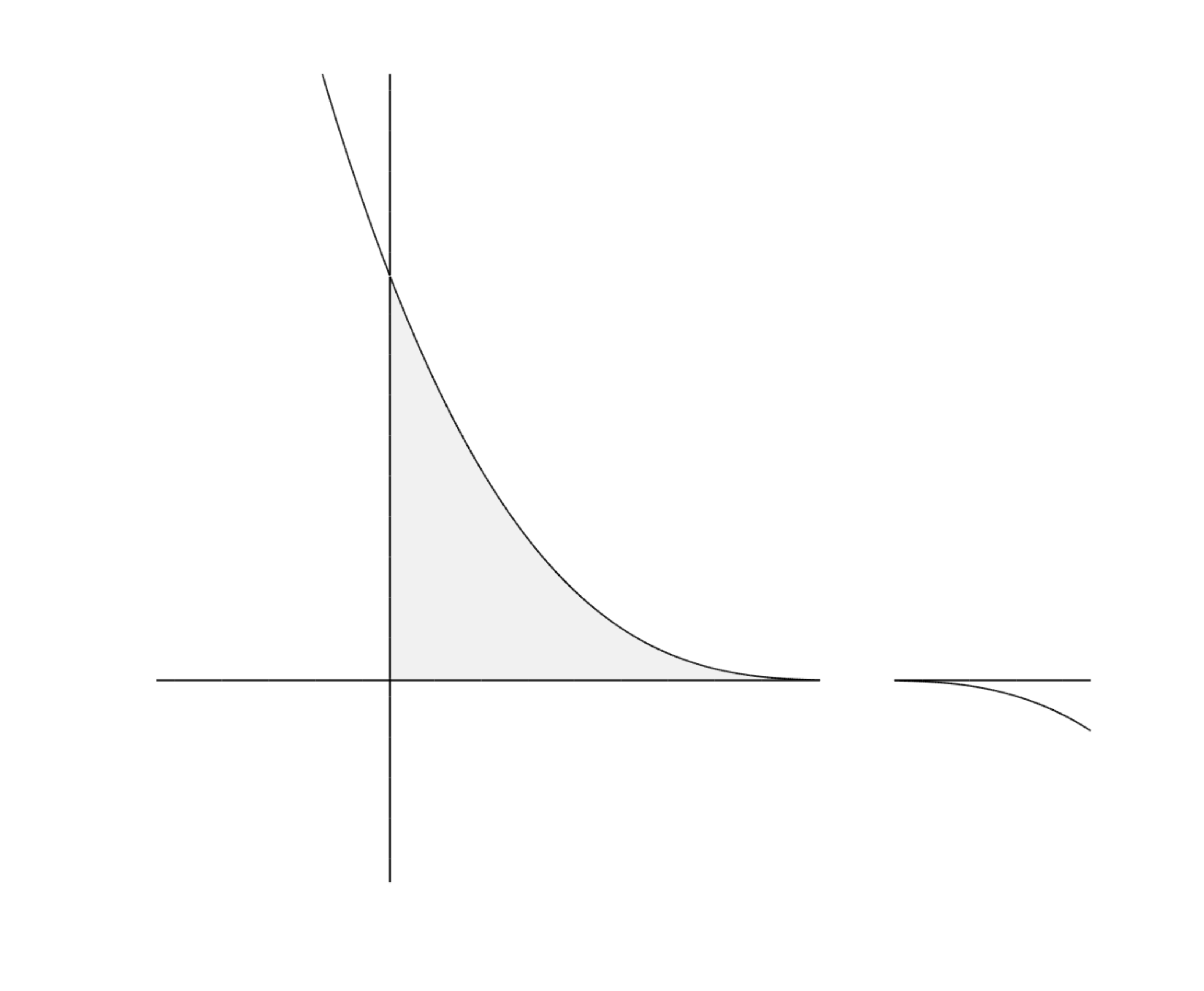}
\caption{The domain $\Omega_6$ for $n=3$}
\end{center}
\end{figure}

We may first add the secondary invariant $\eta_1= zY$. We get a  co-metric
$$\Gamma_7=\bpm 4\theta_1(1-\theta_1)&-4n\theta_1\theta_2&2\eta_1(1-(n+1)\theta_1)\\
&4n^2\theta_2((1-\theta_1)^{n-1}\theta_1-\theta_2)&-2n(n+1)\theta_2\eta_1\\
&& (1-\theta_1)^{n-1}(1+(n^2-1)\theta_1)-(n+1)^2\eta_1-\theta_2
\epm
$$
 The  syzygy relation between $(\theta_1,\theta_2, \eta)$ may be written as 
 $$P(\theta_1, \theta_2, \eta_1)=\theta_1\big((1-\theta_1)^n -\theta_2\big)-\eta_1^2,$$ once again of the form $Q(\theta_1, \theta_2)-\eta^2$, where $Q$ appears in the boundary equation of the corresponding 2-d domain $\Omega_6$.
 
 One may check that for this co-metric, $\theta_2P$ divides $\det(M)$, and moreover that 
 $$\Gamma_7(\theta_1, \log (P))=4(1-(n+1)\theta_1), ~\Gamma_7(\theta_2, \log (P))=-4n(n+1)\theta_2, ~\Gamma_7(\eta_1, \log (P))= -2(n+1)^2\eta_1.$$
 and also
 $$\Gamma_7(\theta_1, \log(\theta_2))=-4n\theta_1, ~\Gamma_7(\theta_2, \log(\theta_2))= 4n^2((1-\theta_1)^{n-1}-\theta_2, \Gamma_7(\eta, \log(\theta_2))=-2n(n+1)\eta_1.$$
 
 This provides a 3-d domain $\Omega_7\subset \bbR^3$ with boundary reduced equation $\theta_2P(\theta_1, \theta_2, \eta_1)=0$, such that 
 $(\Omega_7, \Gamma_7)$ is again a polynomial system.

\begin{figure}[!h]
\begin{center}
\includegraphics[width=13cm]{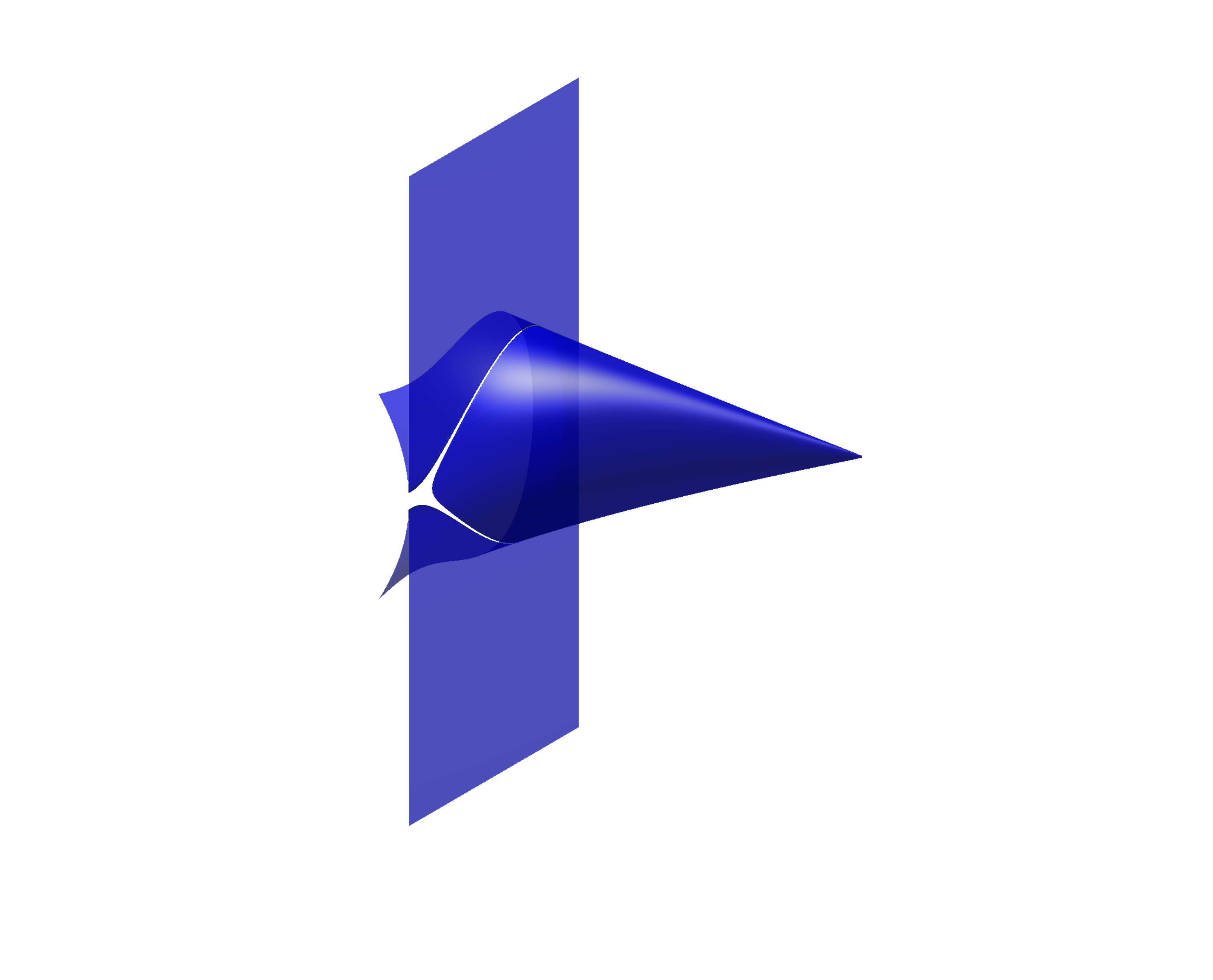}
\caption{The domain $\Omega_7$ for $n=3$}
\end{center}
\end{figure}

Adding the variable  
 $\eta_2= zX$, instead of $\eta_1= zY$, to $\theta_1= z^2, \theta_2=X^2$ leads to the co-metric

 $$\Gamma_8=\bpm 4\theta_1(1-\theta_1)&-4n\theta_1\theta_2&2\eta_2(1-(n+1)\theta_1)\\
&4n^2\theta_2((1-\theta_1)^{n-1}-\theta_2)&2n\eta_2(n(1-\theta_1)^{n-1}-(n+1)\theta_2)\\
&&n^2\theta_1 (1-\theta_1)^{n-1}+ \theta_2-(n+1)^2\eta_2^2\epm
$$
The determinant of this matrix has 3 factors, 2 of them being $P_1=\theta_1\theta_2-\eta_2^2$ and $P_2=\theta_2-(1-\theta_1)^{n}$. $P_1(\theta_1, \theta_2, \eta_2)$ is the syzygy relating $\eta_2$ to $(\theta_1, \theta_2)$. It is still of the form $Q(\theta_1, \theta_2)-\eta_2^2$, where $Q$ appears in the boundary equation of the corresponding 2-d domain $\Omega_6$.

Now, once again, we have
$$\bcas  \Gamma_8(\theta_1, \log(P_1))=4(1-(n+1)\theta_1),\\ \Gamma_8(\theta_2, \log(P_1))=4n((1-\theta_1)^{n-1}-(n+1)\theta_2), \\ 
\Gamma_8(\eta_2, \log(P_1))=-2(n+1)^2\eta_2,\ecas
$$
and
$$\bcas \Gamma_8(\theta_1, \log(P_2))=-4n\theta_1,\\\Gamma_8(\theta_2, \log(P_2))=-4n^2\theta_2, \\\Gamma_8(\eta_2, \log(P_2))= -2n(n+1)\eta_2.\ecas$$
 The third factor  of the determinant does not satisfy the boundary equation.  The domain $\Omega_8$ with boundary reduced equation $P_1P_2=0$ provides a polynomial system $(\Omega_8, \Gamma_8)$.

\begin{figure}[!h]
\begin{center}
\includegraphics[width=6cm]{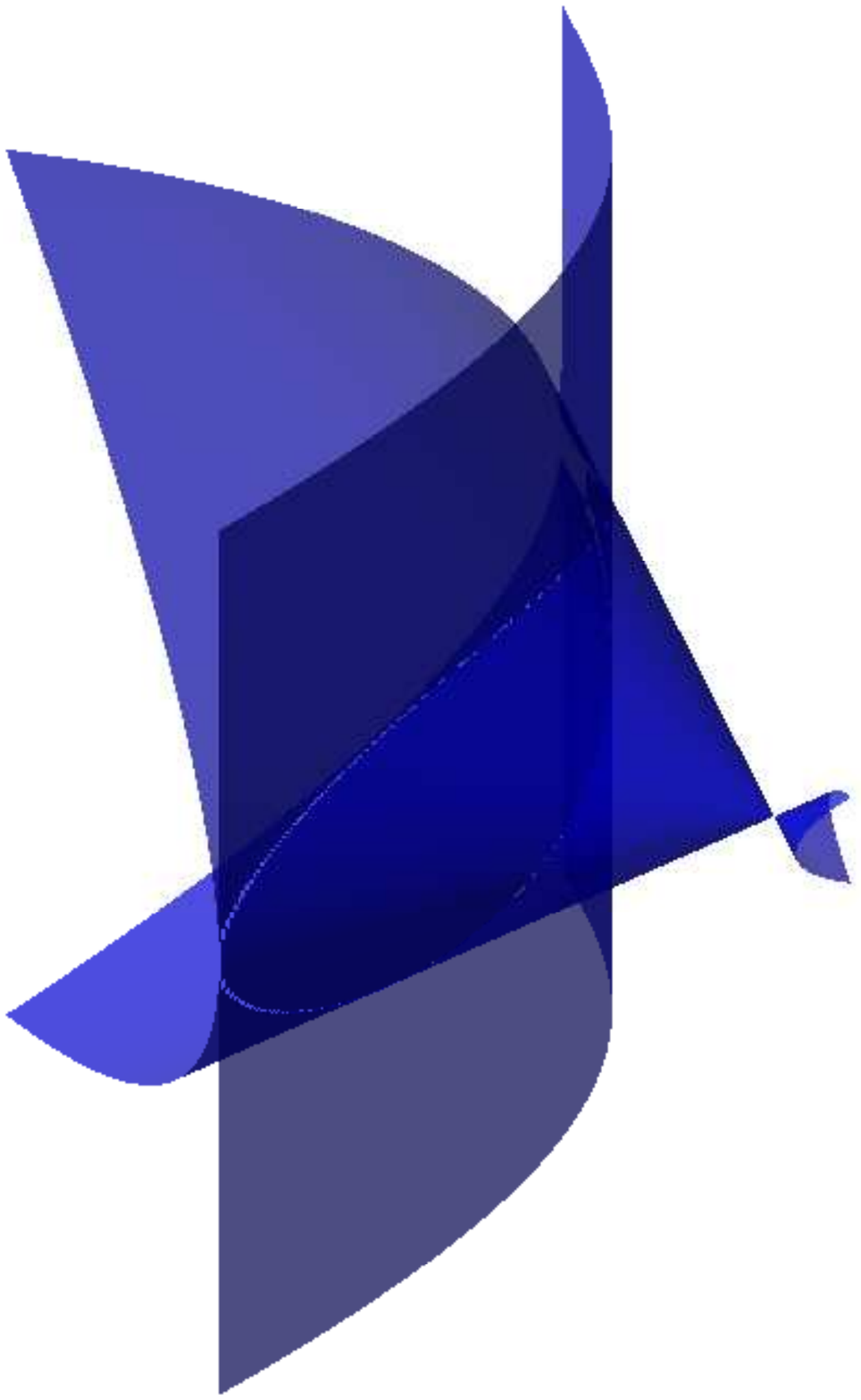}
\caption{The domain $\Omega_8$ for $n=3$}
\end{center}
\end{figure}

 We now may add instead the secondary invariant $\eta_3= X_nY_n$. As already observed, the new system $(z^2= \theta_1, X_n^2= \theta_2, \eta_3= X_nY_n)$ is isomorphic  to the system 
 $(\theta_1= z^2, \theta_2= X_{2n}, \eta= Y_{2n})$ described by the metric $\Gamma_4$ and the domain $\Omega_4$. Observe however that in this presentation, the reduced boundary of the domain is $\theta_1(\theta_2(1-\theta_1)^{n}- \theta_2^2-\eta_3^2)=0$, the second factor being the syzygy.

\begin{figure}[!h]
\begin{center}
\includegraphics[width=10cm]{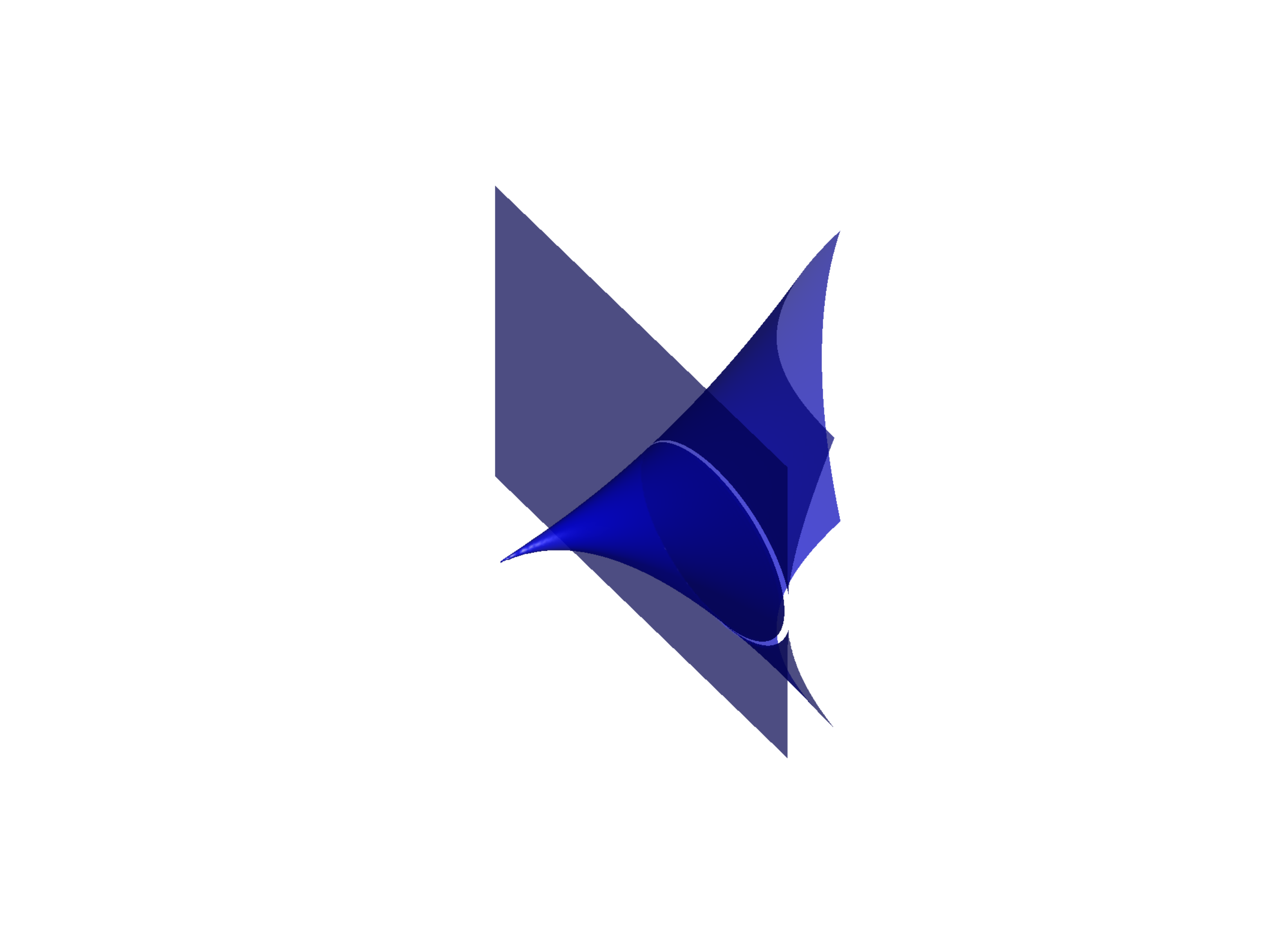}
\caption{The domain $\Omega_9$ for $n=3$. }
\end{center}
\end{figure}

 Finally, one may check that adding two of the secondary invariants $\eta_i$ to $\theta_1= z^2$, $\theta_2= X_n^2$ will not provide any  closed system for $\Gamma$. The system $(\theta_1, \theta_2, \eta_1, \eta_2, \eta_3)$ provides a closed $5$ dimensional  $\Gamma$ operator, but the determinant of the metric vanishes (in $\bbR^5$) and there does not seem to be any polynomial system associated with it.

\section{ Tetrahedron and cube/octahedron\label{sec.cube}}

It makes sense, as we shall see, to treat jointly these cases.  The groups $\cT, \cO$ correspond  to the elements of $SO(3)$ preserving respectively the tetrahedron, the octahedron or its dual the cube. They have respective  cardinality $12$ and $24$. Adding  the central symmetry $J: x\mapsto -x$ to each of them we obtain $\cT_J$ and $\cO_J$. Observe that the first one does not preserve the tetrahedron, while the second one does preserves the cube. We also consider the group $\cT | \cO$ which can be obtained by adding a plane symmetry with respect to the plane symmetry axes of the tetrahedron and which preserves the tetrahedron regardless of orientation.

They are related by the following inclusions diagram : 
$$ 
\begin{array}{rcccl}
&&\cT&&\\
&\swarrow & \downarrow & \searrow & \\ 
\cT_J && \cO && \underline{\cT|\cO} \\
&\searrow & \downarrow & \swarrow & \\
&&\underline{\cO_J}&&
\end{array}
$$

Let $(x,y,z)$ be the standard coordinate system in $\bbR^3$. We put the cube centered at the origin and with faces parallel to the coordinate planes. We put the tetrahedron with edges on the diagonal of the cube. We consider the polynomials $$O_3 = x y z, O_4 = x^4 + y^4 + z^4, O_6= (x^2-y^2)(y^2-z^2)(z^2-x^2),$$ which will play  the same rôle 	as the one played by $(z, X_nY_n)$ in the previous one as basic blocks to construct all the invariants for the various groups concerned in this section.

 We first compute :  
  $$\Gamma = \bpm \Gamma(O_3,O_3)& \Gamma(O_3,O_4)&\Gamma(O_3,O_6)\\ & \Gamma(O_4, O_4)& \Gamma(O_4, O_6)\\ & & \Gamma(O_6,O_6)\epm$$
 which is 
 $$\bpm(1-O_4)/2-9O_3^2&4O_3(1-3O_4)&-18 O_3O_6 \\&8(6O_3^2+3O_4-1-2O_4^2)&8O_6(2-3O_4)\\&&-54 O_3^2O_4+18 O_3^2-3O_4^2+4O_4-1-36 O_6^2\epm.$$

We consider the primary invariants of the Coxeter group $\cT|\cO$ given by $(\theta_1, \theta_2) = (O_3, O_4)$. The determinant of the submatrix $\Gamma_{11}$ given by the first two rows and columns  is 
  $$P(\theta_1,\theta_2) = -108 \theta_1^4 + 20 \theta_1^2 + 2 \theta_2^3- 5 \theta_2^2 +4 \theta_2 - 36 \theta_1^2 \theta_2.$$ 
It provides a domain $\Omega_{11}$ with boundary $P(X,Y)=0$. This corresponds to the  model of the {\it swallow tail} (example~10 in Section~\ref{sec.BOZmodels}).

\begin{figure}[!h]
\begin{center}
\includegraphics[width=10cm]{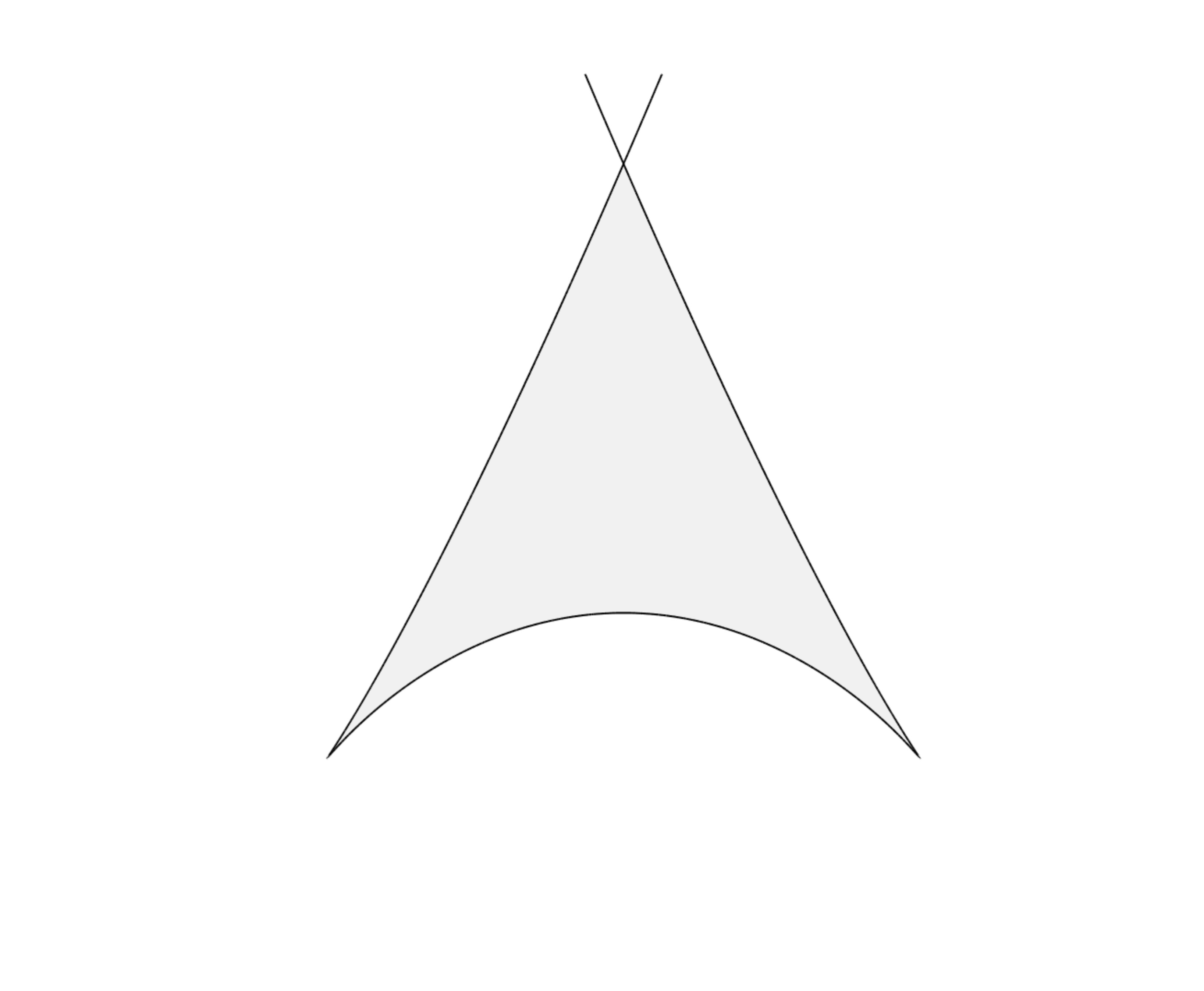}
\caption{The domain $\Omega_{11}$ : swallow tail.}
\end{center}
\end{figure}

\paragraph{For the group $\cT$,} we can choose $\eta = O_6$ as secondary invariant. It is algebraically related to $(\theta_1,\theta_2)  =( O_3, O_4)$ through 
 $$\eta^2 = P(\theta_1, \theta_2).$$
 If we write the matrix 
 $$\Gamma_{12}= \bpm \Gamma(\theta_1,\theta_1)&\Gamma(\theta_1,\theta_2)&\Gamma(\theta_1,\eta)\\
& \Gamma(\theta_2,\theta_2)& \Gamma(\theta_2,\eta)\\
 & & \Gamma(\eta,\eta)\epm,$$ we get 
$$ \left( \begin {array}{ccc} -9\,{\theta_1}^{2}-\theta_2/2+1/2&-12\,\theta_1\theta_2+4\,\theta_1&-18\,\theta_1\eta
\\ \noalign{\medskip}&48\,{\theta_1^2}-16\,{\theta_2}^{2}+24\,\theta_2-8&-24
\,\theta_2\eta+16\,\eta\\ \noalign{\medskip}&&-54\,{\theta_1^2}\theta_2+18
\,{\theta_1^2}-3\,{\theta_2}^{2}-36\,{\eta}^{2}+4\,\theta_2-1\end {array} \right).  $$
The determinant of this matrix factorizes in 
$$D=4\, \left( 3\,\theta_2-1 \right)  \left( 18\,{\theta_1}^{2}+\theta_2-1 \right)  \left( 108
\,{\theta_1}^{4}+36\,{\theta_1}^{2}\theta_2-2\,{\theta_2}^{3}-20\,{\theta_1}^{2}+5\,{\theta_2}^{2}+4\,{\eta}^{2}-4
\,\theta_2+1 \right) $$
It turns out that the factor
$$P_3(\theta_1,\theta_2,\eta)=  108
\theta_1^{4}+36\,\theta_1^{2}\theta_2-2\,\theta_2^{3}-20\,\theta_1^{2}+5\,\theta_2^{2}+4\,\eta^{2}-4
\,\theta_2+1 $$ satisfies  
$$\bcas \Gamma_{12}(\theta_1,\log P_3)= -36 \theta_1, \\\Gamma_{12}(\theta_2,\log(P_3))=-48\theta_2+32, \\\Gamma_{12}(\eta, \log(P_3))=-72 \eta,\ecas $$ so that this provides a new polynomial model $(\Omega_{12}, \Gamma_{12})$ in dimension 3.
(The boundary equation is not satisfied for the two other factors.)

We can check that the complementary of the surface $P_3(X,Y,Z)=0$ has one bounded component in $\bbR^3$ and that the determinant does not vanish inside this component. 
\begin{figure}[!h]
\begin{center}
\includegraphics[width=12cm]{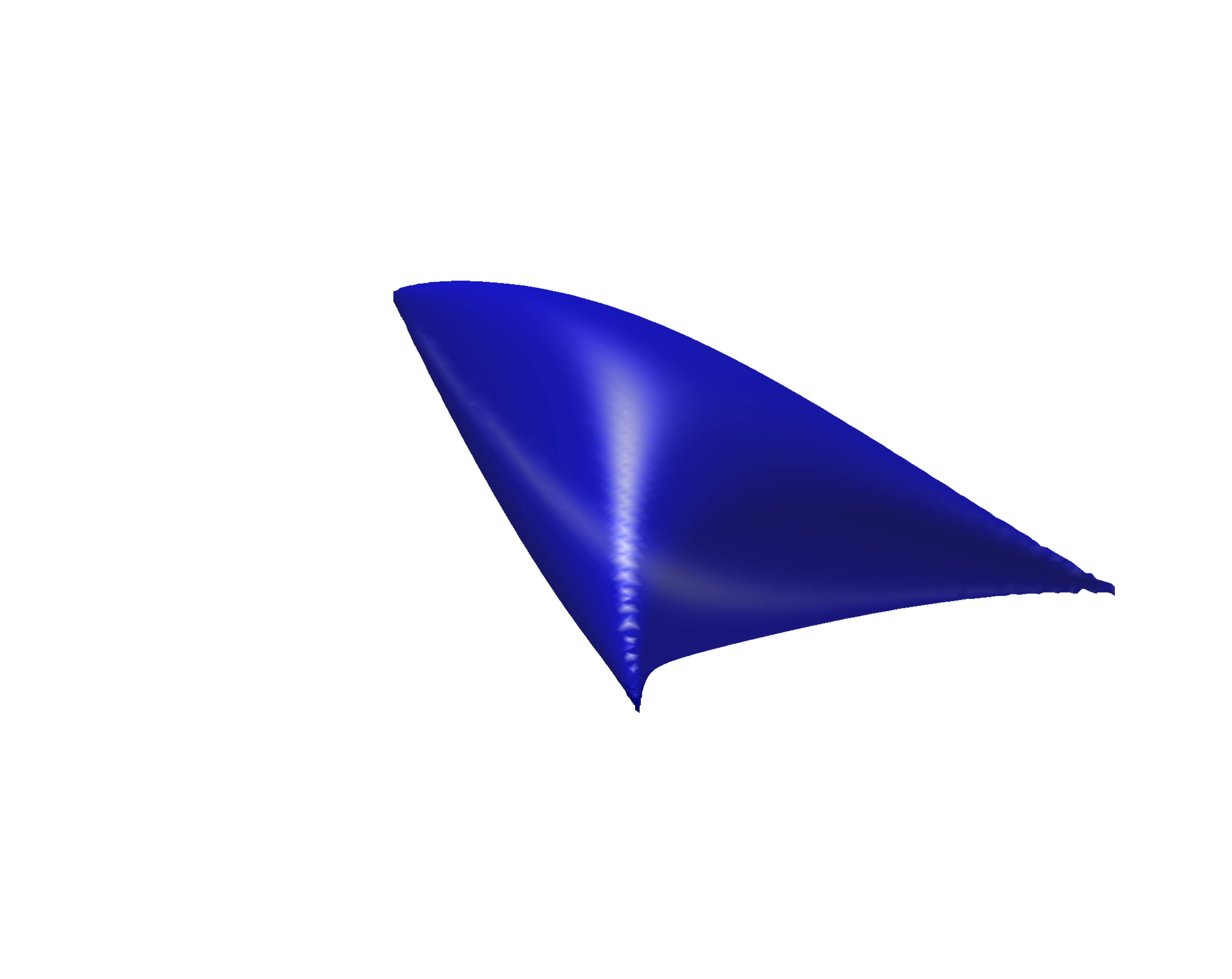}
\caption{Surface $P_3(X,Y,Z)=0$, bounding $\Omega_{12}$. }
\label{fig:tetra2}
\end{center}
\end{figure}

\paragraph{We observe that $\cO_J$} is also a Coxeter group. We can take as primary invariants  $(\theta_1,\theta_2) =( O_3^2, O_4)$ for $\cO_J$.  We get 
$$\Gamma_{13} = \bpm \Gamma(\theta_1,\theta_1)= 4\theta_1((1-\theta_2)/2-9\theta_1) & \Gamma(\theta_1,\theta_2)= 8\theta_1(1-3\theta_2)\\
& \Gamma(\theta_2,\theta_2)= 16(3\theta_1+\frac{3}{2}\theta_2-1/2-\theta_2^2)\epm$$
whose determinant is given by, 
$$Q(\theta_1,\theta_2) = 16 \theta_1 (-108 \theta_1^2 +20 \theta_1 +2 \theta_2^3 - 5 \theta_2^2 + 4 \theta_2 -1 -36 \theta_1 \theta_2) $$
Observe that  $Q(X^2,Y) = X^2 P(X,Y)$.  We recognize that the boundary of the domain $\Omega_{13}$ is the {\it cuspidal cubic with tangent} (model~9 in Section~\ref{sec.BOZmodels}). 
 \begin{figure}[!h]
\begin{center}
\includegraphics[width=8cm]{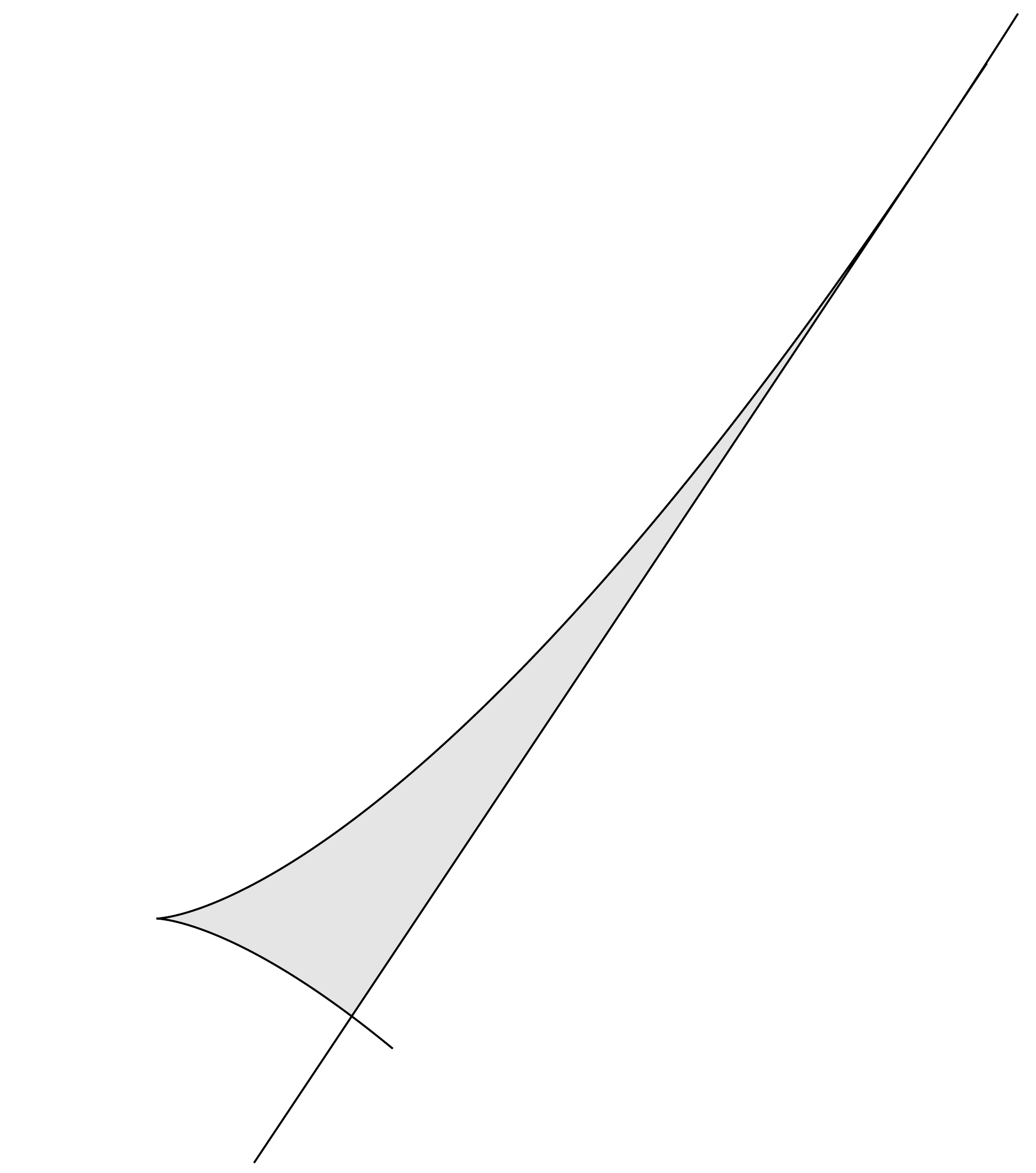}
\caption{The domain $\Omega_{13}$ : cuspidal cubic with tangent.}
\end{center}
\end{figure}

\paragraph{For the group $\cT_J$,} we may add the secondary invariant $ \eta= O_6$. We obtain the co-metric $\Gamma_{14}:$ 
$$ \left[ \begin {array}{ccc} -36\,{\theta_1}^2-2 \theta_1\theta_2+2\theta_1&-24\,\theta_1\theta_2+8\,\theta_1&-36 \theta_1\eta\\ \noalign{\medskip}&48\,{\theta_1}-16\,{\theta_2}^{2}+24\,\theta_2-8&-24
\,\theta_2\eta+16\,\eta\\ \noalign{\medskip}  &&-54\,{\theta_1}\theta_2+18
\,{\theta_1}-3\,{\theta_2}^{2}-36\,{\eta}^{2}+4\,\theta_2-1\end {array} \right] .$$
The determinant of this matrix factorizes as 
$$16\theta_1(3\theta_2-1)(18\theta_1+\theta_2-1)(-2\theta_2^3+108\theta_1^2+36\theta_1\theta_2+5\theta_2^2+4\eta^2-20\theta_1-4\theta_2+1).
$$
Only the two factors $\theta_1$ and 
$$Q_3(\theta_1, \theta_2, \eta)= (-2\theta_2^3+108\theta_1^2+36\theta_1\theta_2+5\theta_2^2+4\eta^2-20\theta_1-4\theta_2+1)$$ satisfy the boundary equation, with
$$\bcas \Gamma_{14}(\theta_1, \log (\theta_1))=-36 \theta_2-2\theta_2+2, \\
\Gamma_{14}(\theta_2, \log (\theta_1))= -24\theta_2+8,\\
\Gamma_{14}(\eta, \log (\theta_1))=-36 \eta, 
\ecas
$$ and
$$\bcas \Gamma_{14}(\theta_1, \log (Q_3))=-72 \theta_1, \\
\Gamma_{14}(\theta_2, \log (Q_3))=-48 \theta_2+ 32,\\
\Gamma_{14}(\theta_3, \log (Q_3))=-73 \theta_3.
\ecas
$$

We observe that this factor writes $\eta^2 -  Q(\theta_1,\theta_2)$.  It appears that the two components bound a domain $\Omega_{14}$ on which the other factors do not vanish.  
 \begin{figure}[!h]
\begin{center}
\includegraphics[width=12cm]{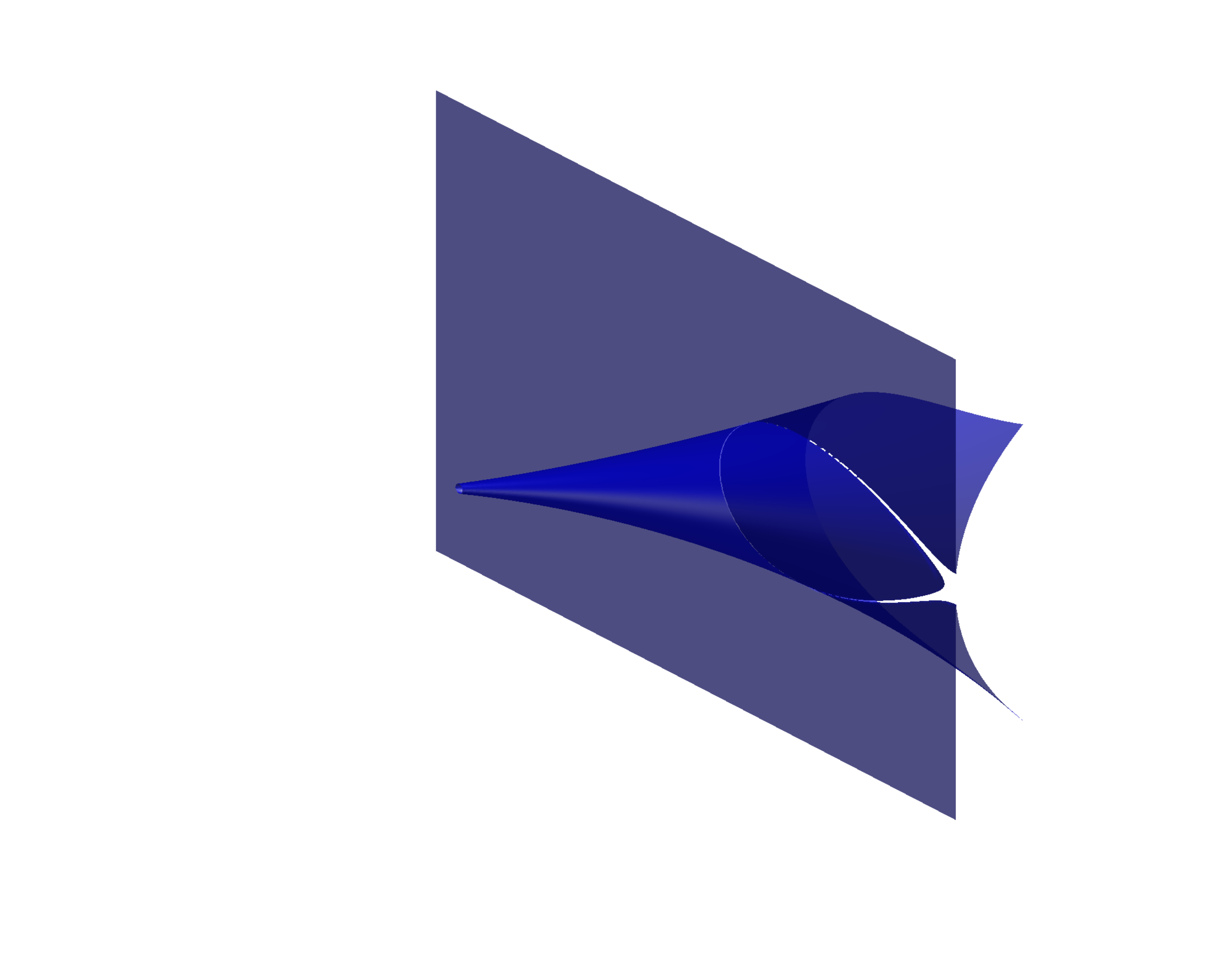}
\caption{Surface $XQ_3(X,Y,Z)=0$ bounding the domain $\Omega_{14}$ }
\label{fig:tetrasym}
\end{center}
\end{figure}

\paragraph{Finally, for the group $\cO$,} we use $(\theta_1,\theta_2,\eta) = (O_3^2, O_4, O_3 O_6)$. 
 We  compute the $\Gamma_{15}$ matrix 
 
 $$ \left[ \begin {array}{ccc} -36\,{\theta_1}^2-2 \theta_1\theta_2+2\theta_1&-24\,\theta_1\theta_2+8\,\theta_1&-54\theta_1\eta-\theta_2\eta+\eta
\\ \noalign{\medskip}&48\,{\theta_1}-16\,{\theta_2}^{2}+24\,\theta_2-8& -36 \eta\theta_2+20 \eta\\ 
\noalign{\medskip}  & & G_{3,3}\end {array} \right],  $$

 with 
 $$G_{3,3}:=-81\,\eta^{2}-{\frac {81\,{\theta_{{1}}}^{2}\theta_{{2}}}{2}}+9/2\,{\theta_{{1}}
}^{2}-3\,\theta_{{1}}\theta_{{2}}+3/2\,\theta_{{1}}+3/2\,\theta_{{1}}{\theta_{{2}}}^{2}-1/4\,{t
_{{2}}}^{4}+{\frac {7\,{\theta_{{2}}}^{3}}{8}}-{\frac {9\,{\theta_{{2}}}^{2}}{8}
}+5/8\,\theta_{{2}}-1/8. 
$$
The determinant of the matrix factorizes as $2Q_1Q_2$, with
$$Q_1=2\,{\theta_{{2}}}^{4}+324\,{\theta_{{1}}}^{2}\theta_{{2}}-12\,\theta_{{1}}{\theta_{{2}}}^{2}-7
\,{\theta_{{2}}}^{3}-36\,{\theta_{{1}}}^{2}+24\,\theta_{{1}}\theta_{{2}}+9\,{\theta_{{2}}}^{2}-
12\,\theta_{{1}}-5\,\theta_{{2}}+1$$
$$Q_2= 2\,\theta_{{1}}{\theta_{{2}}}^{3}+108\,{\theta_{{1}}}^{3}+36\,{\theta_{{1}}}^{2}\theta_{{2}}+5
\,\theta_{{1}}{\theta_{{2}}}^{2}-20\,{\theta_{{1}}}^{2}-4\,\theta_{{1}}\theta_{{2}}+4\,\eta^{2}+\theta_{{1}}.
$$
Only  $Q_2$ satisfies the boundary equation, with
$$\bcas\Gamma_{15}(\theta_1, \log (Q_2))=2-108 \theta_1-2\theta_2\\
\Gamma_{15}(\theta_2, \log (Q_2))= 40-72 \theta_2\\
\Gamma_{15}(\eta, \log (Q_2))=-162 \eta
\ecas
$$
We observe that $(O_3O_6)^2 = O_3^3 Q(O_3^2, O_4)$ so that $\eta^2=\theta_1Q(\theta_1,\theta_2) =: R(\theta_1,\theta_2)$. 
We get then a new domain $\Omega_{15}$ with boundary
$$Q_2(\theta_1, \theta_2, \eta) := \eta^2-R(\theta_1,\theta_2) =0.$$
 \begin{figure}[!h]
\begin{center}
\includegraphics[width=15cm]{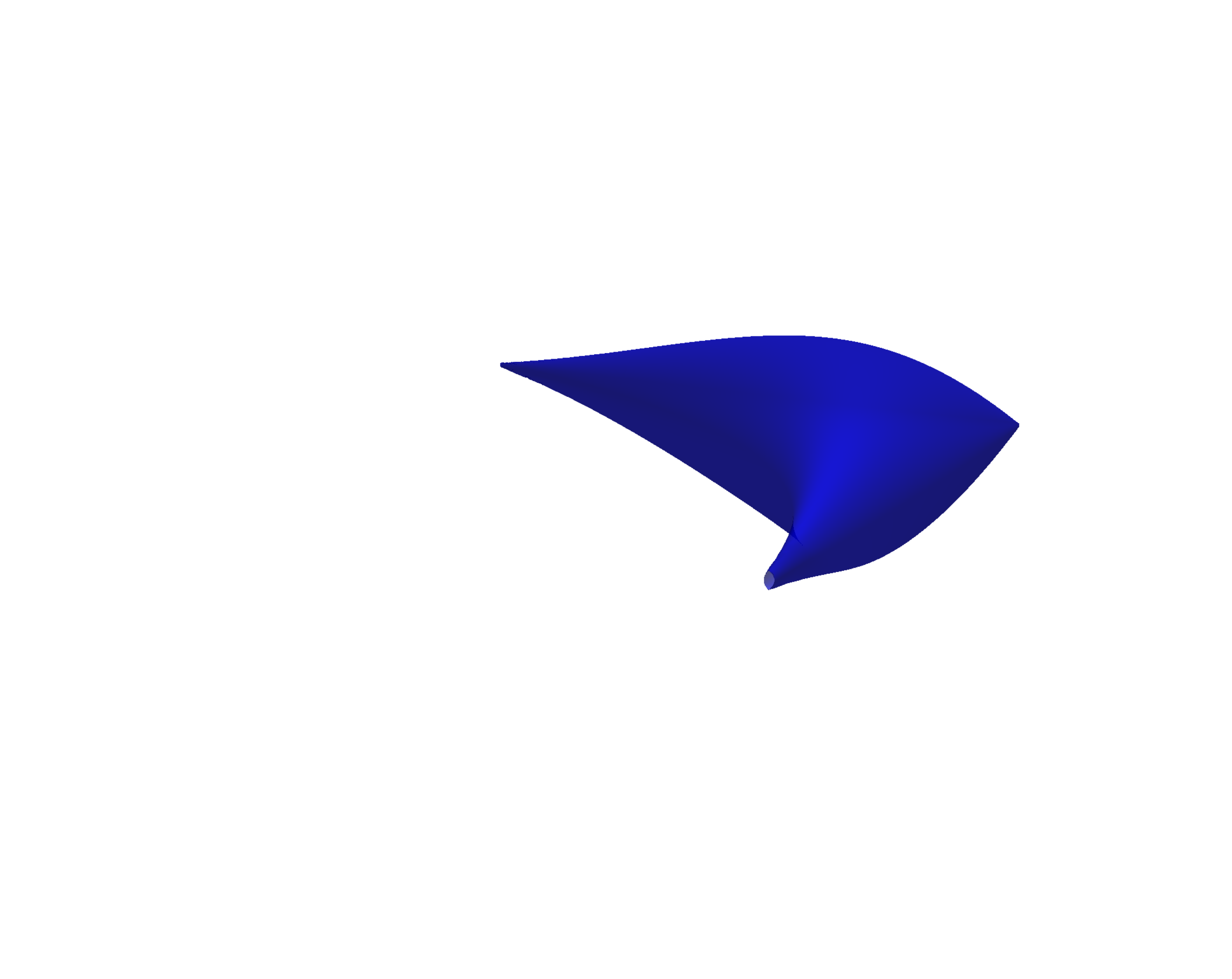}
\caption{Surface $R_3(X,Y,Z)=0$ bounding the domain $\Omega_{15}$. }
\label{fig:cube1}
\end{center}
\end{figure}
We may observe that in all these cases, the 3-dimensional domains can be cut by the plane $\{Z=0\}$ to get a 2-dimensional domain. In the first case, we get the swallow tail, in the other two cases, we get the cuspidal cubic with tangent. They provide various two fold coverings of these two dimensional models. It could be interesting to investigate the shape of the singularities of the boundaries of these 3-dimensional domains. 

\section{Dodecahedron / Icosahedron}
\label{sec:icosahedron}
We finish by the study of the groups $\cI$ and $\cI_J$ of cardinality $60$ and $120$. The computations are a bit more painful but the idea is always the same. Let us introduce as before  three new building blocks. 
With $c= (1+ \sqrt{5})/2$,  
$$\bcas I_6= (c^2x^2-y^2)(c^2y^2-z^2)(c^2z^2-x^2),\\
I_{10}= (x+y+z)(-x+y+z)(x-y+z)(x+y-z)\\
\qquad\qquad(c^{-2}x^2-c^2y^2)(c^{-2}y^2-c^2z^2))(c^{-2}z^2-c^2x^2),\\
I_{15}=xyz(cx+ c^{-1}y+z)(-cx+c^{-1}y+z)((cx- c^{-1}y+z)(cx+ c^{-1}y-z)\\
    \qquad\qquad(x+ cy+c^{-1}z)(-x+ cy+c^{-1}z)((x- cy+c^{-1}z)(x+ cy-c^{-1}z)\\
    \qquad\qquad (c^{-1}x+y+ cz)(-c^{-1}x+y+ cz)(c^{-1}x-y+ cz)(c^{-1}x+y- cz). 
    \ecas
$$
We will use the primary invariants $(\theta_1,\theta_2) = (I_6,I_{10})$ and the secondary invariant $\eta = I_{15}$. Let us do the computations directly with all the invariants in this family. With $(\theta_1,\theta_2,\eta)= (I_6, I_{10}, I_{15})$, we get for the co-metric 
$$\Gamma =\bpm \Gamma(\theta_1,\theta_1)& \Gamma(\theta_1,\theta_2)&\Gamma(\theta_1,\eta)\\& \Gamma(\theta_2,\theta_2)&\Gamma(\theta_2,\eta)\\&&\Gamma(\eta,\eta)\epm$$
with 
$$\bcas \Gamma(\theta_1,\theta_1)= -36\,{{ \theta_1}}^{2}- \left( \sqrt {5}+2 \right)  \left( 7\,{ \theta_1}+5
\,{ \theta_2}+2\,\sqrt {5}{ \theta_2} \right),\\
 \Gamma(\theta_1,\theta_2)= \left( 40-16\,\sqrt {5} \right) {{ \theta_1}}^{2}+ \left( 3\,\sqrt {5}+6
 \right) { \theta_2}+{ \theta_1}\,\sqrt {5}-60\,{ \theta_1}\,{ \theta_2},\\
\Gamma(\theta_2,\theta_2)=\left( 7296-3264\,\sqrt {5} \right) {{ \theta_1}}^{3}+ \left( 96\,\sqrt 
{5}-240 \right) { \theta_1}\,{ \theta_2}\\\qquad\qquad+ \left( -432+192\,\sqrt {5}
 \right) {{ \theta_1}}^{2}-5\,\sqrt {5}{ \theta_2}+ \left( 6-3\,\sqrt {5}
 \right) { \theta_1}-100\,{{ \theta_2}}^{2},\\
 \Gamma(\theta_1,\eta)=-90\,{ \theta_1}\,{ \eta}-2\, \left( \sqrt {5}+2 \right) { \eta}\\
 \Gamma(\theta_2,\eta)=-150\,{ \theta_2}\,{ \eta}+{ \eta}\, \left(  \left( -100+40\,\sqrt {5}
 \right) { \theta_1}-2\,\sqrt {5} \right) \\
 \Gamma(\eta,\eta)=-225\,{{ \eta}}^{2}-
 1/4\, \left( -161+72\,\sqrt {5} \right)  \\\qquad\qquad\left( 
13\,{ \theta_1}\,\sqrt {5}+45\,\sqrt {5}{ \theta_2}+45\,{{ \theta_1}}^{2}+ 4\,
\sqrt {5}+26\,{ \theta_1}+100\,{ \theta_2}+9 \right) \\ \qquad\qquad \left( 30\,{ \theta_1}\,{
 \theta_2}\,\sqrt {5}-3\,{ \theta_1}\,\sqrt {5}-9\,\sqrt {5}{ \theta_2}-19\,{{
 \theta_1}}^{2}+75\,{ \theta_1}\,{ \theta_2}-6\,{ \theta_1}-20\,{ \theta_2} \right) 

 \ecas
 $$

 For the Coxeter group $\cI_J$ we restrict our attention to $(\theta_1, \theta_2)$.  
 Up to some factor, the determinant of the sub-matrix $\Gamma_{21}$ corresponding to $(\theta_1, \theta_2)$ is $S(\theta_1, \theta_2)$, with 
 \beqnas 
 S(\theta_1,\theta_2)=&&688\,\sqrt {5} \theta_1^{4}+6480\,\sqrt {5}\theta_1^{3}{ \theta_2}+
1728\, \theta_1^{5}+364\, \theta_1^{3}\sqrt {5}+6042\,\sqrt {5}
 \theta_1^{2}{ \theta_2}+23400\,\sqrt {5}{ \theta_1}\,\theta_2^{2}
\\ &&+17050\,
\sqrt {5}\theta_2^{3}+1376\,\theta_1^{4}
+14400\, \theta_1^{3}{
 \theta_2}+68\,\theta_1^{2}\sqrt {5}+1288\,{ \theta_1}\,{ \theta_2}\,\sqrt {5
}+1220\,\sqrt {5}\theta_2^{2}\\&&-19520\,\sqrt {5}\eta^{2}
+819\,
 \theta_1^{3}+13515\, \theta_1^{2}{ \theta_2}+52325\,{ \theta_1}\,\theta_2^{2}+38125\, \theta_2^{3}+152\, \theta_1^{2}+2880\,{ \theta_1}\,{ \theta_2
}\\&&+2728\, \theta_2^{2}
\eeqnas
This is a new 2-dimensional domain $\Omega_{21}$.  
\begin{figure}[!h]
\begin{center}
\includegraphics[width=15cm]{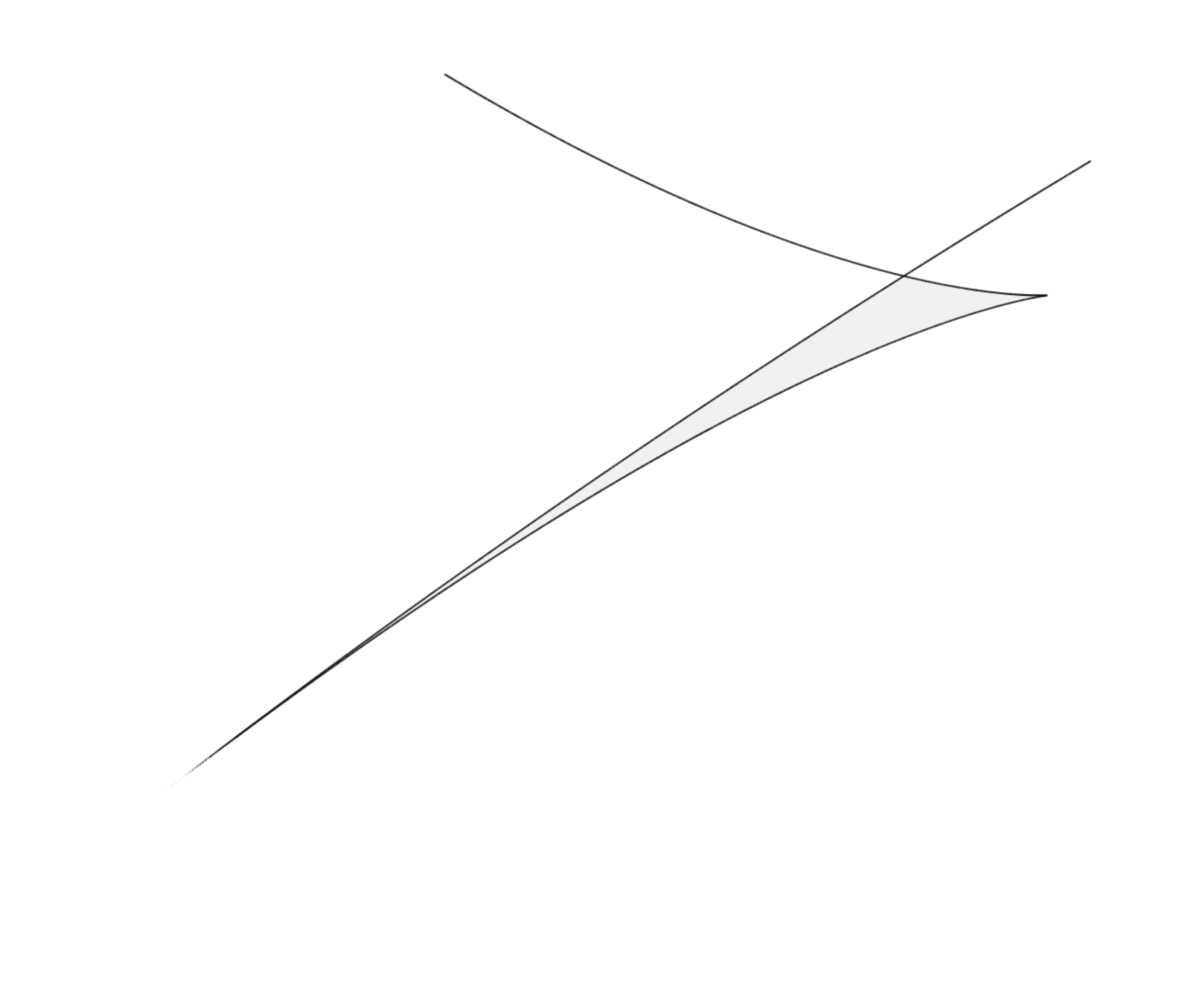}
\caption{Domain $\Omega_{21}$ bounded by  $S(X,Y)=0$.  }
\label{fig:modeleica}
\end{center}
\end{figure}

\paragraph{For the direct subgroup $\cI$,}  the determinant of $\Gamma_{22}$ factorizes (up to some constant) into $S_1S_2S_3$,
with
$$S_1(\theta_1,\theta_2,\eta)=S(\theta_1, \theta_2)-43648\,{{ \eta}}^{2}, $$
$$ S_2(\theta_1,\theta_2)= 13\,{ \theta_1}\,\sqrt {5}+45\,\sqrt {5}{ \theta_2}+45\,\theta_1^{2}+4\,\sqrt {5}+26\,{ \theta_1}+100\,{ \theta_2}+9, $$
and, 
$$S_3(\theta_1,\theta_2)=  30\,{ 
\theta_1}\,{ \theta_2}\,\sqrt {5}-3\,{ \theta_1}\,\sqrt {5}-9\,\sqrt {5}{ \theta_2}-
19\,\theta_1^{2}+75\,{ \theta_1}\,{ \theta_2}-6\,{ \theta_1}-20\,{ \theta_2}. $$
$S_1(\theta_1, \theta_2, \eta)$ is the syzygy relating $\eta$ to $\theta_1$ and $\theta_2$, and the 
 polynomial $S_1$ satisfies the boundary condition for the co-metric $\Gamma_{22}$,  which  once again  provides a new polynomial system with domain $\Omega_{22}$ in dimension 3, since the boundary conditions are satisfied:  
$$\bcas
\Gamma_{22}(\theta_1, \log (S_1))=-4\sqrt{5}-8 -180 \theta_1, \\
\Gamma_{22}(\theta_2, \log (S_1))=4(2\sqrt{5}-5)((30\sqrt{5}+75) \theta_2 + 10 \theta_1+ 2+ \sqrt{5}), \\
\Gamma_{22}(\eta, \log (S_1))=-450 \eta. 
\ecas
$$
Observe that we may as well rescale $\eta$ to get a simpler domain $\eta^2= S(\theta_1, \theta_2)$.
 \begin{figure}[!h]
\begin{center}
\includegraphics[width=12cm]{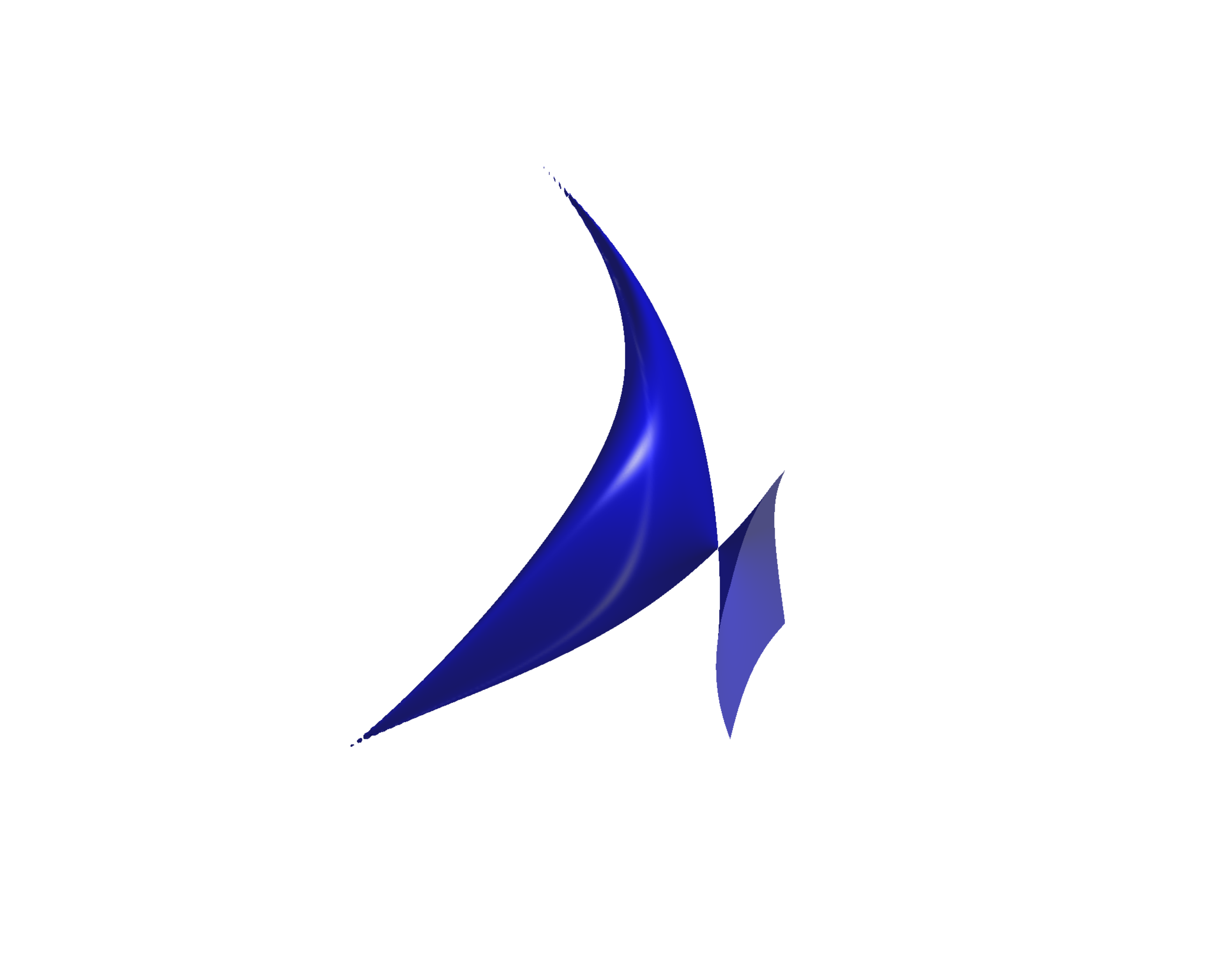}
\caption{Surface $S_1(X,Y,Z)=0$ bounding the domain $\Omega_{22}$.}
\label{fig:ico3D}
\end{center}
\end{figure}

\pagebreak

\section{Summary}
We summarize here  the models detailed above. For the dihedral family, we let $$H_n(X,Y) = (1-X)^n - Y.$$
$$
\begin{array}{|c||c|c|c|c|c|c|}
\hline \hbox{Group} & \theta_1 & \theta_2 & \eta & \Omega &\hbox{Boundary}& \hbox{Picture} \\
\rowcolor{lightgray} \hline \cC_n | \cD_n & z & X_n &  & \Omega_1^{(n)} & H_n(X^2,Y^2) =0& \includegraphics[width=1.85cm]{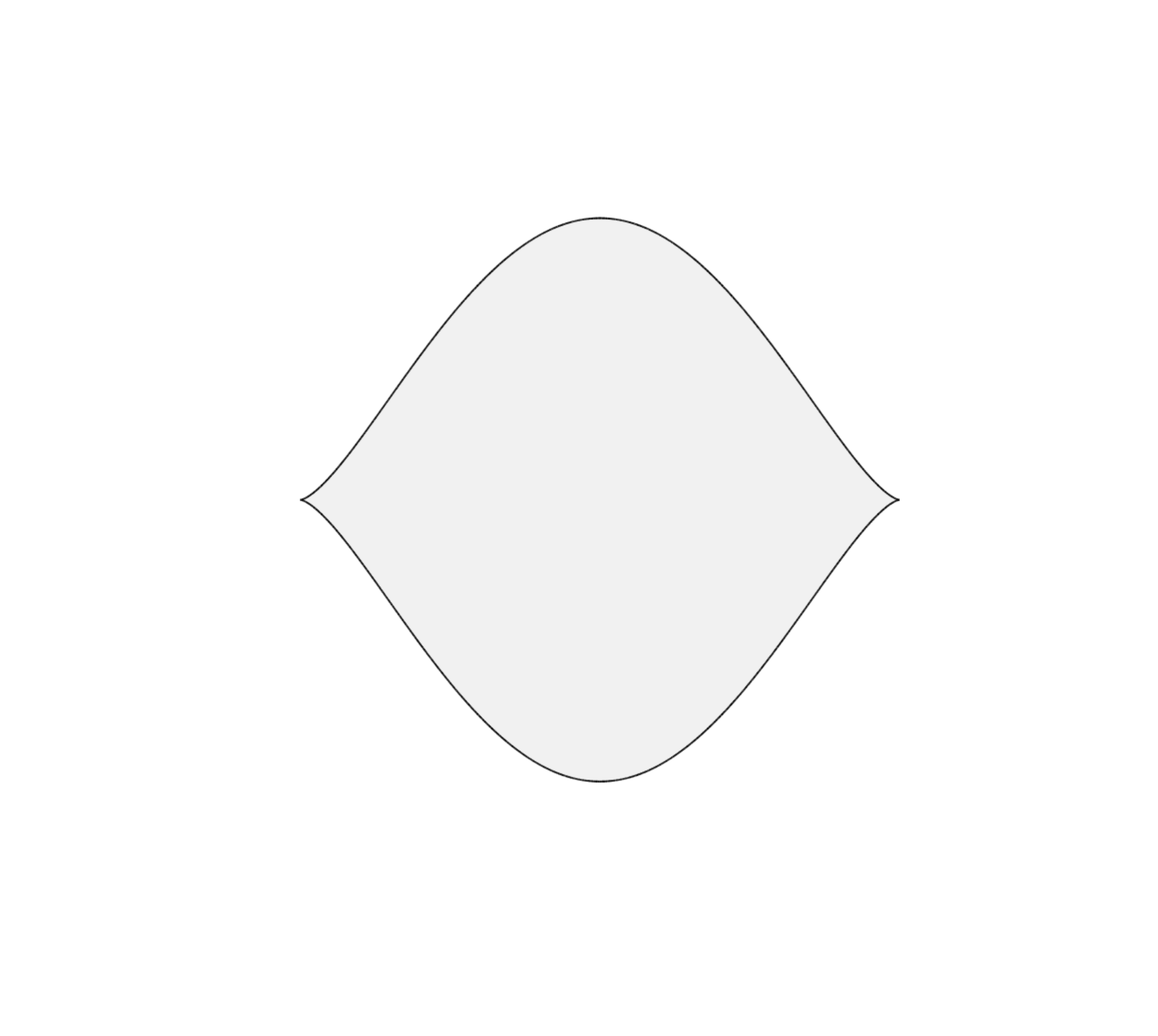}\includegraphics[width=2cm]{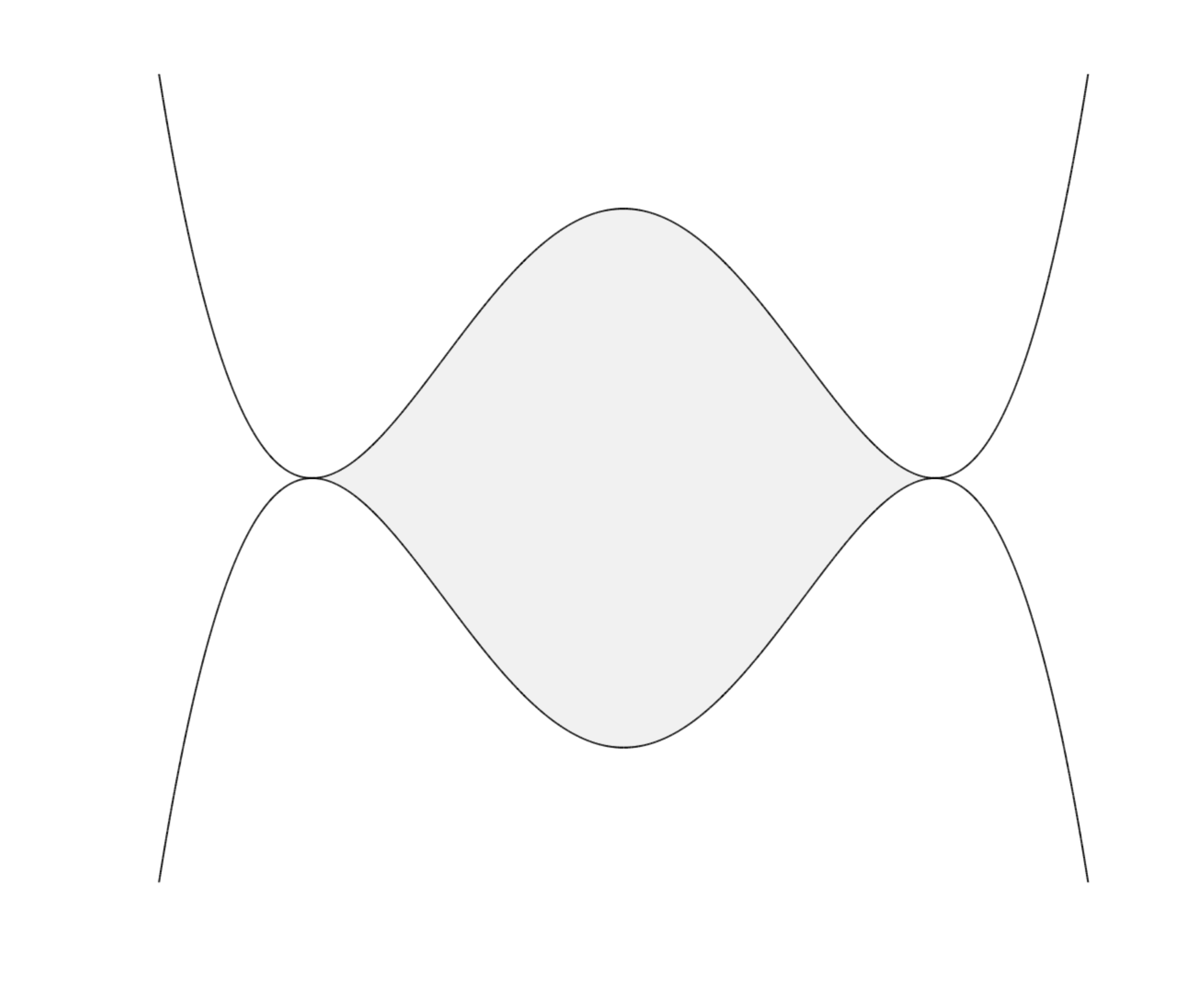}\\
\hline \cC_n & z& X_n & Y_n  & \Omega_2^{(n)} & H_n(X^2,Y^2) - \,{Z}^{2}=0& \includegraphics[width=2cm]{Omega2_3.pdf}   \includegraphics[width=2cm]{Omega2_5.pdf}\\
\rowcolor{lightgray} \hline \cD_{n,J} , \cD_n | D_{2n} & z^2 & X_n &  & \Omega_3^{(n)}&X H_n(X,Y^2) =0& \includegraphics[width=2cm]{Omega3_3.pdf}  \includegraphics[width=2cm]{Omega3_4.pdf}\\
\hline \cC_{nJ} , \cC_n | \cC_{2n} & z^2 & X_n & Y_n& \Omega_4^{(n)} &X(H_n(X,Y^2) -  Z^2) =0& \includegraphics[width=2cm]{Omega4_3.pdf}\\
\hline \cD_n & z^2 & X_n &z Y_n & \Omega_5^{(n)} &X H_n(X,Y^2) - Z^2 =0& \includegraphics[width=2cm]{Omega5_3.pdf}\\
\rowcolor{lightgray} \hline \cD_{2n,J} & z^2 & X_n^2 & &  D_6(\Omega_3^{(n)}) & XY H_n(X,Y) =0 &  \includegraphics[width=2cm]{Omega6_3.pdf}    \\
\hline 
\cD_n | \cD_{2n}, \cD_{n,J} &&&&&&\\
\cC_{n} | \cC_{2n}, \cC_{n,J}  & z^2 & X_n^2 & zY_n & \Omega_7^{(n)} & X H_n(X,Y) - Z^2 = 0 & \includegraphics[width=2cm]{Omega7_3.pdf}  \\
\hline \cC_{n} | \cC_{2n}, \cC_{n,J}& z^2 & X_n^2 & zX_n & \Omega_8^{(n)} &(XY - Z^2) H_n(X,Y) = 0 & \includegraphics[width=1.2cm]{Omega8.pdf}  \\
\hline \cC_{n} | \cC_{2n}, \cC_{n,J} & z^2 & X_n^2 & X_n Y_n & D_9(\Omega_4^{(n)}) &  X ( Y H_n(X,Y) - Z^2) = 0
&\includegraphics[width=2.5cm]{Omega9_3.pdf}
\\
\hline 
\end{array}
$$
When there are two groups, the first one is for $n$ odd, the second is for $n$ even. 

\pagebreak

For the tetrahedron / cube / octahedron family, we let  $H(X,Y) = 108 X^2 - 20 X - 2 Y^3 + 5 Y^2 - 4 Y + 36 X Y, $ and for the dodecahedron / icosahedron family, we set $S$ be as defined in Section \ref{sec:icosahedron}. 

$$
\begin{array}{|c||c|c|c|c|c|c|}
\hline \hbox{G} & \theta_1 & \theta_2 & \eta & \Omega & \hbox{Boundary}& \hbox{Picture} \\
\hline 
\rowcolor{lightgray}\hline \cT|\cO & O_3 & O_4 &  & \Omega_{11} & H(X^2,Y)=0& \includegraphics[width=2cm]{swallowtail.pdf}\\
\cT & O_3 & O_4 & O_6  & \Omega_{12} &  H(X^2,Y) - 4\,{Z}^{2}=0& \includegraphics[width=2cm]{Omega12.pdf}\\
\rowcolor{lightgray} \hline \cO_J & O_3^2 & O_4 &  & \Omega_{13} &X H(X,Y) =0 & \includegraphics[width=2cm]{cubcusptgt.pdf} \\
\hline \cT_J & O_3^2 & O_4 & O_6 & \Omega_{14} &X(H(X,Y) - 4 Z^2) =0& \includegraphics[width=2cm]{Omega14.pdf}\\
\hline \cO & O_3^2 & O_4 & O_3 O_6  & \Omega_{15} &Z^2-X H(X,Y) =0&  \includegraphics[width=2cm]{Omega15.pdf}\\
\rowcolor{lightgray} \hline \hline \cI_J & O_6 & O_{10} &   & \Omega_{21} & S(X,Y)=0&  \includegraphics[width=2cm]{icos.pdf}\\
\hline \cI & O_6 & O_{10} & O_{15}  & \Omega_{22} & Z^2 = S(X,Y)&  \includegraphics[width=2cm]{Omega22.pdf}\\
\hline 
\end{array}
$$

\paragraph{Covers} If $H$ is a subgroup of $G$, then polynomials that are invariant by $G$ are also invariant by $H$. It follows that we can express the $G$ invariants (primary and secondary) as polynomials in terms of the primary and secondary $H$-invariants (modulo $x^2+y^2+z^2-1$). These polynomials define a mapping from the $H$-domain onto the $G$-domain which is a $h$-covering, $h$ being the index of $H$ in $G$. 

Let us see what happens  for a few examples. This is specially easy for the $2$-covers.  We go from 3D to 2D models by simple projection (forgetting variable $z$). In other cases the effect is that of adding a new symmetry for instance : 
$$
\begin{array}{ccc}
\Omega_1^{(n)}    & \rightarrow & \Omega_3^{(n)} \\
(x,y)   & \mapsto & (x^2, y)
\end{array}
$$
$$
\begin{array}{ccc}
\Omega_2^{(n)}    & \rightarrow & \Omega_5^{(n)} \\
(x,y,z)   & \mapsto & (x^2, y,xz)
\end{array}
$$
$$
\begin{array}{ccc}
\Omega_2^{(n)}    & \rightarrow & \Omega_4^{(n)} \\
(x,y,z)   & \mapsto & (x^2, y,z)
\end{array}
$$

The case of $\cC_3|\cD_3$ as subgroup of $\cT|\cO$ is a little bit more tedious because we did not choose the same coordinates for the representations. We can see what happens taking for  $\cC_3|\cD_3$ the variables $(x+y+z)/{\sqrt{3}}$ and $xyz$ to have a common invariant.  It is a linear computation to get  $a$, $b$, $c$ and $d$ such that 
$$x^4+y^4+z^4 = a (x+y+z)^4 + b (x+y+z)^2(x^2+y^2+z^2) + c (x+y+z) xyz + d (x^2+y^2+z^2)^2. $$
We obtain a map of the form : 
$$
\begin{array}{ccc}
\Omega_1^{(3)}    & \rightarrow & \Omega_{11}\\
(x,y)   & \mapsto & (y, a x^4 + b x^2 + c xy + d). 
\end{array}
$$

\section{ The bounded two dimensional models of~\cite{BOZ2013} \label{sec.BOZmodels}}

We provide here for completeness the complete list of models in dimension 2 described in~\cite{BOZ2013}.  With the restriction that the valuation is the usual one, they are the only ones which may occur up to affine transformations.  We indicate the (scalar) curvature when it is  constant ($+$ when it is  a positive constant, $0$ otherwise). When no curvature is indicated, it comes from the fact that the metric is not unique  (models~2 and~3), in which case there exist at least one metric for which the curvature is constant and positive), or it is not constant (model 7).  Up to isomorphism, one may replace one parabola by an horizontal line in model 4, so that this changes the degree in the boundary (in this particular case however, the co-metric is no longer unique).

 Up to isomorphism, we have
$$\bcas (2) \simeq \Omega_1^{(1)}, ~(3)\simeq \Omega_3^{(2)},  ~(4)\simeq \Omega_1^{(2)}\simeq \Omega_3^{(1)},\\ (5) \simeq \Omega_6^{(2)}\simeq \Omega_4^{(4)}, ~
(8) \simeq \Omega_2^{(3)},~ (9) \simeq \Omega_13, ~(10) \simeq \Omega_{11}\ecas
$$
 $$
\begin{array}{|c||c|c|c|c|c|c|}
\hline \sharp  &\hbox{Curv.}& \hbox{d($\Omega$)}  & \hbox{Boundary}& \hbox{Picture} \\
\hline 1 & 0 & 4 &   (1-X^2) (1-Y^2) =0& \includegraphics[width=1.5cm]{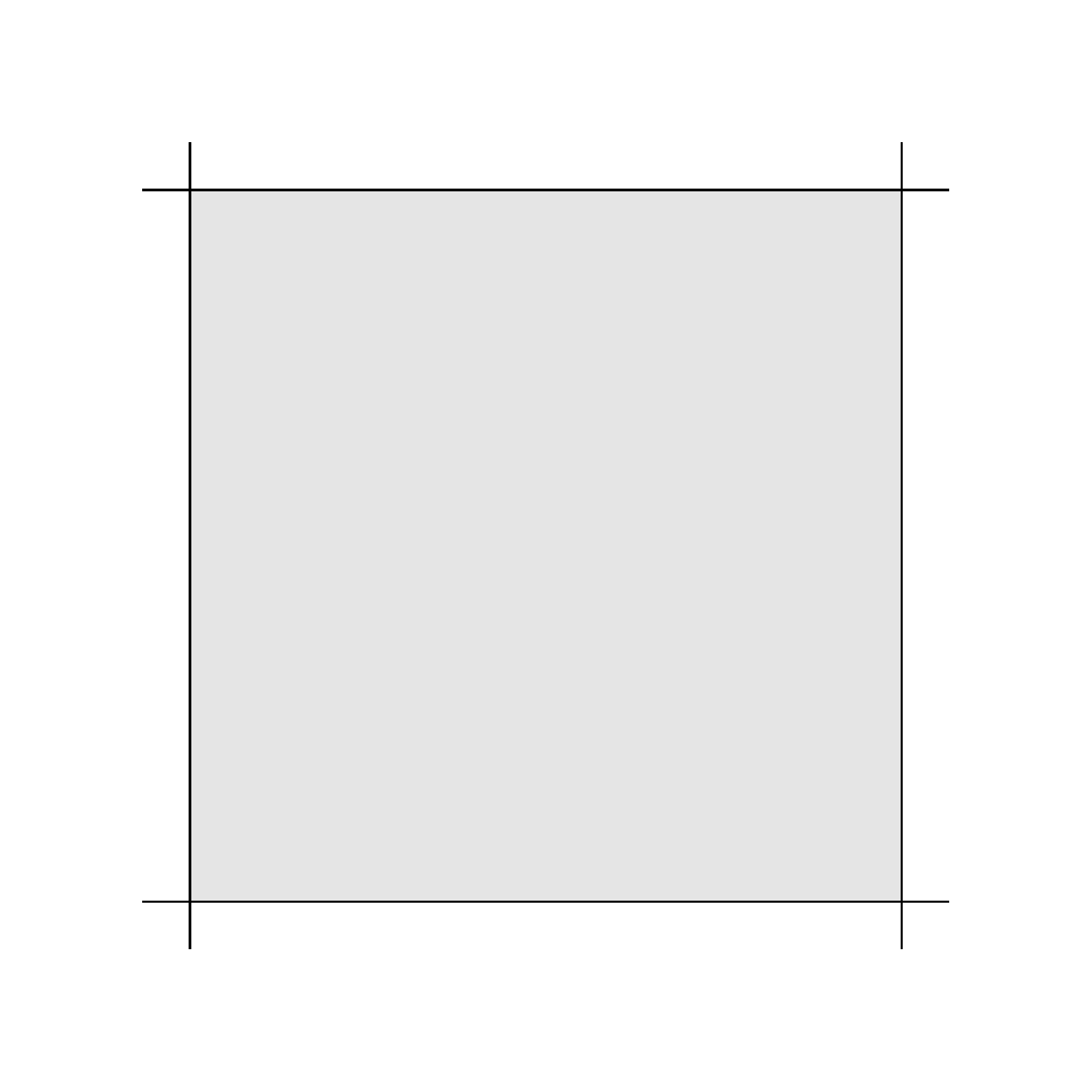}\\
\hline 2 & \simeq & 2 &   1 - X^2 - Y^2=0& \includegraphics[width=1.5cm]{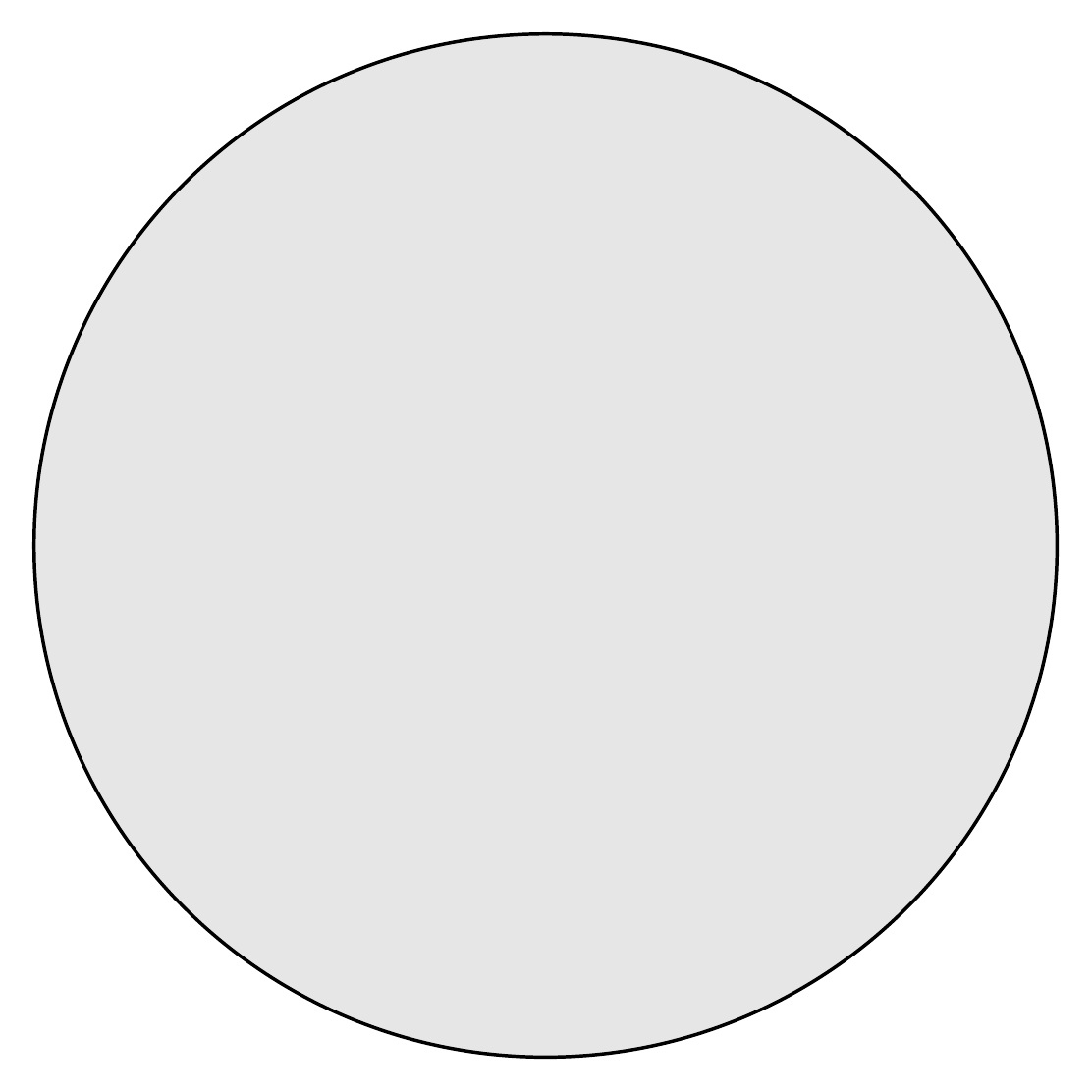}\\
\hline 3 & \simeq   & 3  & X Y (1 - X - Y) =0&\includegraphics[width=1.5cm]{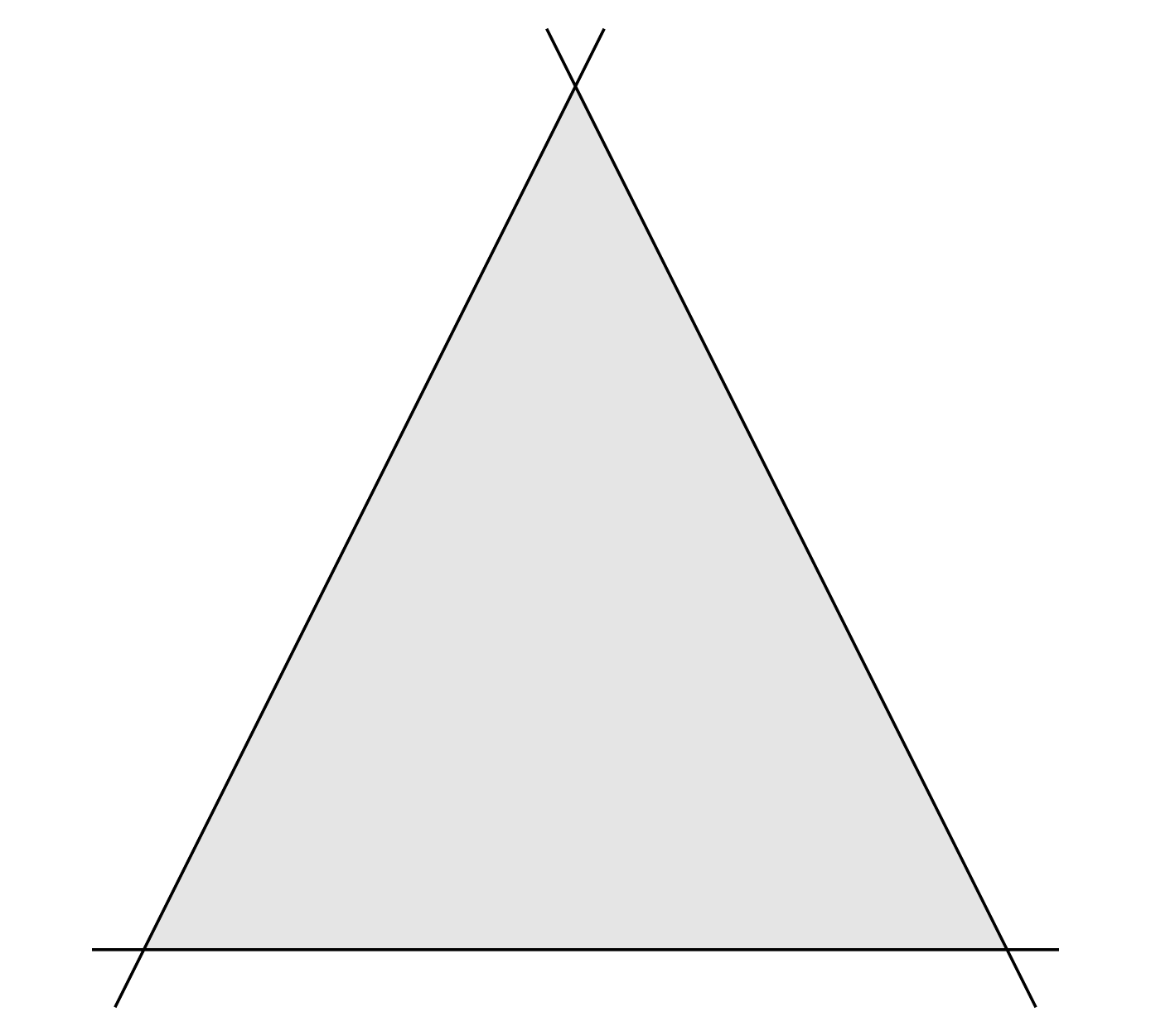}\\
\hline 4& + & 4,3  &  (1 - X^2)^2 - Y^2 = 0& \includegraphics[width=1.5cm]{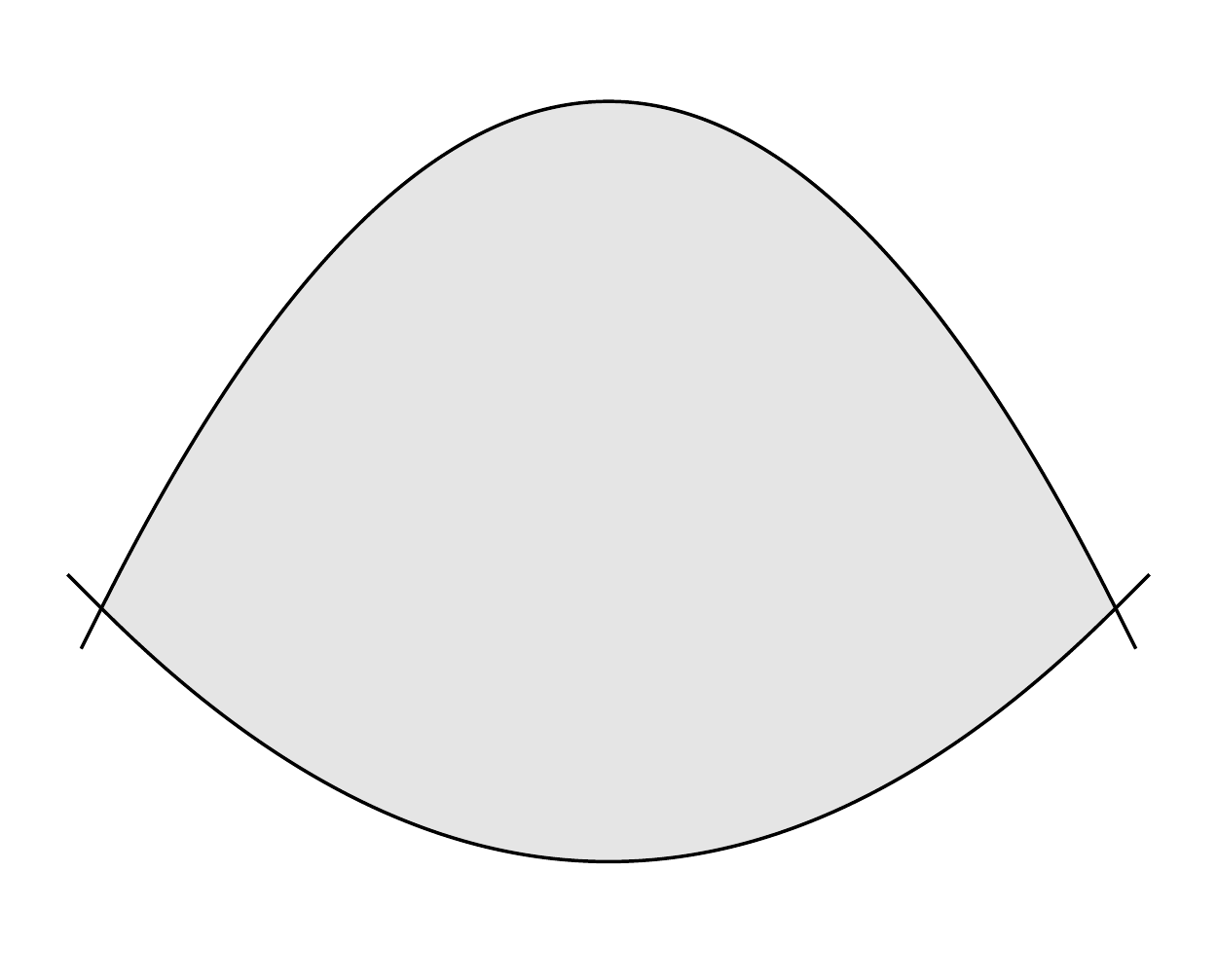}  \\
\hline5& +  & 4 & Y (1-X) (X^2-Y)&\includegraphics[width=1.3cm]{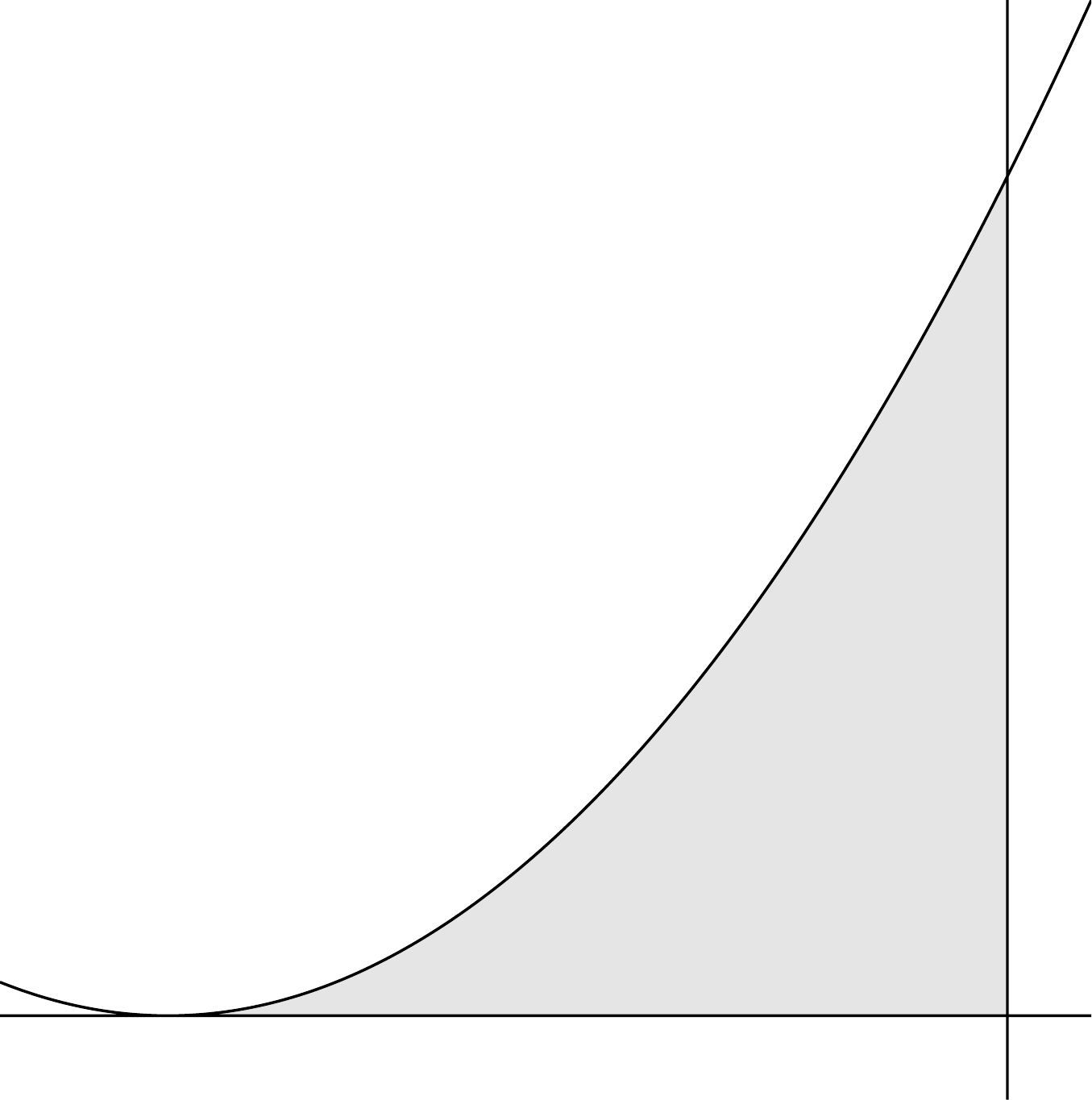} \\
\hline 6 & 0  & 4  &(Y-X^2) ((Y+1)^2 - 4X^2) = 0&\includegraphics[width=1.5cm]{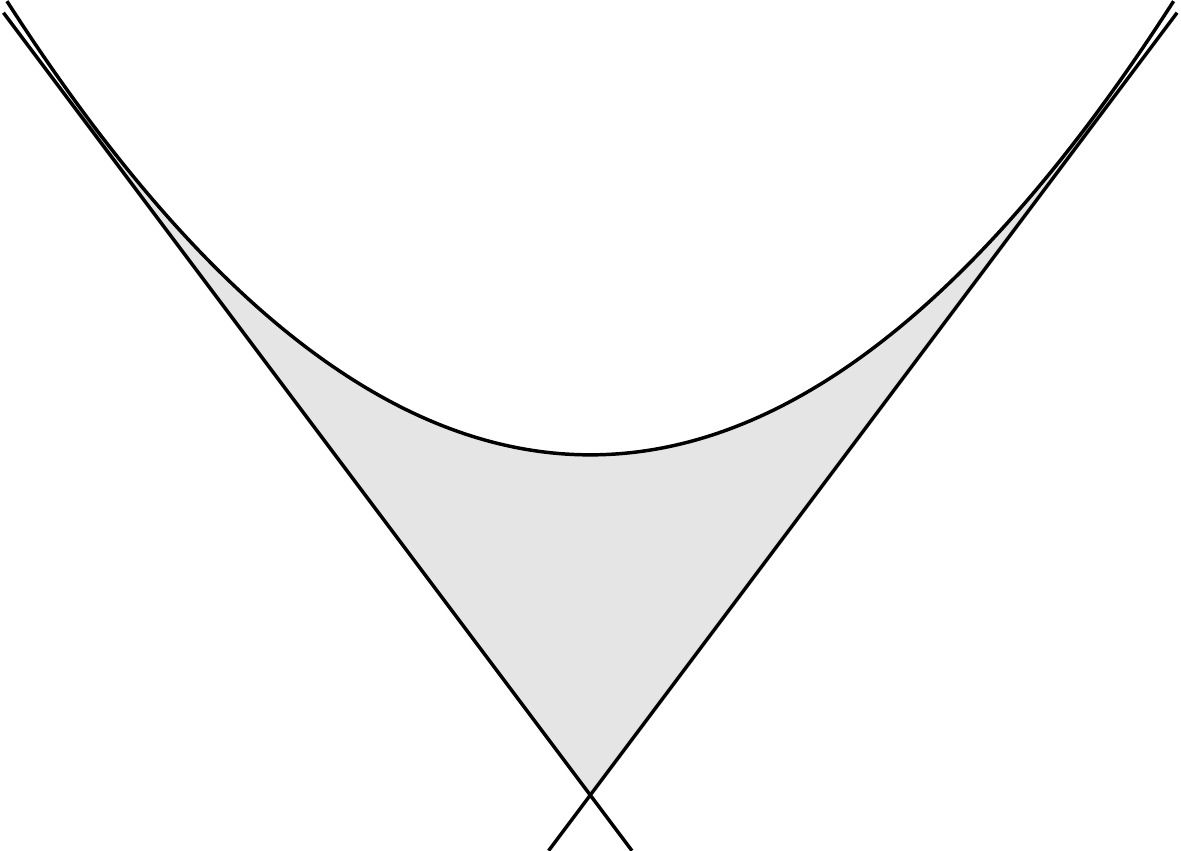}\\
\hline7& \simeq  & 3 &Y^2 - X^2(1-X) =0&\includegraphics[width=1.5cm]{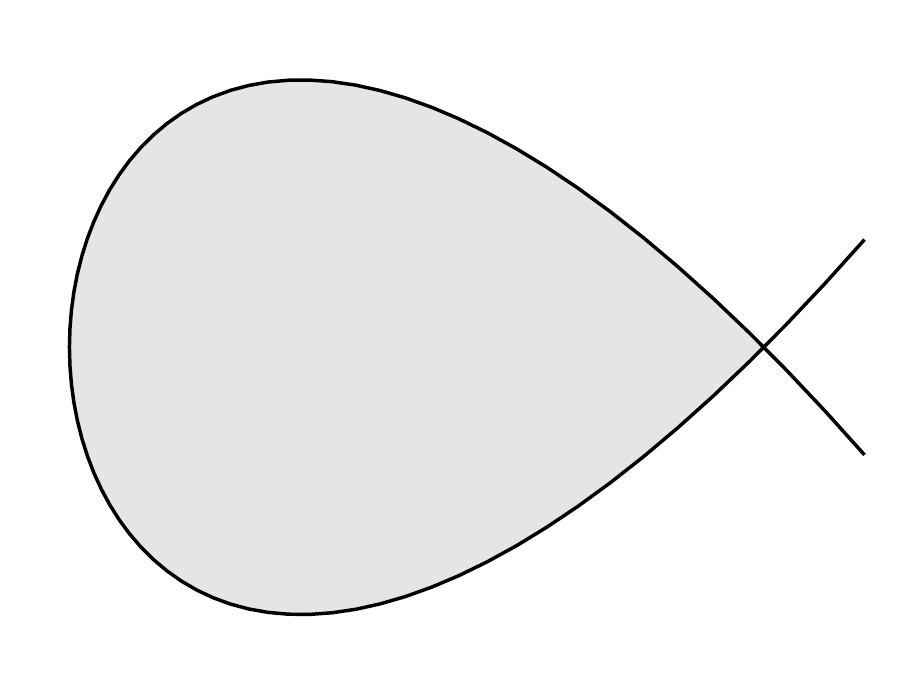}\\
\hline 8 & + & 4 & (Y^2 - X^3) (X-1) =0& \includegraphics[width=1.5cm]{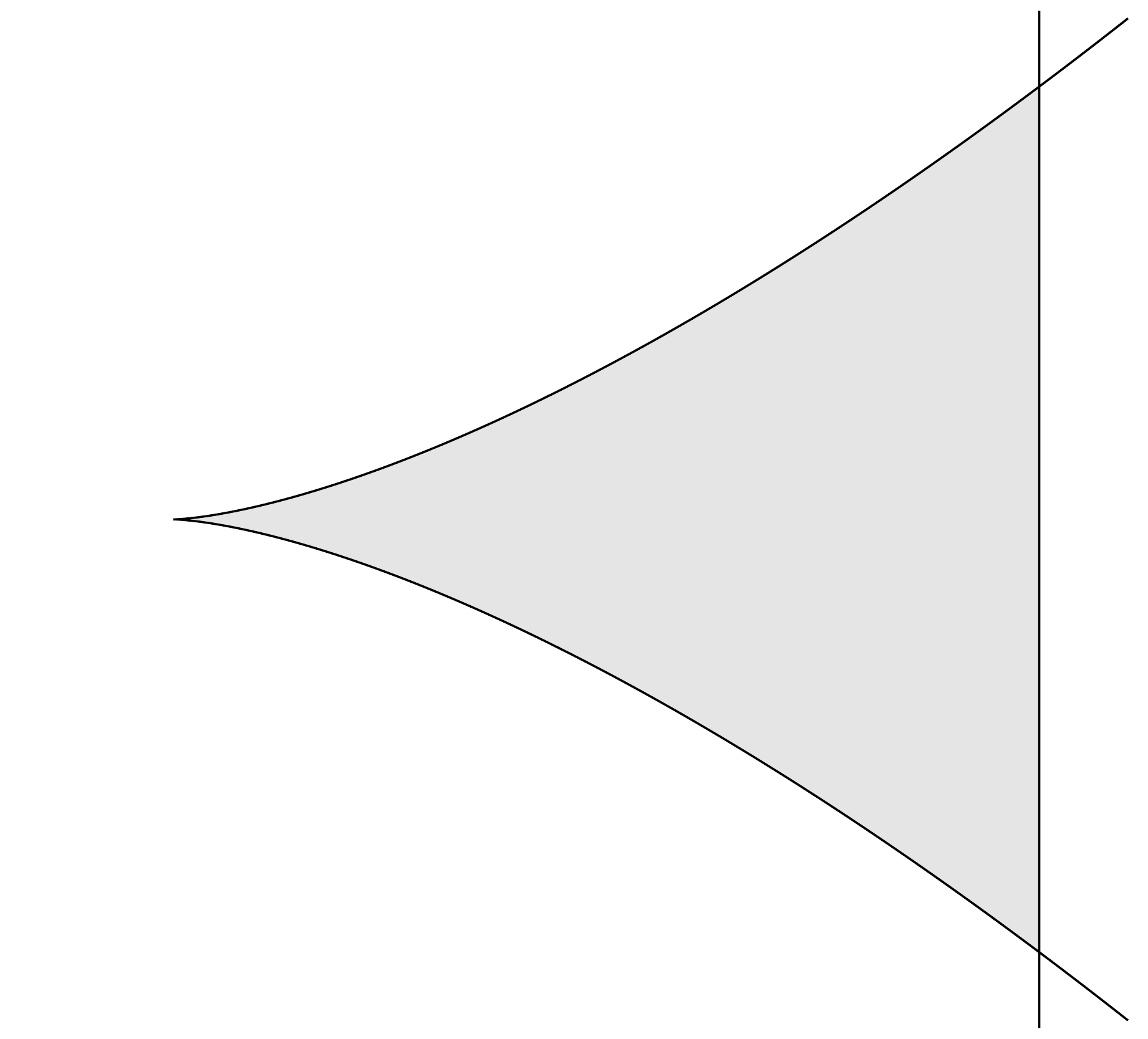}\\
\hline 9 & +  & 4&(Y^2  - X^3) (2(Y -1)-3(X -1)) =0& \includegraphics[width=1.5cm]{cubcusptgt.pdf}\\
\hline 10& +  & 4& 4 X^2 - 27 X^4 +16 Y -128Y^2 - 144 X^2 Y +256Y^3 =0 &\includegraphics[width=1.5cm]{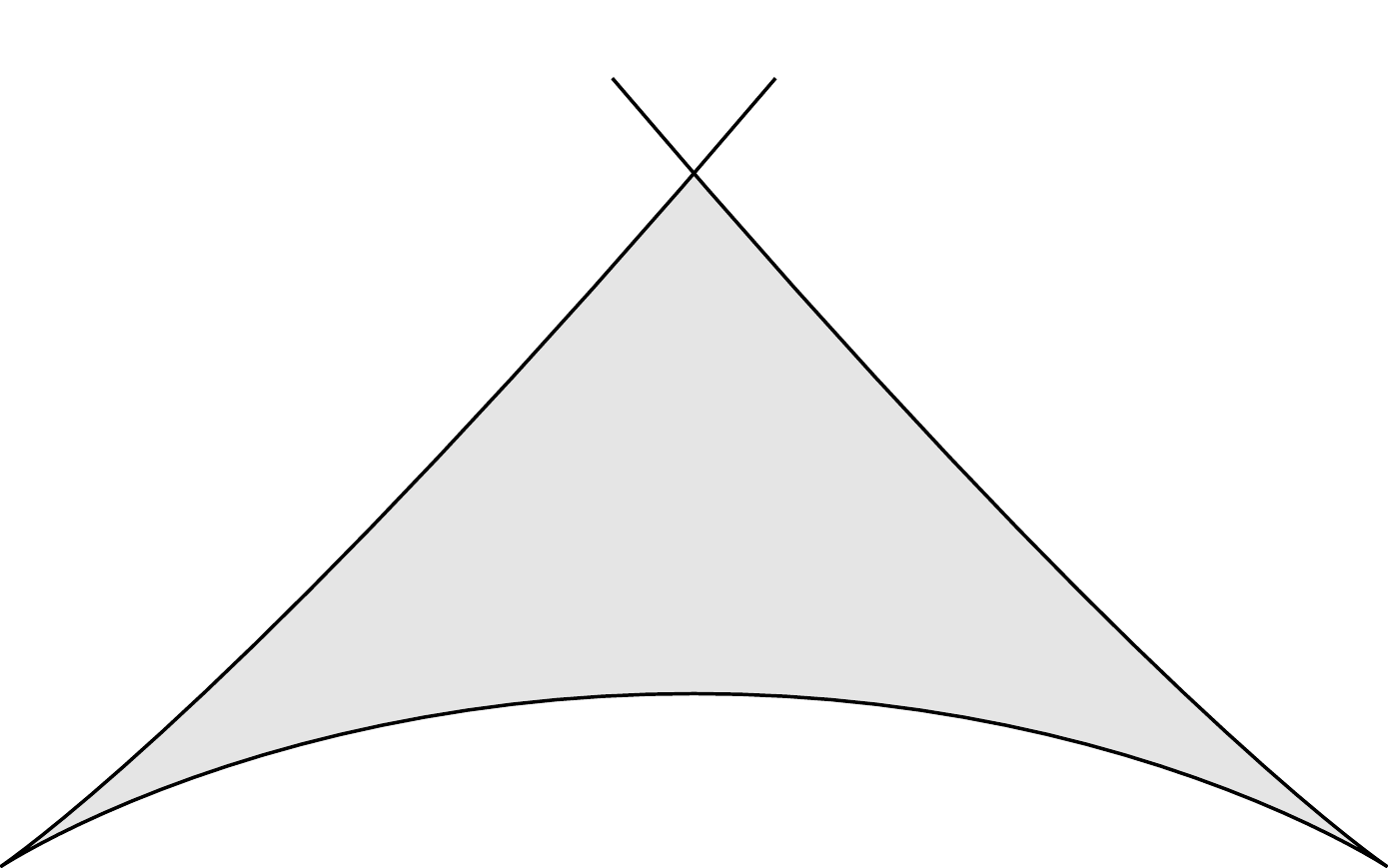}\\
\hline 11& 0  &4  &   (X^2 +Y^2)^2 +18(X^2 +Y^2)-8X^3 +24XY^2 -27=0   &\includegraphics[width=1.5cm]{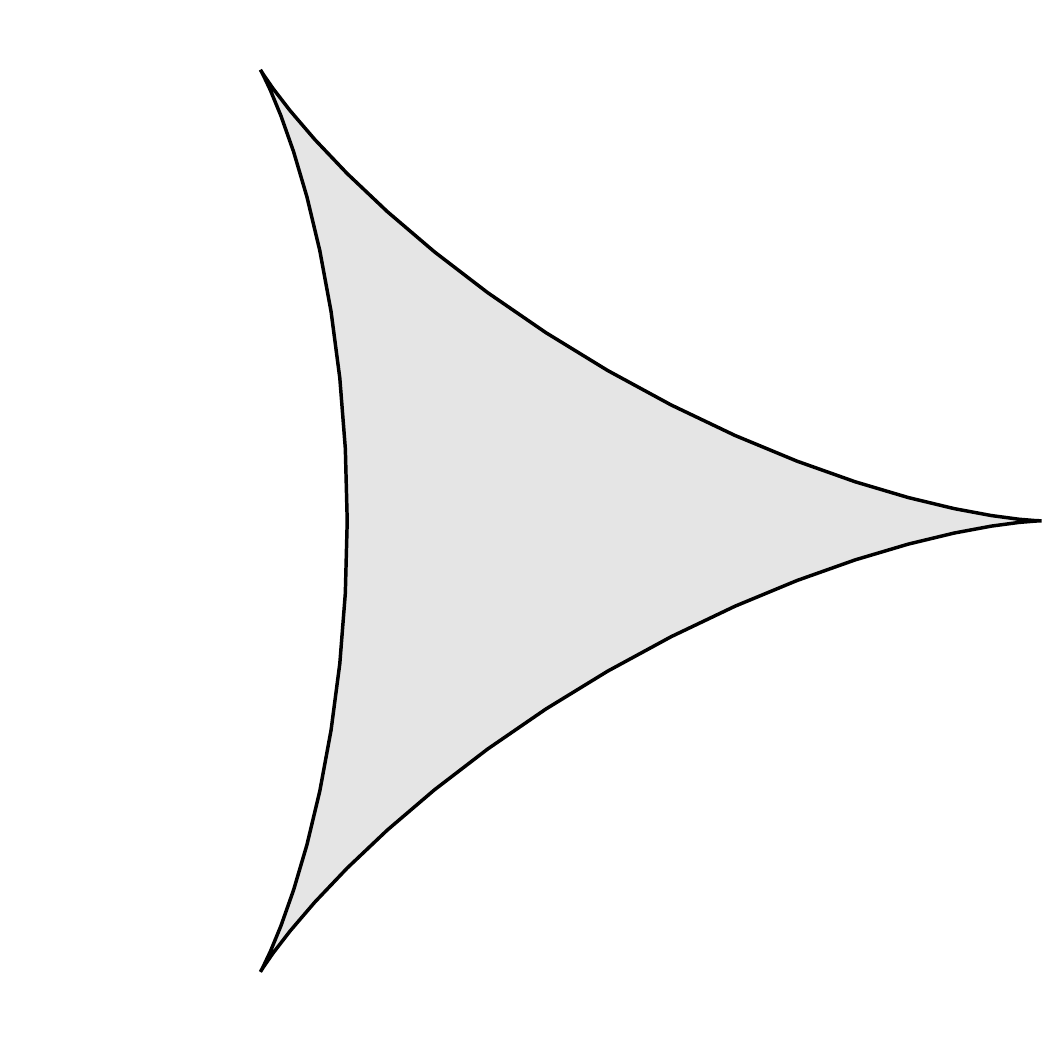} \\
\hline 
\end{array}
$$

\section{Further remarks\label{sec.further}}

In the various models presented here, it happens that the primary invariants provide a closed system. The reason is that all these groups are subgroups of finite Coxeter groups, for which the invariants are our primary invariants.  It is not true that this is always the case. Here is an example provided by Y. Cornulier of a group in dimension 4 which is not a subgroup of a finite Coxeter group. Let $M_p$ the matrix of a rotation with angle $2\pi/p$ in $\bbR^2$, where $p$ is an odd prime. Let $\cI_2$ be the $2\times 2$ identity matrix. Then, we consider the group generated by 
$$N_1= \bpm M_p&0\\0& \cI_2\epm, N_2= \bpm \cI_2 & 0\\ 0& M_p\epm, J= \bpm 0 & \cI_2\\ \cI_2& 0 \epm.$$
This group has $2p^2$ elements, $N_1^{k_1}N_2^{k_2}, JN_1^{k_1}N_2^{k_2}$, $0\leq k_1,k_2\leq p-1$, and, thanks to Mollien's formula~\eqref{eq.mollien}, the Hilbert sum is easily computed
$$F(t)= \frac{1+t^p+2t^{p+2}+ 2t^{2p}+ t^{2p+2}+ t^{3p+2}}{(1-t^2)(1-t^4)(1-t^{2p})(1-t^p)}.$$

This leads to the description of primary and secondary invariants when restricted on the unit sphere in $\bbR^4$. Following Section~\ref{sec.cyclic}, we identify $\bbR^4\simeq \bbC^2$, and for a pair $(z_1,z_2)$, consider 
$z_j^p = X_j+ i Y_j$ and $R_j = |z_j|^2$,  $j= 1,2$. Then we chose 
$$\theta_1= X_1+X_2, \theta_2= X_1X_2, \theta_3= R_1R_2$$ as primary invariants, and the secondary invariants may be chosen as
\beqnas &&\eta_1=Y_1+Y_2, 
\eta_2=(Y_1-Y_2)(R_1-R_2),   \eta_3=(X_1-X_2)(R_1-R_2),\\
&& 
\eta_4= Y_1Y_2,  \eta_5=(X_1-X_2)(Y_1-Y_2), 
 \eta_6=(X_1Y_1-X_2Y_2)(R_1-R_2), 
\eta_7=  \eta_3\eta_4
\eeqnas
It turns out that $\theta_1, \theta_2, \theta_3$ is not closed for $\Gamma$, where $\Gamma$ is the square field operator on the unit sphere in $\bbR^4$. For example 
 $$\Gamma(\theta_1,\theta_2)= p^2(T_{p-1}-2\theta_1\theta_2),$$
where $T_k= R_1^kX_2+ R_2^kX_1$, and $T_{p-1}$ may be expressed as 
$$T_{p-1} = \theta_1 Q_1(\theta_3)+ T_1Q_2(\theta_3),$$ where $Q_i$ are polynomials, and $T_1= (\theta_1+ \eta_3)/2$, so that $\Gamma(\theta_1, \theta_2)$ is not a polynomial of $(\theta_1, \theta_2, \theta_3)$.

In this example, one may observe that indeed $(\theta_1, \theta_2, \theta_3, \eta_3)$ form a closed system for $\Gamma$, but this does not provide any model in $\bbR^4$ (the boundary equation is not satisfied).

A final remark is that our various 3-d models constructed from 2-d one could appear as provided by Coxeter groups in dimension 4. If such would be the case, the natural Ricci curvature carried by the associated cometric $\Gamma$ would be constant (since it would locally be the spherical co-metric seen through a diffeomorphism). One may easily check that this is not the case.

\bibliographystyle{amsplain}   
\bibliography{Bib.Models2D}
\end{document}